\documentclass[11pt]{amsart}
\usepackage{amsmath,amssymb,amsthm,amscd}

\def\p{\partial}
\def\R{\mathbb{R}}
\def\C{\mathbb{C}}
\def\N{\mathbb{N}}

\def\l{\lambda}

\def\cB{{\mathcal B}}

\def\cE{{\mathcal E}}

\def\cG{{\mathcal G}}

\def\cH{{\mathcal H}}

\def\cJ{{\mathcal J}}

\def\cR{{\mathcal R}}
\def\cS{{\mathcal S}}
\def\cT{{\mathcal T}}

\def\cV{{\mathcal V}}
\def\cW{{\mathcal W}}

\def\tr{\text{Tr}}

\newtheorem{thm}{Theorem}
\newtheorem{prop}{Proposition}[section]
\newtheorem{lemma}[prop]{Lemma}
\newtheorem{cor}[prop]{Corollary}
\newtheorem{rmk}[prop]{Remark}
\newtheorem{exam}[prop]{Example}
\newtheorem{defn}[prop]{Definition}

\newtheorem{pl}[prop]{Problem}

\numberwithin{prop}{section}

\begin{document}
\title{Energy minimizing harmonic almost complex structures}
\author{Weiyong He}

\email{whe@uoregon.edu}
\address{University of Oregon, Eugene, OR, 97403}

\begin{abstract}
Consider all compatible almost complex structure with a fixed almost Hermitian metric $(M, g)$ and we study energy-minimizing almost complex structures, with the energy functional $E(J)=\int_M |\nabla J|^2dv$. C. Wood  studied this problem in 1990s and he named a critical point  \emph{harmonic almost complex structure}, satisfying the equation $[J, \Delta J]=0$. Since then there are considerate interest to study these objects but in general harmonic almost complex structures are not well understood. 

We introduce the notion of \emph{admissible} almost complex structure as a natural generalization of smooth almost complex structure to $W^{1, 2}$ setting. The main goal of the paper is to study the regularity of \emph{admissible} energy-minimizing almost complex structures on a compact Riemannian (almost Hermitian) manifold $(M, g)$ of even dimension $m=2n$. Our starting point is to observe that a critical point satisfies the equation $\Delta J-J\nabla J\nabla J=0$. This is a semi-linear elliptic system of tensor-valued functions which resembles the elliptic system of harmonic maps (vector-valued functions) in many ways. We study the regularity of energy-minimizing almost complex structures and prove that, all major results for energy-minimizing harmonic maps hold in our setting, notably the seminal work of Schoen-Uhlenbeck and recent improvement of regularity by Cheeger-Naber. We use comparison almost complex structures as was done by Schoen-Uhlenbeck. As in Schoen-Uhlenbeck's paper, finding comparison almost complex structures which satisfies the constraints is a major technical difficulty.
The major difference is that due to the compatibility with the background metric, one cannot apply the method of Schoen-Uhlenbeck directly. Since the background metric gives  strong restrictions on compatible almost complex structures and hence it has to play a very important role in the regularity theory.  This is the main point that differs from the theory of  harmonic maps.  Our construction of comparison almost complex structures relies on a new input  that given an arbitrary almost complex structure, we can construct a unique (canonical) almost complex structure  which is compatible with a given metric. 

Smooth maps between manifolds are not necessarily dense in $W^{1, 2}$ (Sobolev maps), which is a well-known phenomenon in the theory of harmonic maps. Similarly an \emph{admissible} almost complex structures is not necessarily an almost  complex structure.  For example manifolds such as $S^4$ do not admit any smooth almost complex structure, but $W^{1, 2}$ admissible almost complex structures exist. In particular, \emph{admissible} energy-minimizing almost complex structures exist on $S^4$ (for any Riemannian metric) and form a compact set. These objects rely partly on the smooth structure of the underlying manifolds and it might be interesting objects to study further.  
 \end{abstract}

\maketitle
\tableofcontents

\section{Introduction}
This is the first paper in a series in which the author will study almost Hermitian geometry, from a point view of geometric analysis.  
Almost complex manifolds contain central objects in the modern theory of differential geometry, such as complex manifolds, symplectic manifolds and K\"ahler manifolds. K\"ahler geometry, complex geometry and symplectic geometry have gained prominence for long time, while the study of general almost complex manifolds (almost Hermitian geometry) remains largely unexplored, even though there is a long history to study these objects.  There are attempts recently to study almost Hermitian structure using the modern theory of geometric analysis, in particular by the means of curvature flows. 

We will adopt a different approach in this paper.  One natural question in almost Hermitian geometry is the following, given an almost Hermitian manifold $(M, g)$ of dimension $m=2n$ with compatible almost complex structures, what is the ``best" almost complex structure? 
This problem dates back to Calabi-Gluck \cite{CG} and C. Wood \cite{Wood} in 1990s using the theory of twistor bundles. In particular, C. Wood \cite{Wood1} came up with the notion of the \emph{harmonic almost complex/Hermitian structure} by considering minimizing the energy functional, for all compatible almost complex structures, 
\begin{equation}
E(J)=\int_M |\nabla J|^2 dv.\end{equation}
Clearly a K\"ahler structure $\nabla J=0$ gives an absolute minimizer of the energy functional. But there are  various very interesting examples which are  absolute minimizers of the energy functional but are not K\"ahler, see for example \cite{BorSalvai}. In this sense one can view the (energy-minimizing) harmonic almost Hermitian structures as a natural generalization of K\"ahler structures. 
The harmonic almost complex structure has gained considerate interest and we refer the readers to the recent survey paper \cite{Davidov} for the background, history and results in this subject. A main interest was to construct and to verify examples of energy-minimizers (see \cite{BorSalvai}) and to study some properties of harmonic almost complex structure, for example, whether the harmonic almost complex structure, viewed as a map from the manifold to its twistor space, is a genuine harmonic map.
 
We study the harmonic almost complex (Hermitian) structure from the point of view in geometric analysis in this paper. The Euler-Lagrangian equation of the energy functional is given by $[J, \Delta J]=0$, which was first written down by C. Wood \cite{Wood, Wood1}. Our starting point is to rewrite this equation in the following form,
\begin{equation}\label{H1}
\Delta J-J\nabla_p J \nabla_p J=0,
\end{equation}
where $\Delta$ is the rough Laplacian defined by the metric $g$. An advantage is that \eqref{H1} is of the form of a semilinear (strictly) elliptic system, which has very similar structure as the harmonic map equation. In some special cases, harmonic almost complex structures are harmonic maps from the manifold to its twistor space (see \cite{Davidov} for example); hence in these cases, the deep results in harmonic map can be directly quoted. However, the theory of harmonic almost complex structures has substantial differences with its own interest. Hence it should have an independent theory. For example, constant harmonic maps are of no interest at all with zero energy. But when there exists a harmonic almost complex structure with zero energy on a given almost Hermitian manifold  ($\nabla J=0$, that is K\"ahler structure) is a very delicate and important problem: the K\"ahler condition gives very strong restrictions on the manifold and the metric.   \\

We introduce \emph{admissible} almost complex structures, which are sections of $\Gamma (T^*M\otimes TM)$ with $W^{1, 2}$ coefficients, and satisfy the constraint conditions $J^2=-id, g(J\cdot, J\cdot)=g$ almost everywhere. We denote the space of all smooth compatible almost complex structures  by $\cJ_g$ and the space of admissible almost complex structures by $W^{1, 2}(\cJ_g)$. The later is similar as  the $W^{1, 2}$ Sobolev maps studied in  the theory of harmonic maps. The existence of an energy-minimizing \emph{weakly harmonic} almost complex structures in $W^{1, 2}(\cJ_g)$ is straightforward by a standard argument in the theory of calculus of variations. The main results in the paper are to demonstrate that similar regularity results as in the theory of energy-minimizing harmonic maps hold for energy-minimizing \emph{weakly harmonic} almost complex structures.  The arguments follow the general strategy in the regularity theory of energy minimizing harmonic maps by Schoen-Uhlenbeck \cite{SU}.  We should emphasize instantly that there are substantial differences, due to the constraint condition $g(J\cdot, J\cdot)=g$. Namely the background metric $g$ plays an essential role.  In particular the method of construction of comparison maps in \cite{SU} does not work directly in our case. As a comparison, it is well-known that the background metric $(M, g)$ plays few very rule in the theory of harmonic maps $f: M\rightarrow N$ (rather the target manifold $(N, h)$ is essential). Our construction of comparison almost complex structures uses a very basic observation in linear algebra, that there is a unique way to project an almost complex structure to a new one, which is compatible with the metric, see Proposition \ref{matrix}. Combining this elementary fact and deep analysis in energy-minimizing harmonic maps, we can prove the first main results in the paper,
\begin{thm}[Section 4]\label{mainthm1}Given an almost Hermitian manifold $(M, g)$ with compatible almost complex structures in $\cJ_g$, there exists an energy minimizing weakly harmonic almost complex structure $J\in W^{1, 2}(\cJ_g)$, which satisfies the equation $\Delta J-J\nabla_p J\nabla_p J=0$ in the weak sense. Moreover,  there exists $\epsilon=\epsilon(n, M, g)$ such that if, for some $r\in (0, 1)$ we have
\[
r^{2-m}\int_{B_p(r)} |\nabla J|^2 dv<\epsilon,
\]
then $J$ is smooth in $B_{r/4}$. The singular set $\cS$ has Hausdorff dimension at most $m-3$. If $p\in \cS$, then any tangent almost complex structure, is an energy-minimizing weakly harmonic almost complex structure on $\R^{m}=\R^{2n}$ compatible with Euclidean metric and it is of homogeneous degree zero. In particular, there is a natural stratification of the singular set 
\[
\cS^0\subset \cS^1\subset\cdots\subset \cS^{2n-3}=\cS,
\]
such that the Hausdorff dimension $\text{dim}\,\cS^k\leq k, k=1, \cdots, 2n-3$, 
where $\cS^k$ denotes the singular set where any tangent almost complex structure is at most $k$-homogeneous. 
\end{thm}
We should mention that  a tangent almost complex structure can be identified with the energy-minimizing harmonic maps in $W^{1, 2}(\R^{2n}, SO(2n)/U(n))$ of homogeneous degree zero. In particular, when $n=2$, these are  energy-minimizing harmonic maps in $W^{1, 2}(\R^{4}, S^2)$ of homogeneous degree zero. Indeed almost all results in energy-minimizing harmonic maps have their companion in our setting, including an improved version of $\epsilon$-regularity result in \cite{SU}, and compactness results for energy-minimizing harmonic maps by Luckhaus \cite{Luckhaus1, Luckhaus2} etc. 
In particular, by applying the method in the recent work of Cheeger-Naber \cite{CN, CN1}, we have the following \emph{quantitative stratification}, 

\begin{thm}[Section 5]\label{mainthm2}Let $J$ be an energy-minimizing almost complex structure in $W^{1, 2}(\cJ_g)$ on $(M, g)$. Then for $\eta>0, r\in (0, 1)$, the $k$-th quantitative stratum $\cS^k_{\eta, r}$ has the following volume estimate,
\[
\text{Vol}(\cT_r(\cS^k_{\eta, r}))\leq C(M, g, m, \eta) r^{m-k-\eta}.
\] 
Let $r^J: M\rightarrow [0, \infty)$ be the regularity scale defined by
\[
r^J(p)=\max\{r\in [0, 1]: \sup_{B_r(p)} r|\nabla f|+r^2|\nabla^2 f|\leq 1\}.
\]
Then for any $p\in [2, 3)$, there exists a uniform constant $C_0=C_0(M, g, m, p)$ such that  
\[
\int_M |\nabla J|^pdv, \int_{M} |\nabla^2J|^{p/2}dv\leq \int_M (r^J)^{-p}dv\leq C_0.
\]
\end{thm} 

Harmonic almost complex structures should be viewed as a tensor-valued companion of harmonic maps. The theory of harmonic maps is vast and it has tremendous applications. There are thousands of articles and many nice books, survey papers and we mention a few \cite{EL1, EL2, HW, LW, SY, Simon} and the references therein for more information.

We organize the paper as follows. We introduce the notion of \emph{admissible almost complex structure} and prove the existence of energy-minimizer in Section 2. The main regularity results are proved in Section 3, 4, and 5.  We establish the monotonicity formula for energy-minimizing almost complex structures in Section 3, where we have constructed comparison almost complex structures satisfying the constraints. We prove the Schoen-Uhlenbeck type regularity for energy-minimizing  almost complex structures in Section 4 and establish  the quantitative stratification and the improved regularity results in Section 5, applying the method of Cheeger-Naber. In Section 6 we discuss various examples of energy-minimizing almost complex structures and discuss energy-minimizing problem on other Sobolev spaces of almost complex structures. In Appendix, we include some general discussions of harmonic almost complex structure, and derive that a continuous weakly harmonic almost complex structure is smooth.\\

{\bf Acknowledgment:} the author is supported partly by NSF, award no. 1611797.
 
\numberwithin{equation}{section}
\numberwithin{thm}{section}

\section{Admissible energy-minimizing almost complex structure}
Let $(M, g)$ be a compact Riemannian manifold with a compatible almost complex structure, of real dimension $2n=m$.  We denote $\cJ_g$ to be the space of smooth almost complex structures which are compatible with $g$.
The compatible conditions give the constraints, for any $p\in M$, 
\begin{equation}\label{C2}
J^2_p=-id, g(J_px, J_py)=g(x, y), \forall x, y\in T_pM
\end{equation}
Locally  the constraints \eqref{C2} imply the poinwise bound $|J|^2=g^{ij}g_{kl}J^k_iJ^l_j=2n$.

Let $S: TM\rightarrow TM$ be an endomorphism.  Its $W^{1, 2}$ norm is given by
\[
\|S\|_{W^{1, 2}}=\left(\int_M \left(|S|^2+|\nabla S|^2\right)dv\right)^{1/2}.
\]
We define $W^{1, 2}(M, \text{End}(TM))$ to be the completion of smooth endomorphisms  with respect to $W^{1, 2}$-norm.  Note that $\text{End}(TM)$ can be identified with $T^*M\otimes TM$, and hence $W^{1, 2}(M, \text{End}(TM))$ are simply the sections of $T^*M\otimes TM$ (or rather, the completion of spaces of sections of $T^*M\otimes TM$ with respect to $W^{1, 2}$-norm.)
Locally, if we choose a local coordinate, these are matrix-valued functions with $W^{1, 2}$ entries.  We can also define other $L^p$ and $W^{k, p}$ spaces of endomorphisms, or spaces of sections of tensor bundles. Note that these spaces all have linear structures, hence automatically have a local smooth structure. 

\begin{defn}\label{admissible}We say $J$ an \emph{admissible} or a \emph{weak} almost complex structure which is compatible with $g$ if $J\in W^{1, 2}(M, \text{End}(TM))$ and  moreover $J_p: T_pM\rightarrow T_pM$ is a well-defined endomorphism almost everywhere such that \eqref{C2} holds a.e. \end{defn}

If $J$ is an \emph{admissible} almost complex structure, $J$ satisfies \eqref{C2} a.e. In particular, this gives that $|J|^2=2n$ a.e, hence it implies that $J\in L^\infty(M, \text{End}(TM))$ automatically. 
We shall use the notation $W^{1, 2}(\cJ_g)$ to denote all admissible almost complex structures. We also use $L^\infty(\cJ_g)$ to denote $L^\infty$ endomorphisms of $TM$ such that \eqref{C2} holds almost everywhere. Note that as a consequence of \eqref{C2} $L^\infty(\cJ_g)\subset W^{1, 2}(\cJ_g)$ as a set, but the topology induced by the $L^\infty$-norm and that by  $W^{1, 2}$-norm are evidently different. \\

Before we prove our main results in this section, we present a well-known fact about almost complex structures.
\begin{prop}\label{pm}Given an almost complex structure  $J_0: \R^{2n}\rightarrow \R^{2n}$. Then for $\epsilon$ small enough,  the $L^\infty$ neighborhood of $J_0$, given by $\{J: \|J-J_0\|_{L^\infty}<\epsilon,\}$ can be parametrized by matrices in a small ball centered at $S=0$ inside the linear subspace  $\{S: J_0S+SJ_0=0\}$.
If in particular $J_0$ is compatible with the Euclidean metric, then the $L^\infty$ neighborhood of $J_0$, given by $\{J: \|J-J_0\|_{L^\infty}<\epsilon, J+J^T=0\}$ can be parametrized by skew-symmetric matrices in a small ball such that $\{S: J_0S+SJ_0=0, S+S^T=0\}$.
\end{prop}

\begin{proof}Note that for any $J$ near $J_0$, we have that $J+J_0$ is invertible. Denote $S=(J+J_0)^{-1}(J-J_0)$, then $S$ is in a small ball around zero such that $SJ_0+J_0S=0$. On the other hand, for any $S$ in a small ball centered at $S=0$ such that $J_0S+SJ_0=0$, then $J_0(id+S)(id-S)^{-1}$ gives an almost complex structure near $J$. The two maps $J\rightarrow (J+J_0)^{-1}(J-J_0)$ and $S\rightarrow J_0(id+S)(id-S)^{-1}$ are invertible to each other, and gives a parametrization of almost complex structures near $J_0$. If in addition, we require that $J_0$ and $J$ are compatible with the Euclidean metric, then we will need to require $J, J_0, S$ are all skew-symmetric. 
\end{proof}

Proposition \ref{pm} can be directly applied to get the following results for almost complex structures on a manifold $J: TM\rightarrow TM$. Note that we only require for $L^\infty$ almost complex structures and we work in the $L^\infty$ norm. 

\begin{prop}\label{var1}Let $J_0: TM\rightarrow TM$ be an $L^\infty$ almost complex structure on $M$. Then almost complex structures in a small neighborhood of $J_0$ (in $L^\infty$ norm) can be parametrized by endomorphisms near zero such that $SJ_0+J_0S=0$. If almost complex structures are required to be compatible with a Riemannian metric $g$, then almost complex structures in a small neighborhood of $J_0$ can be parametrized by endomorphisms near zero such that $SJ_0+J_0S=0$, $g(S\cdot, \cdot)+g(\cdot, S\cdot)=0$.\end{prop}

Now we are ready to prove the following, 
\begin{thm}\label{equation1}A critical point of the energy functional $E(J)$ will satisfy the equation
\[
\Delta J-J\nabla_p J\nabla_p J =0
\]
in the following sense, for any $T\in L^\infty(M, \text{End}(M))\cap W^{1, 2}(M, \text{End}(M))$, 
\begin{equation}\label{H00}
\int_M \langle \nabla J, \nabla T\rangle  dv+\int_M \langle J\nabla_p J\nabla_p J , T\rangle dv=0
\end{equation}
Such a weak almost complex structure will be called a \emph{weakly harmonic almost complex structure}. 
\end{thm}
\begin{proof}The main difficulty is that $W^{1, 2}(\cJ_g)$ is not a linear subspace, hence does not inherit the local smooth structure automatically. To overcome this difficulty, first we consider the linearization of the constraints \eqref{C2} when $J$ is smooth. Taking a variation gives 
\begin{equation}\label{C3}
JS+SJ=0, g(Sx, y)+g(x, Sy)=0. 
\end{equation}
It is straightforward to see that, given $J_0\in \cJ_g$, the tangent space $T_{J_0}\cJ_g$ is given by the endomorphisms satisfying the constraints \eqref{C3}, which are clearly necessary conditions. To see this fact in a simple way, suppose $S$ is a smooth endomorphism satisfying \eqref{C3} such that $\|S\|_{L^\infty}$ is sufficiently small, then $id\pm S$ are both invertible. Let $J=J_0 (id+S) (id-S)^{-1}=(id-S)J_0(id-S)^{-1}$, then $J$ gives a compatible almost complex structure. 

Suppose $J_0\in W^{1, 2}(\cJ_g)$ is a critical point of $E(J)$. Let $S\in L^\infty\cap W^{1, 2}$ satisfying \eqref{C3}. For $|t|$ sufficiently small such that $id-tS$ is invertible, consider $J(t)=J_0(id+tS)(id-tS)^{-1}$. Note that $J(t)\in L^\infty\cap W^{1, 2}$, hence by Proposition \ref{var1} $J(t)\in W^{1, 2}(\cJ_g)$ for $|t|$ small.  At $t=0$, we have, 
\[
\frac{d}{dt} E(J(t))=0. 
\]
A direct computation gives
\begin{equation}\label{H0}
\int_M \langle \nabla (J_0S), \nabla J_0\rangle =0
\end{equation}
If $J_0$ is smooth, it follows that 
\[
\int_M \langle S, J_0\Delta J_0\rangle =0  
\]
Since $SJ_0+J_0S=0$ and $g(Sx, y)+g(x, Sy)=0$, we have 
\[
\int_M \langle S, J_0\Delta J_0-\Delta J_0 J_0\rangle =0
\]
Note that $J_0\Delta J_0-\Delta J_0 J_0$ satisfies the constraints \eqref{C3}. It follows that 
\[
[\Delta J_0, J_0]=0.
\]
An equivalent form, by applying $\Delta (J_0^2)=0$, is given by
\[
\Delta J_0-J_0\nabla_pJ_0\nabla_pJ_0=0,
\]
where $J\nabla_p J\nabla_p J=J_j^a\nabla_p J_a^b\nabla_q J_b^k g^{pq}=J^k_b\nabla_pJ^b_a\nabla_qJ^a_j g^{pq}$ in local coordinates. 
When $J_0$ is not smooth (we write $J=J_0$ for simplicity), for any $T=(T^j_i)\in L^\infty\cap W^{1, 2}$, we consider $S=(S^j_i)$ given by
\[
S^j_i=(T^j_i+J^a_iT^b_aJ^j_b)-g_{il}\left(T^l_k+J^a_kT^b_aJ^l_b\right)g^{kj}.
\]
Then $S\in L^\infty\cap W^{1, 2}$ and satisfies the constraints \eqref{C3}. By taking $g_{ij}=\delta_{ij}$ at one point, we see that
$\langle (\nabla J) S, \nabla J\rangle $=0. It follows from \eqref{H0} that
\[
\int_M \langle J\nabla S, \nabla J\rangle=\int_M \langle \nabla S, J\nabla J\rangle=0
\]
Again take $A=T+JTJ$ and $g_{ij}=\delta_{ij}$ at a point, we have $S^j_i=A^j_i-A^i_j$. It follows that (note that $J^j_i+J^i_j=0$ at the point), 
\[
\int_M \langle \nabla A, J\nabla J\rangle=0
\]
Take $\tilde T=TJ$, then $A=J\tilde T-\tilde TJ$. By a direct computation using the above, we get
\[
\int_M \langle \nabla J, \nabla \tilde T\rangle  dv+\int_M \langle J\nabla_p J\nabla_p J , \tilde T\rangle dv=0
\]
It completes the proof. 
\end{proof}

In this statement, even though we choose the testing matrix $S$ satisfying the constraints \eqref{C3}, but \eqref{H00} holds for any matrix $T\in L^\infty\cap W^{1, 2}$. We can also consider the minimizing problem $E(J)$ for all almost complex structures $J$, which are not necessarily compatible with $g$. In particular if $J$ is a critical point of $E(J)$ among all almost complex structures, and if $J$ is compatible with $g$, then the Euler-Lagrangian equation is also \eqref{H00}. 

Next we prove the existence of a minimizer of $E(J)$ on any compact almost Hermitian manifold. This follows from the standard argument in calculus of variations. 
\begin{thm}\label{existence1}Let $(M, g)$ be a compact almost Hermitian manifold with a compatible almost complex structure. Then there exists at least one minimizer of $E(J)$ over $W^{1, 2}(\cJ_g)$. 
In particular it is a weak solution of \eqref{H1} in the sense of \eqref{H00}. Moreover, for any sequence of energy-minimizers $\{J_i\}$, there exists a subsequence $\{J_{k_i}\}$ and an energy-minimizing harmonic almost complex structure $J$ such that $J_{k_i}\rightarrow J$ strongly in $W^{1, 2}$. \end{thm}

\begin{proof}Let 
\[
E_0=\inf\{E(J): J\in W^{1, 2}(\cJ_g)\}
\]
Let $\{J_k\}$ be a minimizing sequence. Then $J_k$ has bounded $W^{1, 2}$ norm. By the weak compactness of $W^{1, 2}(M, \text{End}(TM))$, there exists a subsequence, which we still denote by $\{J_k\}$ for simplicity, such that $J_k$ converges weakly  to $J_0$ in $W^{1, 2}$. By Rellich's compactness theorem, 
\[
\lim_{k\rightarrow \infty}\|J_k-J_0\|_{L^2}=0.
\]
In particular $J_k$ converges to $J_0$ almost everywhere. Since $J_k$ satisfies the constraints \eqref{C3} a.e, it follows that $J_0$ also satisfies the constraints \eqref{C3} a.e. It follows that $J_0\in W^{1, 2}(\cJ_g)$. In particular, $E(J_0)\geq E_0$. On the other hand, by the lower semicontinuity, 
\[
E(J_0)\leq \lim_{k\rightarrow \infty} E(J_k)=E_0.
\] 
Therefore we have $E(J_0)=E_0$ and $J_0$ is a minimizer. To prove the strong convergence, we consider
\[
\int_M |\nabla J_k-\nabla J_0|^2dv=\int_{M}\left(|\nabla J_k|^2+|\nabla J_0|^2\right)-2\int_M \langle \nabla J_k, \nabla J_0\rangle
\]
By taking $k\rightarrow \infty$, and using the fact that $J_k$ converges to $J_0$ weakly in $W^{1, 2}$, it follows that 
\[
\int_M  |\nabla J_k-\nabla J_0|^2dv\rightarrow 0.
\]
Hence $J_k$ converges to $J_0$ strongly in $W^{1, 2}$. If $J_k$ is a sequence of minimizers, the above arguments directly apply to conclude that there exists a subsequence $J_{k_i}$ converges to a minimizer $J$ and the convergence is strong in $W^{1, 2}$. In other words, the minimizing harmonic almost complex structures form a compact subset in $W^{1, 2}$. 
\end{proof}

\begin{defn}
We introduce two notations, 
\begin{equation}
\begin{split}
E_0(g)&=\min \{E(J): J\in W^{1, 2}(\cJ_g)\}, \\
 E^s_0(g)&=\inf\{E(J): J\in \cJ_g\}.
\end{split}
\end{equation}
\end{defn}

In other words, $E_0(g)$ is the minimal energy level for all \emph{admissible} almost complex structures and $E^s_0(g)$ is the minimal energy level for all \emph{smooth} compatible almost complex structures. 
Clearly we have $E_0(g)\leq E^s_0(g)$ and $E_0(g)$ is always obtained. We should emphasize that in general it is not clear whether $E_0(g)=E^s_0(g)$, since $\cJ_g$ might  not be dense in $W^{1, 2}(\cJ_g)$. If $E_0(g)<E^s_0(g)$, this is called \emph{energy gap}. 
 In the case of harmonic map, energy gap is a well-known phenomenon, see \cite{SacksU, HL1, GMS, Lin1} and in particular the survey paper \cite{Hardt1}. We would like to ask,
\begin{pl}When is there no energy gap, that is $E_0(g)=E^s_0(g)$?
In particular, if there exists a smooth almost complex structure such that it is a minimizer among smooth almost complex structures, is it a minimizer in $W^{1, 2}(\cJ_g)$?
\end{pl}

It seems plausible if there exists a smooth minimizer among smooth almost complex structures, then there is no energy gap. We can show that (given the regularity results we will establish later) this is the case for  the energy-minimizing examples considered in \cite{BorSalvai}. 
As a trivial example, if $E_0(g)=0$, we have the following simple observation. 
\begin{cor} $E_0(g)=0$ if and only if $g$ is K\"ahler with respect to an energy minimizer. 
\end{cor}

\begin{proof}If $g$ is K\"ahler, the clearly $E_0(g)=E^s_0(g)=0$. If  $E_0(g)=0$, then there exists $J_0\in W^{1, 2}(\cJ_g)$ such that $E(J_0)=0$. It follows that $\nabla J_0=0$ almost everywhere. Hence $J$ is in particular Lipschitz. This implies that $J_0$ is indeed smooth. Hence $(g, J_0)$ defines the required K\"ahler structure. 
\end{proof}
This simple corollary gives a gap of energy if $g$ is not a K\"ahler metric (with respect to any compatible (almost) complex structures). 
\begin{cor} Let $(M, g)$ be a compact almost Hermitian manifold. Suppose $g$ is not a K\"ahler metric, then there exists $\epsilon_0>0$ such that $E_0(g)\geq \epsilon_0$. 
\end{cor}

\begin{rmk}It is also straightforward to talk about the manifolds with boundary and minimizing problem with fixed boundary data, as the classical Dirichlet's problem. 
We shall mainly focus on the compact case. 
\end{rmk}

An immediate question is that the regularity of an energy-minimizing weakly harmonic almost complex structures in $W^{1, 2}(\cJ_g)$.
It is well-known that a continuos weakly harmonic map is smooth; certainly this is the case for weakly harmonic almost complex structures. We will derive this result in the appendix. 
In general a weakly harmonic almost complex structure is not necessarily smooth.
We study the regularity theory of energy-minimizing almost complex structures in the following sections, starting  with a monotonicity formula for energy minimizing almost complex structures. 

\section{Monotonicity formula for energy minimizing almost complex structures}
Monotonicity formula plays an essential role in regularity theory of harmonic maps.  It was first proved by Schoen-Uhlenbeck \cite{SU} for minimizing harmonic maps, using comparison maps directly. Later Price \cite{P} proved the monotonicity formula which works for more general harmonic maps by considering domain variations, called \emph{stationary harmonic maps}. However, both methods in harmonic maps do not work in our case directly, due to the restriction that  almost complex structures are compatible with the metric $g$.  

To derive the monotonicity formula for energy-minimizing harmonic almost complex structures, we need to construct various comparison almost complex structures. An important ingredient is an elementary matrix decomposition for matrices satisfying $N^2=-id$, see Proposition \ref{matrix}. This new ingredient is essential for us to modify the approach in \cite{SU} for construction of various comparison almost complex structures. 
We should mention that if the objects are smooth (harmonic map or harmonic almost complex structure), the monotonicity formula follows directly from the structural equation in a straightforward way. But the main point  is to establish monotonicity formula for non-smooth objects. 
As in Schoen-Uhlenbeck \cite{SU}, the key is to choose various comparison almost complex structures accordingly. However the cone-type  of $J_r(x)=J(rx/|x|)$ as in \cite{SU} doesn't work in our situation since in general it is not compatible with the metric $g(x)$.

The theory of stationary harmonic maps plays an important role in harmonic map theory. But we do not have a direct companion for harmonic almost complex structure. We will briefly discuss plausible notion of \emph{stationary weakly harmonic almost complex structures} in appendix.

Since the discussions are mainly local, we assume that the injectivity radius of the $g$ is larger than $1$, by a definite scaling if necessary. Hence for any $p$, the exponential map $\exp_p: B=B(1)\subset \R^{2n}\rightarrow B_1(p)\subset M$ is a diffeomorphism. When we work locally inside $B_1(p)$, we can then consider $g\circ \exp=\exp^* g, J\circ \exp=\exp^*J$ on $B$, which we will still denote as $g$ and $J$ for brevity. 
We can assume that on $B$, $g_{ij}(0)=\delta_{ij}$, $\p g (0)=0$, and 
\[ (1-\delta |x|^2)\delta_{ij}\leq g(x)\leq (1+\delta |x|^2)\delta_{ij}, |\p g|(x)\leq \delta |x|, |\p^2 g|\leq \delta\] for some small but fixed constant $\delta$.  We can assume this holds for all points on $M$ by a definite scaling if necessary.  When we need to work on different scalings $(x\rightarrow Rx)$, we can assume the scaling invariant bound,
\[
(1-\delta R^{-2}|x|^2)\delta_{ij}\leq g(x)\leq (1+\delta R^{-2}|x|^2)\delta_{ij}, R |\p g|(x)\leq \delta  R^{-1}|x|, R^{2}
|\p^2 g|\leq \delta \] 

The following is straightforward, comparing the energy measured by $g$ and by Euclidean metric $\delta_{ij}$ over $B$. 
We fix some notations. We use  $\nabla$ to denote the metric covariant directive, and $D, \p$ to denote derivatives with respect to local coordinate. We use $\int f dv$ for the integration with the metric measure, and  $\int f:=\int f dx$ for the integration with Euclidean measure. Denote $B_r\subset B$ and $B_r(p)\subset B_1(p)$ for $0<r\leq 1$. Note that $\exp_p: B_r\rightarrow B_r(p)$.  We denote $c_n$ to be a uniformly bounded dimensional constant.
\begin{prop}\label{c1}We have the following,
\begin{equation}\label{comparison1}
\begin{split}
\int_{B_r} |\nabla J|^2dv\geq &(1-c_n \delta r^2)\int_{B_r} |D J|^2 -c_n\delta^2 r^m \\
\int_{B_r} |\nabla J|^2dv\leq &(1+c_n \delta r^2)\int_{B_r} |D J|^2 +c_n\delta^2 r^m
\end{split}
\end{equation}
\end{prop}
\begin{proof}It is straightforward to get
\begin{equation}\label{comparison2}
\begin{split}
\int_{B_r} |\nabla J|^2dv\geq &(1-c_n \delta r^2)\int_{B_r} |D J|^2 -c_n\delta r\int_{B_r} |D J|-c_n\delta^2 r^m \\
\int_{B_r} |\nabla J|^2dv\leq &(1+c_n \delta r^2)\int_{B_r} |D J|^2 +c_n\delta r\int_{B_r} |D J|+c_n\delta^2 r^m,
\end{split}
\end{equation}
 using the fact that $\nabla J=D J+\Gamma*J$ in $B$ and $g$ is close to the Euclidean metric, in the sense that
\[
(1-\delta r^2)\delta_{ij}\leq g_{ij}\leq (1+\delta r^2)\delta_{ij}, |\p g|\leq \delta r. 
\]
Applying a Cauchy-Schwarz inequality to the term  $\int |D J|$, we can get \eqref{comparison1}. 
\end{proof}
 
First we consider a special case, that the background metric $g$ is locally conformally flat. We assume for any given point $p$, there exists a small neighborhood $U$ of $p$, such that there exists a positive function $u$, $g_{ij}=e^u \delta_{ij}$ over $U$. We assume that $B_1(p)\subset U$. The assumption guarantees that the comparison almost complex structure $J_r(x)$ is compatible with $g(x)$.  Then the argument as in \cite[Proposition 2.4]{SU} can be directly applied.

\begin{prop}\label{conformalflat}Suppose $g$ is locally conformally flat. Let $J$ be an energy minimizing almost complex structure, then we have the following monotonicity formula, for $r<R\leq 1$
\begin{equation}
r^{2-m}\int_{B_r}|\nabla J|^2 dv\leq c_0 R^{2-m}\int_{B_R}|\nabla J|^2dv+c_n\delta^2 R^2,
\end{equation}
where $c_0=1+c_n\delta R$.
\end{prop}

\begin{proof} For any $0<r\leq 1$, denote
\[
J_r(x)=\left\{\begin{matrix}
&J(x), |x| \geq r\\
&J(\frac{rx}{|x|}), |x|<r
\end{matrix}
\right.
\]
Then clearly $J_r(x)\in L^\infty\cap W^{1, 2}$. Since $g=e^u\delta_{ij}$ over $B_1(p)$, it follows that $J^2=-id$ and $J+J^t=0$. Hence $J_r$ is then $g$-compatible. Since $J$ is energy minimizing ($J=J_r$ outside $B_r$), we have
\begin{equation}\label{m21}
\int_{B_r} |\nabla J|^2 dv\leq \int_{B_r} |\nabla J_r|^2 dv
\end{equation}
Using \eqref{comparison1}, we have
\begin{equation}\label{m1}
\begin{split}
\int_{B_r} |\nabla J_r|^2 dv\leq &(1+c_n \delta r^2)\int_{B_r} |D J_r|^2 +c_n\delta^2 r^m\\
\int_{B_r} |\nabla J|^2 dv\geq &(1-c_n \delta r^2)\int_{B_r} |D J|^2 -c_n\delta^2 r^m
\end{split}
\end{equation}
We compute
\begin{equation}\label{m2}
\begin{split}
\int_{B_r} |D_x J_r|^2dx =&\int_0^r d\rho \int_{\p B_\rho}|D_x J_r|^2 dS_\rho(x)\\
=&\int_0^rd\rho \int_{\p B_r} (|D_y J|^2-|\p_r J|^2) dS_r(y) r^{3-m} \rho^{m-3}\\
=&\frac{r}{m-2}\int_{\p B_r}\left(|D J|^2-|\p_r J|^2\right)dS_r,
\end{split}
\end{equation}
where we have used the fact that $(y=rx/|x|)$\[
\rho^{1-m}dS_\rho(x)=r^{1-m}dS_r(y),\;\; \rho^2 \sum_i|\p_{x_i} J_r|^2=r^2 \left(\sum_i|\p_{y_i} J|^2-|\p_r J|^2\right)
\]
Using \eqref{m21} together with \eqref{m1} and \eqref{m2} we obtain, 
\begin{equation}
\int_{B_r} |DJ|^2\leq (1+c_n\delta r^2) \frac{r}{m-2}\int_{\p B_r} (|DJ|^2-|\p_rJ|^2)dS_r+c_n \delta^2 r^m
\end{equation}
We rewrite this as follows,
\begin{equation}\label{m22}
\int_{\p B_r} |\p_rJ|^2dS_r\leq \int_{\p B_r} |DJ|^2dS_r-\frac{m-2}{r(1+c_n\delta r^2)}\int_{B_r}|DJ|^2+c_n\delta^2 r^{m-1}.
\end{equation}
If $\delta=0$ (Euclidean case), this leads to the monotonicity formula
\[
0\leq r^{2-m}\int_{\p B_r} |\p_r J|^2 dS_r\leq \frac{\p}{\p r}\left(r^{2-m}\int_{B_r}|DJ|^2\right)
\]
In general, we can get that, 
\begin{equation}\label{m23}
0\leq e^{c_n\delta r} r^{2-m}\int_{\p B_r} |\p_r J|^2 dS_r\leq \frac{\p}{\p r}\left(e^{c_n\delta r} r^{2-m}\int_{B_r}|DJ|^2+e^{c_n\delta r} c_n \delta^2 r^2\right)
\end{equation}
A direct integration with \eqref{comparison1} leads to the desired monotonicity (assuming $\delta<<1$), 
\begin{equation}
r^{2-m}\int_{B_r}|\nabla J|^2 dv\leq c_0 R^{2-m}\int_{B_R}|\nabla J|^2dv+c_n\delta^2 R^2,
\end{equation}
where $c_0=1+c_n\delta R$.  
\end{proof}

Next we deal with the general situation. We need the following elementary fact.

\begin{prop}\label{matrix}Let $N$ be an invertible matrix such that $N^2=-id$. We can construct a unique $\bar N$ out of $N$ such that $\bar N^2=-id$ and $\bar N$ is skew symmetric.
\end{prop}

\begin{proof}We write $N=S+A$,
where $S$ is the symmetric part of $N$ and $A$ is the skew symmetric part. In other words,
\[
2S=N+N^t, 2A=N-N^t. 
\]
Then we have $AS+SA=0$ and $S^2+A^2=-id$. In particular $A^2=-(id+S^2).$  Note that $id+S^2$ is positive definite symmetric, there exists a unique positive definite symmetric matrix $Q$ such that $Q^2=id+S^2$. In other words, $Q$ is the square root of $id+S^2$. Since $A$ commutes with $id+S^2$, $A$ commutes with $Q$ ($Q$ can be written as a polynomial of $id+S^2$). Denote $\bar N=Q^{-1}A$, then $\bar N^2=-id$ and $\bar N$ is skew-symmetric. Such $\bar N$ is evidently unique. 
\end{proof}

Proposition \ref{matrix} is a pointwise construction. It has an interesting consequence, which seems not to be noticed in literature. 

\begin{prop}\label{matrix1}For an almost complex manifold, any Riemannian metric is an almost Hermitian metric with some compatible almost complex structure. 
\end{prop}

\begin{proof}The proof is local in nature. We need a notion for $g$-symmetric matrix in local coordinates. Fix a local coordinate, and we identify any endomorphism of $TM$ with its local matrix representation. For a matrix $N$, we write $N^t$ to be its transpose and $N^t_g=gN^tg^{-1}$ to be its $g$-transpose. In particular, $N^t_g$ can be defined by
\[
g(N\cdot, \cdot)=g(\cdot, N^t_g \cdot). 
\] 
Then $N$ is called $g$-symmetric if $N^t_g=N$ ($gN^tg^{-1}=N$) and $N$ is called $g$-skew symmetric if $N^t_g+N=0$.  Now suppose $N^2=-id$. 
The same argument in Proposition \ref{matrix} then applies, using $g$-symmetric and $g$-skewsymmetric to replace symmetric and skew-symmetric. To be more precise, we can write $g=GG^t$ for a nondegenerate matrix $G$. Then $g$-transpose $N_g^t$ satisfies
\[
G^{-1}N^t_g G=(G^{-1}NG)^t
\]
In particular $N$ is $g$-symmetric if and only if $G^{-1}NG$ is symmetric. 
We write $N=S_g+A_g$, where $2S_g=N+N^t_g$ is the $g$-symmetric part and $A_g$ is the $g$-skew symmetric part. Denote $S=G^{-1}S_gG$ and $A=G^{-1}A_gG$. Then
\[
2S=G^{-1}NG+G^{-1}N^t_gG=G^{-1}NG+(G^{-1}NG)^t
\]
Hence $S$ is the symmetric part of $G^{-1}NG$ and $A$ is the skew-symmetric part. Applying Proposition \ref{matrix} to $G^{-1}NG$, we can get a symmetric positive matrix $Q$ such that $Q^2=id+S^2=-A^2$. Hence $Q^{-1}A$ is a skew-symmetric matrix such that $(Q^{-1}A)^2=-id$. Denote $\bar N_g=GQ^{-1}AG^{-1}$, then $\bar N_g$ satisfies $\bar N_g^2=-id$ and $\bar N_g$ is $g$-skewsymmetric.
If we denote $Q_g=GQG^{-1}$, then $Q_g$ is the square root of $id+S_g^2$ and it is $g$-symmetric. Then $\bar N_g=Q_g^{-1}A_g$. 
Clearly if $N$ is a matrix-valued function, then $\bar N_g$ is smooth if $N$, $g$ are both smooth. 

\end{proof}

If $J^2=-id$ is an almost complex structure, we can then construct a unique almost complex structure $\bar J_g$ using the process above. Then $\bar J_g$ is $g$-compatible almost complex structure. 
Note that if $J\in L^\infty\cap W^{1, 2}$, then $\bar J\in W^{1, 2}(\cJ_g)$.

Let $J$ be an energy-minimizer. We consider $J_r(x)$ as above, which might not be $g(x)$-compatible in general. We construct $\bar J_r(x)$ by Proposition \ref{matrix1} out of $J_r(x)$ for each $x$. We will use $\bar J_r(x)$, which is $g$-compatible, as the comparison almost complex structure. The main point is that $J_r(x)$ might not be $g(x)$-compatible, but it is $g_r(x):=g(rx/|x|)$ compatible. While over $B(1)$, we assume that $g$ is close to an Euclidean metric, in particular $g(x)$ is closed to $g_r(x)$. In other words, $J_r(x)$ is compatible with $g_r(x)$, hence it is \emph{almost} compatible with $g(x)$. 
This would imply that $\bar J_r$ is almost identical to $J_r$.  
To be precise, 
\begin{prop}\label{c3}We have the following estimate
\begin{equation}\label{comparison3}
\begin{split}
\int_{B_r} |\nabla J_r|^2dv\geq &(1-c_n \delta r^2)\int_{B_r} |\nabla \bar J_r|^2 dx-c_n\delta^2 r^m \\
\int_{B_r} |\nabla J_r|^2dv\leq &(1+c_n \delta r^2)\int_{B_r} |\nabla \bar J_r|^2 dx+c_n\delta^2 r^m
\end{split}
\end{equation}
\end{prop}

\begin{proof} 
Consider $J_r=J(rx/|x|)$. We apply Proposition \ref{matrix1} to obtain $\bar J_r(x)$ (for $x\neq 0$), which is $g(x)$-compatible. 
We want to prove that $\bar J_r$ satisfies \eqref{comparison3}. 

Let $S_g$ be $g$-symmetric part of $J_r$, and $A_g$ be $g$-skewsymmetric part. And let $Q_g$ be the square root of $id+S_g^2$, in the sense of Proposition \ref{matrix1}. Then $\bar J_r=Q_g^{-1}A_g$.

First it is easy to check that $\bar J_r\in W^{1, 2}(\cJ_g)$ if $J\in W^{1, 2}(\cJ_g)$. 
Denote $g_r(x)=g(rx/|x|)$ for $x\in B_r, x\neq 0$. 
Note that $J_r(x)=J(rx/|x|)$ is $g_r(x)$ compatible, hence we have $J_r+g_rJ_r^t g_r^{-1}=0$. 
We compute
\[
2S_g=J_r+gJ_r^t g^{-1}=gJ_r^tg^{-1}-g_r J_r^t g_r^{-1}=(g-g_r)J_r^tg^{-1}+g_r J^t_r (g^{-1}-g^{-1}_r).
\]
Hence it follows that $|S_g|\leq c_n \delta r^2$ and $|\nabla S_g|\leq c_n \delta r(1+r|\nabla J_r|)$. Note that $A_g=J_r-S_g$, $|Q_g^2-id|\leq |S_g|^2$.
We also have
 \begin{equation}\label{matrix018}|\nabla Q_g|\leq |S_g||\nabla S_g|\end{equation}
Since $|S_g|\ll 1$, we can write
\[
Q_g=(id+S_g^2)^{1/2}=\sum_{n\geq 0} \binom{1/2}{n} S_g^{2n}.
\]
Hence we get 
\[
|\nabla Q_g|=\sum_{n\geq 0} 2n \binom{1/2}{n} |S_g|^{2n-1} |\nabla S_g|=\frac{|S_g||\nabla S_g|}{1+|S_g|^2}\leq |S_g| |\nabla S_g|,
\]
which establishes \eqref{matrix018}.
We can then estimate (pointwise)
\[
|\nabla J_r|(1-c_n \delta r^2)-c_n \delta r\leq |\nabla \bar J_r|\leq |\nabla J_r|(1+c_n \delta r^2)+c_n \delta r.
\]
This completes the proof of \eqref{comparison3} by direct integration. 
\end{proof}

Now we are ready to state and prove the main theorem in this section. 
\begin{thm}\label{monotonicity}Let $J\in W^{1, 2}(\cJ_g)$ be an energy minimizing harmonic almost complex structure, then we have the following monotonicity formula, for $r<R\leq 1$, 
\begin{equation}
r^{2-m}\int_{B_r}|\nabla J|^2 dv\leq c_0 R^{2-m}\int_{B_R}|\nabla J|^2dv+c_n\delta^2 R^2,
\end{equation}
where $c_0=1+c_n\delta R$. 
\end{thm}
\begin{proof}By energy minimizing, we have the following,
\[
\int_{B_r} |\nabla J|^2 dv\leq \int_{B_r} |\nabla \bar J_r|^2dv
\]
Now applying Proposition \ref{c1} and Proposition \ref{c3}, we argue exactly as in Proposition \ref{conformalflat} to get \eqref{m23}. This completes the proof. 
\end{proof}

Denote 
\begin{equation}\label{density5}
\tilde \Theta_{J}(p, r)=e^{c_n\delta r}\left(r^{2-m}\int_{B_r(p)}|DJ|^2+c_n\delta^2 r^2\right)
\end{equation}
By \eqref{m23}, we know that $\tilde \Theta_J(p, r)$ is increasing in $r$, 
hence  $\lim_{r\rightarrow 0} \tilde \Theta_J(p, r)$ exists. 
We define two normalized energy quantities, for $r\in (0, 1]$, 
\begin{equation}\label{density6}
\begin{split}
\Theta^e_J(p, r)=&r^{2-m}\int_{B_r(p)}|DJ|^2\\
\Theta_J(p, r)=& r^{2-m}\int_{B_r(p)}|\nabla J|^2 dv
\end{split}
\end{equation}
It is then straightforward to derive the following, using \eqref{comparison1}
\begin{equation}\label{density3}
\lim_{r\rightarrow 0}\tilde \Theta_{J}(p, r)=\lim_{r\rightarrow 0} \Theta^e_J(p, r)=\lim_{r\rightarrow 0} \Theta_J(p, r)
\end{equation}
The monotonicity formula leads to the definition of the following density function
\begin{defn}For any $p\in M$, the density function is defined to be
\begin{equation}\label{density2}
\Theta_J{(p)}=\lim_{r\rightarrow 0} \Theta_J(p, r)=\lim_{r\rightarrow 0} r^{2-m}\int_{B_r(p)}|\nabla J|^2dv.
\end{equation}
\end{defn}
By \eqref{density3}, the density function is independent of the metric $g$. Since in a very small scale, the metric $g$ is almost identical to the Euclidean metric (after scaling), but the normalized energy is scaling invariant. Hence the density function can be defined using the Euclidean metric, as indicated by \eqref{density3}.  The density function $\Theta_J(p)$ will play an important role in the study of properties of singular set in the following sections. 
 
\begin{rmk}\label{remark1}For a weakly harmonic almost complex structure, the main difficulty to derive a monotonicity formula lies in the fact that for $J\in L^\infty\cap W^{1, 2}$, and a smooth diffeomorphism $\phi_t$ with $\phi_0=id$, the one-parameter family $\phi_t^* J$ is not a small perturbation of $J$ in $L^\infty$ norm. 
For energy-minimizing harmonic map, Schoen-Uhlenbeck  constructed $u_r=u(rx/|x|)$ as a comparison map. Note that $u_r(x)$ might not be a small perturbation of $u$ in general, no matter how small $r$ is.  Price \cite{P}  consider the variation for $({\phi^{-1}_t}^*g, J)$ instead of $(g, \phi_t^* J)$. Note that the metric $g$ is smooth and ${\phi^{-1}_t}^*g$ is a smooth family of $g$. Price's approach  leads to the definition of stationary harmonic map. In our situation, both methods in Schoen-Uhlenbeck and Price do not apply directly. The difficulty is that the almost complex structure $J_r$  is not $g(x)$-compatible in general, and is not a small perturbation of $J$ in $L^\infty$ norm either, no matter how small $r$ is. To overcome the difficulty, we construct $\bar J_r$ using Proposition \ref{matrix1}. The key point is $\bar J_r$ is $g$-compatible and  it is  a small perturbation of $J_r$ in $L^\infty\cap W^{1, 2}$. 
\end{rmk}

\section{Regularity of energy-minimizing harmonic almost complex structure}

In this section we derive regularity result for energy-minimizing harmonic almost complex structures. Our argument follows closely the regularity theory of energy-minimizing harmonic map by Schoen-Uhlenbeck\cite{SU} (we also follow Simon's book \cite{Simon} and Lin-Wang's book \cite{LW}). The similar structure of the elliptic system (together with the monotonicity formula) makes the argument very similar to \cite{SU}.  As in \cite{SU}, the key is to construct various comparison almost complex structures. 
We should also emphasize there are several notable differences and we need extra care to construct various comparison almost complex structures. We will present the arguments in detail.

We fix some notations. Let $g_{ij}$ be a Riemannian metric on the unit ball $B$ in $\R^{2n}$ which is close to an Euclidean metric in $C^3$ in the following sense: there exists a small positive constant $\delta>0$ such that, for $x\in B$, 
\begin{equation}\label{delta}(1-\delta |x|^2)\delta_{ij}\leq g_{ij}(x)\leq (1+\delta |x|^2)\delta_{ij}, |\p g|(x)\leq \delta |x|, |\p^2 g|+|\p^3 g|\leq \delta.\end{equation}
Here we should emphasize that $\delta<<1$ is a dimensional constant.  
In this section $C$ denotes a positive uniformly bounded dimensional constant (if not mentioned otherwise) which might differ line by line. We view $J=(J^j_i)$ both as a tensor when we need to work in geometric terms, such as $\nabla J$. We also view $J$ as a matrix-valued function  in $W^{1, 2}\cap L^\infty$ satisfying the constraint conditions $J_i^kJ^j_k=-\delta_{ij}, g_{kl}J_i^kJ^l_j=g_{ij}$. 

\subsection{A priori estimate for a smooth harmonic almost complex structure}In this section we consider smooth harmonic almost complex structures and derive some a priori estimates. These estimates would be useful for the study of moduli of harmonic almost complex structures. In the smooth case the arguments are almost identical to the case for harmonic maps. We include details for completeness. 
We have the following, 

\begin{thm}There exists $\epsilon$ such that if, for $r<1$, 
\[
r^{2-m}\int_{B_r} |\nabla J|^2<\epsilon,
\]
then  we have
\begin{equation}\label{estimate1}
r^2\sup_{B_{r/4}} |\nabla J|^2\leq c_n r^{2-m}\int_{B_r} |\nabla J|^2.
\end{equation}
Moreover, for any $k\geq 1$, we have the estimates, 
\begin{equation}\label{estimate2}
r^{2k}\sup_{B_{r/8}} |\nabla^kJ|^2\leq C r^{2-m}\int_{B_r} |\nabla J|^2.
\end{equation}

\end{thm}

\begin{proof}First we prove \eqref{estimate1}.  The proof is almost identical as in harmonic map case. We roughly present the proof here for completeness. Denote $u(x)=|\nabla J|^2(x)$. Take $r_1=r/2$. There exists $\sigma_0\in [0, r_1)$ such that
\begin{equation}\label{point1}
(r_1-\sigma_0)^2 \sup_{B_{\sigma_0}} u=\max_{\sigma\in [0, r_1]} (r_1-\sigma)^2 \sup_{B_\sigma} u.
\end{equation}
Moreover, let $x_0\in \bar{B}_{\sigma_0}$ such that 
\[
e_0=u(x_0)=\sup_{B_{\sigma_0}} u
\]
Set $\rho=(r_1-\sigma_0)/2$. Then by the choice of $\sigma_0$, we have $B_\rho(x_0)\subset B_{\rho+\sigma_0}$, and 
\[
\sup_{B_\rho(x_0)}u\leq \sup _{B_{\rho+\sigma_0}} u\leq 4e_0. 
\]
If we have $e_0\leq 4\rho^{-2}$, then we have, by \eqref{point1},
\[
\max_{\sigma\in [0, r_1]} (r_1-\sigma)^2 \sup_{B_\sigma} u\leq 16. 
\]
By taking $\sigma=r_1/2$, we have proved the desired result. Otherwise consider \[\tilde J(y)=J\left(\frac{y}{\sqrt{e_0}}+x_0\right)\]
Now let $v: B_{r_0}\rightarrow \R$, with $r_0=\rho\sqrt{e_0}$, by
$v(y)=|\nabla \tilde J|^2$. Note that $r_0\geq 2$ in this case. 
Then we have
\[
v(0)=1, \sup_{B_{r_0}}v \leq 4
\]
By the Bochner formula \eqref{bochner} proved in Appendix, we have
\[
\Delta u\geq -C(u^2+\sqrt{u})
\]
Denote $\tilde g(y)=g(y/\sqrt{e_0}+x_0)$ and $\tilde \Delta$ to be the corresponding Laplacian operator. Then we have the following, 
\[
\tilde \Delta v\geq -C(v^2+e_0^{-3/2}\sqrt{v})\geq -C(v+e_0^{-3/2}\sqrt{v}).
\]
By the standard elliptic estimate (see Theorem 8.17 \cite{GT}), we have
\[
1=v(0)\leq C\int_{B_1} v+Ce_0^{-3/2}.
\]
We compute 
\[
\int_{B_1} v=\sqrt{e_0}^{n-2}\int_{B_{{e^{-1/2}_0}}(x_0)} u
\]
Now we apply the monotonicity formula for $B_{{e^{-1/2}_0}(x_0)}\subset B_{r/2}(x_0)$,
\[
\sqrt{e_0}^{n-2}\int_{B_{{e^{-1/2}_0}}(x_0)} u\leq c_0 2^{m-2} r^{2-m}\int_{B_{r/2}(x_0)} u+c_n \delta^2 r^2\leq C \epsilon+c_n\delta^2r^2. 
\]
It then  follows that
\[
1=v(0)\leq C\epsilon+C\delta^2r^2+Ce_0^{-3/2}. 
\]
By our choice of $\epsilon$ and $\delta$, we can assume that $C\epsilon+C\delta^2r^2\leq 1/2$. Hence it follows that $e_0\leq C$. Again by \eqref{point1} and take $\sigma=r_1/2$, it completes the proof for $k=1$.  For $k\geq 2$, \eqref{estimate2} follows from the standard linear elliptic regularity and boot-strapping argument. 
 \end{proof}

\subsection{An $\epsilon$-Regularity of energy-minimizing almost complex structure}
In this section we derive the regularity theory for energy-minimizing harmonic almost complex structures. Our argument is a direct adaption of Schoen-Uhlenbeck \cite{SU} for energy-minimizing harmonic maps (see the presentations in Simon's book \cite{Simon}, Chapter 2).  There are a few differences which will need modifications. For example, the compatible condition of $J$ with the background metric $g$ plays an important role. As an tensor-valued function, there are a few places the curvature and its derivatives would appear when we commute the covariant derivatives of $J$. This will introduce extra terms. More importantly, $\nabla u=0$ means $u$ is a constant map, while $\nabla J=0$ simply means $J$ is a compatible almost complex structure which defines a K\"ahler structure with $g$. This will give some differences for the Poincare inequality involved with $J$ and $\nabla J$. Moreover, as a tensor-valued function $J$, sometimes we need to treat it as a matrix-valued function, and there are differences when we take derivaives. For example, 
We need to  solve a Laplacian equation of matrix-valued function (for each component $K^\beta_\alpha$),
\[
\Delta K=\left(g^{ij}\frac{\p^2}{\p x^i\p x^j}+\Gamma^k_{ij}\p_k\right) K^\beta_\alpha=0,
\]
where $K=(K^\beta_\alpha)$ is viewed as a matrix-valued function. We need also to deal with a tensor-valued function $\tilde K=K^\beta_\alpha \p_\beta dx^\alpha$. Even though $\tilde K$ and $K$ has the same components, they are not exactly the same. For example we have
\[
\Delta \tilde K=\Delta K+\Gamma*\p K+\p\Gamma *K
\]
where the lower order terms with connection and curvature will appear. Note that in general, the system (with the Dirichlet boundary condition)
\[
\Delta \tilde K=0
\]
might not even be solvable due to the zeroth order terms. On the other hand, these lower order terms (the connection $\Gamma$ and its derivatives $\p \Gamma$) are assumed to be uniformly small, which would not cause essential difficulties. 
Of course the major difference is in the construction of various comparison almost complex structures. 
We present the arguments for completeness and indicate the necessary changes.
 
\begin{thm}There exists $\epsilon$ such that, for $B_r(x)\subset B$, 
\[
r^{2-m}\int_{B_r(x)}|\nabla J|^2dv<\epsilon, 
\]
then $J\in C^\infty(B_{r/4}(x))$ and for $k\geq 1$, we have the estimates
\begin{equation}\label{regularity}
r^{2k}\sup_{B_{r/8}} |\nabla^kJ|^2\leq C r^{2-m}\int_{B_r} |\nabla J|^2.
\end{equation}
\end{thm}
The key point is to prove that $J\in C^\alpha(B_{r/4}(x))$ for some $\alpha\in (0, 1)$. Since if $J$ is continuous then $J$ is smooth. While for smooth harmonic almost complex structures,  we have derived  a priori estimates \eqref{regularity} for $k=1$, and for $k\geq 2$, this follows from the standard linear elliptic theory. Hence we need only to prove the following,

\begin{thm}\label{holder}There exists $\epsilon$ such that, for $B_r(x)\subset B$, 
\[
r^{2-m}\int_{B_r(x)}|\nabla J|^2dv<\epsilon, 
\]
then $J\in C^{\alpha}(B_{r/4}(x))$.
\end{thm}

The key is to prove the following energy decay estimates on small balls. Denote \[E_r(J)=\int_{B_r} |\nabla J|^2 dv.\] 
\begin{lemma}\label{energydecay}Let $E_1(J)<\epsilon$  be sufficiently small and $\delta<\epsilon$. There exists a dimensional constant $\theta\in (0, 1)$ such that
\begin{equation}\label{decay1}
\theta^{2-m}E_\theta(J)\leq \frac{1}{2}\left(E_1(J)+\epsilon^{1/2}\delta\right),
\end{equation}
where $\delta$ is the dimensional constant in \eqref{delta}.
\end{lemma}

Given the energy decay estimates in Lemma \ref{energydecay}, the proof of H\"older estimate in Theorem \ref{holder} follows a rather standard process. We shall skip the details. For more detailed arguments, see \cite{SU}[the proof of Theorem 3.1]. Now we prove  Lemma \ref{energydecay}. 

\begin{proof}[Proof of Lemma \ref{energydecay}]The proof is similar in nature to Schoen-Uhlenbeck \cite{SU}[Proposition 3.4]. The main difference is due to various constructions of comparison almost complex structures. Our choice of comparison almost complex structures is similar to the construction in Section 3. 

In general for $J\in W^{1, 2}(\cJ_g)$, it is even not clear $J$ can be approximated by smooth almost complex structures. The monotonicity formula and the small energy $E_1(J)\leq\epsilon$ will imply that $J$ has small oscillation in average, by Poincare inequality; this would in particular imply that $J$ can be approximated by smooth almost complex structures.  Denote \[\l_{x, \rho}=\int_{B_\rho(x)} J(y) dy.\]
The standard Poincare inequality, together with the energy comparison \eqref{comparison1} implies that
\[
\int_{B_1}|J-\l_{0, 1}|^2 dx\leq C_1\int_{B_1}|\p J|^2dx\leq C_2 (E_1(J)+\delta).
\]
This implies that $J$ has small oscillation in average, if $E_1(J)$ and $\delta$ are sufficiently small.
We will also need weighted average. Let $\phi$ be a smooth radial nonnegative function with support $\text{supp}(\phi)\subset B$, such that $\int_B\phi=1$. Define
\begin{equation}\label{wave}
J_{\rho}(x):=\int_{B_\rho(x)}\phi_\rho(y-x)J(y)dy=\int_{B_1(0)}\phi(z)J(x+\rho z) dz,
\end{equation}
where we denote
\[
\phi_\rho(x)=\rho^{-m}\phi\left(\frac{x}{\rho}\right)
\]
Then we have a version of modified Poincare inequality, 
\begin{equation}\label{poincare}
\rho^{-m}\int_{B_\rho(x)}|J(y)-J_{\rho}(x)|^2dy\leq C\rho^{2-m}\int_{B_\rho(x)}|\p J|^2 dy\leq C\left(\rho^{2-m}E_{B_\rho(x)}(J)+\delta^2\rho^2\right).
\end{equation}
Now for $x\in B_{1/2}$, $\rho\in (0, 1/4]$, we have by the monotonicity formula (Theorem \ref{monotonicity}),
\begin{equation}\label{monotone}
\rho^{2-m}E_{B_\rho(x)}(J)\leq c_1 \left(2^{m-2}\int_{B_{1/2}(x)}|\nabla J|^2+\delta^2\right)\leq C (E_1(J)+\delta)
\end{equation}

When $\rho>0$, $J_{\rho}(x)$ is a smooth function of $x$. We shall emphasize that in \eqref{wave}, $J_{\rho}(x)=J(x)$ is well-defined when $\rho=0$. Moreover, for $x$ fixed, $\rho$ can be taken as a constant depending on $x$, hence $\rho$ can be a nonnegative function of $x$. This is one of the most crucial points in \cite{SU} to construct comparison maps since one needs that the comparison maps have the same boundary data on some small balls.

To proceed, we denote $\bar \rho=\epsilon^{1/4}, \tau=\epsilon^{1/8}$, and $\theta\in (\tau, 1/4)$ which will be specified below. We choose a smooth function nonincreasing $\rho(|x|)=\rho(x)$ to be defined by
\[
\rho(r)=\bar \rho,\;\text{for}\;r\leq \theta, \rho(\theta+\tau)=0\; \text{and}\; |\rho^{'}(r)|\leq 2\tau
\]
Then we set, as in \eqref{wave}, that, for $\rho(x)>0$,
\begin{equation}\label{wave1}
J_{\rho(x)}(x):=\int_{B_\rho(x)}\phi_\rho(y-x)J(y)dy=\int_{B_1(0)}\phi(z)J(x+\rho z) dz;
\end{equation}
while for $\rho(x)=0$, as noted above, $J_{\rho(x)}(x)=J(x)$. 
Again the Poincare inequality holds for $\rho=\rho(x)>0$ as in \eqref{poincare}. 
\begin{lemma}
We have the following estimates, for $\rho\in (0, \bar \rho]$ and $x\in B_{1/2}$,
\begin{equation}\label{estimate5}
\int_{B_{1/2}}|\nabla J_\rho|^2\leq CE_1(J);\; |\nabla J_{\bar \rho}|^2\leq C\epsilon^{1/2}
\end{equation}
It follows that
\begin{equation}\label{osc}
|\text{osc}_{B_{1/2}} J_{\bar \rho}|\leq C\epsilon^{1/4}.
\end{equation}
\end{lemma}

\begin{proof}The first integral estimate in \eqref{estimate5} is by a simple direct computation. The second follows from the monotonicity as follows,
\[
\begin{split}
|\nabla J_{\bar \rho}|^2=&\left(\int_B \phi(z)\nabla J(x+\bar \rho z)dz\right)^2\\
\leq& \int_B\phi(z)|\nabla J|^2(x+\bar \rho z)dz\\
=&\bar\rho^{-m}\int_{B_{\bar \rho}} \phi\left(\frac{x-y}{\bar\rho}\right)|\nabla J|^2(y) dy\\
\leq &C\bar\rho^{-2} \bar \rho^{2-m}E_{\bar \rho}(J)\\
\leq & C\bar\rho^{-2}(E_1(J)+c_n\delta^2)\\
\leq & C\epsilon^{1/2}
\end{split}
\]
This completes the proof. 
\end{proof}
Now we construct a comparison almost complex structure using $J_\rho(x)$ when $\rho(x)>0$. Clearly $J_\rho(x)$ might  not be an almost complex structure which is $g$-compatible.
But by \eqref{wave1} and \eqref{poincare}, $J_\rho(x)$ is almost an almost complex structure and is almost $g$-compatible. By \eqref{wave}, we have
 \[
 J_\rho(x)+ gJ^t_\rho(x) g^{-1}=\int_{B}\phi(z) \left(g(x)J^t(x+\rho z)g^{-1}(x)- g(x+\rho z)J^t(x+\rho z)g^{-1}(x+\rho z) \right)dz
 \]
We compute  
\[
\begin{split}
g(x)J^t(x+\rho z)g^{-1}(x)- g(x+\rho z)J^t(x+\rho z)&g^{-1}(x+\rho z)=\left(g(x)-g(x+\rho z)\right) J^t(x+\rho z) g^{-1}(x)\\
&+g(x+\rho z)J^t(x+\rho z)(g^{-1}(x)-g^{-1}(x+\rho z))
\end{split}
\]
It follows that
\begin{equation}\label{a-0}
|J_\rho(x)+g(x) J^t_\rho(x) g^{-1}(x)|\leq C\delta \rho^2. 
\end{equation}
By \eqref{poincare}, there exists some $y\in B_\rho(x)$, 
\[
|J(y)-J_\rho(x)|\leq C\sqrt{\epsilon}
\]
In particular it follows that
\begin{equation}\label{a-1}
|J_\rho(x) J_\rho(x)+id|\leq C\sqrt{\epsilon}.
\end{equation}
Let $S_g$ and $A_g$ be the $g$-symmetric part and $g$-skewsymmetric part of $J_\rho$ respectively. 
Then we have that  ($\delta<\epsilon$), using \eqref{a-0} and \eqref{a-1}, 
\begin{equation}\label{a-2}
|A^2_g+id|\leq C\sqrt{\epsilon}, \;|S_g|\leq C\delta \rho^2
\end{equation}
Since $A_g$ is $g$-skewsymmetric, $A^2_g$ is $g$-symmetric. By \eqref{a-2}, we conclude that $-A^2_g$ is a $g$-symmetric positive definite matrix, provided $\epsilon$ is small enough. Hence there exists a unique $g$-symmetric positive matrix such that $Q^2=-A^2_g$ and $QA_g=A_gQ$ (see Proposition \ref{matrix1}). Denote
\begin{equation}
\bar J_\rho(x)=Q^{-1}A_g, 
\end{equation}
Then $\bar J_\rho(x)$ is a $g$-compatible almost complex structure. Note that when $\rho=0$, or $|x|\geq \theta+\tau$, $J_\rho(x)=J(x)$ and hence we also have $\bar J_\rho (x)=J_\rho(x)=J(x)$ for  $|x|\geq \theta+\tau$.
We compute,
\begin{equation}\label{comparisonj}
|\nabla \bar J_\rho|\leq C|\nabla A_g|\leq C ((1+c_n\delta \rho) |\nabla J_\rho|+c_n\delta)
\end{equation}
We also compute
\[
|\nabla J_\rho|\leq \int_B\phi(z) (|\nabla J|+|\rho^{'}||\nabla J|)(x+\rho z)dz
\]
It then follows that
\begin{equation}\label{e1}
\begin{split}
\int_{B_{\theta+\tau}\backslash B_\theta}|\nabla \bar J_\rho|^2\leq& C\int_{B_{\theta+\tau}\backslash B_\theta}|\nabla J_\rho|^2+C\delta^2\\
\leq &C \int_{B_{\theta+\tau}\backslash B_{\theta}}\left(\int_{B}\phi(z)|\nabla J|(x+\rho z) dz\right)^2+C\delta^2\\
\leq & C\int_{B_{\theta+2\tau}\backslash B_{\theta-\tau}} |\nabla J|^2+C\delta^2
\end{split}
\end{equation}
Since $\bar J_\rho (x)$ agrees with $J(x)$ for $|x|\geq \theta+\tau$, then by the minimizing property of $J$ we have
\begin{equation}\label{e2}
E_{\theta+\tau}(J)\leq E_{\theta+\tau}(\bar J_\rho)
\end{equation}
Now we want to estimate $E_\theta (\bar J_\rho)$. By \eqref{comparisonj}, we have
\begin{equation}\label{e3}
E_\theta (\bar J_\rho)\leq CE_\theta (J_\rho)+C\delta^2
\end{equation}
To estimate $E_\theta(J_\rho)$, we need a harmonic approximation of $J_\rho$. Note that $\rho=\bar \rho$ when $|x|\leq \theta$. We solve the following harmonic equation with the boundary data $J_{\bar\rho}$,
\begin{equation}\label{laplacian} \Delta K=0, \;\text{in}\; B_{1/2};\;
K=J_{\bar \rho}, \;\text{on}\; \p B_{1/2}
\end{equation}
Here $K$ is a matrix-valued function and by \eqref{laplacian} we solve the equation for each entry of the matrix. We also denote $\tilde K=K^\beta_\alpha \p_\beta\otimes dx^\alpha$ as the corresponding tensor. 
Note that we have lower order difference for derivatives of $K$ and $\tilde K$, interpreted as a matrix-valued function and tensor-valued function,
\[
\begin{split}
\nabla \tilde K=&\nabla K+\Gamma*K\\
\Delta \tilde K=&\Delta K+\Gamma*\p K+\p \Gamma*K
\end{split}
\]
Nevertheless, the difference is uniformly small (note that if $g$ is Euclidean, the difference is zero) under our assumption. 

\begin{lemma}We have the following estimate,
\begin{equation}\label{comparison6}
E_\theta(J_{\bar\rho})\leq C\epsilon^{1/4}E_1(J)+C\theta^m E_1(J),
\end{equation}
for any $\theta\in (0, 1/4)$.
\end{lemma}
\begin{proof}
Let $(J^\beta_\alpha)$ be the matrix-valued function which corresponds to the tensor $J_{\bar\rho}=J^\beta_\alpha \p_\beta\otimes dx^\alpha$.
We will need to distinguish $(K, \tilde K)$ and $(J^\beta_\alpha, J_{\bar \rho})$. As we have mentioned, under our assumption, the differences will be under uniform control. 
For example,
\[
\nabla \tilde K=\nabla K+\Gamma* K, \nabla J_{\bar\rho}=\nabla J^\beta_\alpha+\Gamma* J^\beta_\alpha.
\]
We compute
\begin{equation}\label{a1}
E_\theta(J_{\bar\rho})\leq 2\int_{B_\theta} |\nabla J_{\bar\rho}-\nabla \tilde K|^2+\int_{B_\theta}|\nabla \tilde K|^2
\end{equation}
We can compute
\[
\Delta |\nabla K|^2=2|\nabla^2 K|^2+Ric(\nabla K, \nabla K)\geq-C\delta  |\nabla K|^2.
\]
By the Harnack inequality, we have
\[
\sup_{B_{1/4}} |\nabla K|^2\leq C\int_{B_{1/2}} |\nabla K|^2
\]
Since $K$ solves \eqref{laplacian}, we have (using \eqref{estimate5})
\[
\int_{B_{1/2}} |\nabla K|^2\leq \int_{B_{1/2}}|\nabla J^\beta_\alpha|^2\leq C(E_1(J)+\delta)
\]
It follows that (since $\delta<E_1(J)=\epsilon$)
\begin{equation}\label{a2}
\int_{B_\theta}|\nabla \tilde K|^2\leq CE_1(J) \theta^m
\end{equation}
While we have (using \eqref{osc})
\[|\text{osc}_{B_{1/2}} K |\leq |\text{osu}_{B_{1/2}}J_{\bar\rho}|\leq C\epsilon^{1/4}\]
Then we compute
\[
\int_{B_\theta} |\nabla J_{\bar\rho}-\nabla \tilde K|^2=-\int_{B_\theta}\Delta J_{\bar \rho} (J_{\bar\rho}-\tilde K)+C\delta \theta^m \leq C\epsilon^{1/4}\int_{B_{1/2}} |\Delta J_{\bar \rho}|+C\delta \theta^m
\]
We compute, by \eqref{wave1} and the fact that $J$ is a weakly harmonic almost complex structure,
\begin{equation}\label{part}
\begin{split}
\Delta J_{\bar \rho}=&\int_{\C^m} \Delta_x \phi_{\bar\rho}(x-y)J(y)dy\\
=&\int_{\C^m}\Delta_y \phi_{\bar\rho}(x-y)J(y) dy\\
=&\int_{\C^m}\phi_{\bar \rho}(J\nabla_p J\nabla_pJ+\p J*\p g+J*\p^2 g) dy,
\end{split}
\end{equation}
where the last step needs some explanation. First we write
\[
\Delta=g^{i j}\p_{x_i}\p_{x_j}+\p g*\p,
\]
where $\p g *\p$ denotes the first order terms. Integration by part, we get for any $\phi$ with compact support,
\begin{equation}\label{part1}
\int \Delta \phi J=\int \phi (\Delta J+\p J*\p g+J*\p^2 g)
\end{equation}
We use the structural equation $\Delta J=J\nabla_pJ \nabla_p J$ in the weak sense that, 
\[
\int\Delta \phi J=\int (-g^{ij}\p_i\phi \p_j J+\phi J*\p^2g+\phi \p J*\p g)=\int \phi (J\nabla J\nabla J+\p J*\p g+J*\p^2 g),
\]
which proves \eqref{part}. 
Hence it follows that
\begin{equation}\label{a3}
\int_{B_{1/2}} |\Delta J_{\bar \rho}|\leq C(E_1(J)+\delta)
\end{equation}
The claim \eqref{comparison6} follows by combining \eqref{a1}, \eqref{a2} and \eqref{a3}. 
\end{proof}
The remaining arguments are then the same as in \cite{SU}[Proposition 3.4]. We sketch the details for completeness. 
Now let $\gamma_m\in (0, 1/16]$ be a number depending only on $m$ and set $\bar\theta=\epsilon^{\gamma_m}$. Let $p$ be the greatest integer that is less than (or equal to) $\bar\theta/3\tau$ (recall that $\tau=\epsilon^{1/8}$) and write the interval
\[
[\bar\theta, \bar \theta+3p\tau]=\cup_i I_i,
\]
where $I_i=[\bar \theta+(i-1)\tau, \bar \theta+i\tau]$, $i=1, 2, \cdots, 3p$. Since 
\[
\sum_i\int_{|x|\in I_i}|\nabla J|^2=E_{\bar \theta+3p\tau}(J)-E_{\bar\theta}(J)\leq E_1(J),
\]
there exists some $j$ such that for $|x|\in I_j$, 
\[
\int_{|x|\in I_j}|\nabla J|^2\leq p^{-1}E_1\leq C\epsilon^{1/16} E_1(J)
\]
Let $\theta$ be the number such that $I_j=[\theta-\tau, \theta+2\tau]$.  We compute,
\begin{equation}
\begin{split}
E_{\theta+\tau}(J)\leq& E_{\theta+\tau}(\bar J_\rho)\\
=&E_{\theta}(\bar J_\rho)+\int_{B_{\theta+\tau}\backslash B_{\theta}} |\nabla \bar J_\rho|^2\\
\leq & C(\epsilon^{1/4}+\theta^m) E_1(J)+C\epsilon^{1/16}E_1(J)+C\delta^2. 
\end{split}
\end{equation}
where we have used \eqref{e2}, \eqref{e1} and \eqref{comparison6}. Hence we compute ($\delta<\epsilon$)
\[
\theta^{2-m}E_\theta(J)\leq C\left(\epsilon^{1/16} \theta^{2-m}+\theta^2\right)E_1(J)+C\theta^{2-m} \delta^2.
\]
Since $\theta\in [\bar \theta, 2\bar\theta]=[\epsilon^{\gamma_m}, 2\epsilon^{\gamma_m}]$, we have
\begin{equation}
\theta^{2-m}E_\theta(J)\leq C\left(\epsilon^{1/16+(2-m)\gamma_m}+\epsilon^{2\gamma_m}\right)E_1(J)+C\theta^{2-m}\delta^2.
\end{equation}
We choose $\gamma_m=\min (1/64, 1/32(m-2))$ and take $\epsilon$ sufficiently small such that $C\epsilon^{2\gamma_m} \leq 1/2$, then we arrive at the conclusion
\[
\theta^{2-m}E_\theta(J)\leq \frac{1}{2}\left(E_1(J)+\epsilon^{1/2}\delta\right).
\]
 This completes the proof of the energy decay estimate. 
\end{proof}

This $\epsilon$-regularity implies the following, 
\begin{cor}Then $p\in \cR$ ($p$ is regular) if and only if $\Theta_J(p)=0$. 
\end{cor}

It is not hard to see that for $p_j\in M$ such that $p_j\rightarrow p$, then we have
\begin{cor}$\Theta_J(p)$ is upper-semicontinuous with respect to $p$.
\[
\Theta_J(p)\geq \limsup_j\Theta_J(p_j).
\]
\end{cor}

We also have the following control of singular set, 
\begin{cor}If $J\in W^{1, 2}({\cJ_g})$ is energy-minimizing, then $\cH^{m-2}(\cS)=0$; namely the singular set $\cS$ has $m-2$ Hausdorff measure zero. 
\end{cor}
\begin{proof}
This is standard and we include the details for completeness. 
We need to show that, for any $\varepsilon>0$, there exists a countable collections of balls $B_{\rho_j}(p_j)$ such that $\cS\subset \bigcup_{j} B_{\rho_i}(p_j)$, and $\sum \rho_j^{m-2}<\varepsilon$. 
For $p\in \cS$, there exists a uniformly $\epsilon_0>0$, such that $\Theta_J(p)\geq \epsilon_0$. By the monotonicity, we then have, for any $r\in (0, 1]$, 
\[
\int_{B_r(p)}|\nabla J|^2 dv\geq \epsilon_0 r^{m-2}. 
\]
Fix $\rho\in (0, 1)$, we can choose a maximal collection of balls $B_{\rho}(p_i)$, $p_i\in \cS$, for $i=1, \cdots, I$, such that $B_\rho(p_i)\cap B_{\rho}(p_j)=\emptyset$ for all $i\neq j$, and such that $I$ is the maximum integer such that such a collection exists. 
Then $\cS\subset \bigcup B_{2\rho}(p_i)$. Otherwise, there exists $p\in \cS\backslash \bigcup B_{2\rho}(p_i)$, then we can add $B_\rho(p)$ to the collection. This contradicts the maximality of $I$. 
It then follows that, 
\begin{equation}\label{hausdorff}
I \rho^{m-2}\leq \epsilon_0^{-1} \int_{\cup B_\rho(p_i)} |\nabla J|^2 dv
\end{equation}
In particular, we have
\[
I \rho^{m}\leq \epsilon_0^{-1} \rho^2 \int_M |\nabla J|^2 dv
\]
Letting $\rho\rightarrow 0$, this implies that $\cS$ has Lebesgue measure zero. Again, by letting $\rho\rightarrow 0$, then by the Lebesgue dominated convergence theorem,
\[
\int_{\cup B_\rho(p_i)}|\nabla J|^2dv\rightarrow 0.
\]
Hence \eqref{hausdorff} implies that $\cH^{m-2}(\cS)=0$. 
\end{proof}

\subsection{An improved $\epsilon$-regularity}

We derive an improved $\epsilon$-regularity for energy-minimizing harmonic almost complex structures. For energy minimizing harmonic maps, an improved $\epsilon$ regularity is due to Schoen-Uhlenbeck \cite{SU}, and the compactness theorem is due to Luckhaus \cite{Luckhaus1, Luckhaus2} (partial results have been obtained by  Schoen-Uhlenbeck \cite{SU} and Hardt-Lin \cite{HL}).  Our treatment uses an extension lemma by Luckhaus and we follow the presentations in L. Simon's book \cite{Simon} (see Sections 2.6-2.9).

The following corollary of Luckhaus' lemma (see Lemma 1, Corollary 1 in \cite{Simon}[Section 2.6, Section 2.7]) will play an important role. 

\begin{cor}\label{extension}Suppose $N$ is a smooth compact manifold embedded in $\R^p$ and $\Lambda>0$. There exists $\delta_0=\delta_0(n, N, \Lambda)$ and $C=C(n, N, \Lambda)$ such that the following hold:
\begin{enumerate}
\item If we have $\epsilon \in (0, 1]$ and if $u\in W^{1, 2}(B_\rho(y); N)$ with $\rho^{2-n}\int_{B_\rho(y)}|D u|^2\leq \Lambda$ and
$\rho^{-n}\int_{B_\rho(y)} |u-\lambda_{x, \rho}|^2\leq \delta_0^2\epsilon^{2n}$, then there is $\sigma\in (3\rho/4, \rho)$ such that there exists a function $w=w_\epsilon \in W^{1, 2}(B_\rho(y); N)$ which agrees with $u$ in a neighborhood of $\p B_\sigma (y)$ and which satisfies
\begin{equation}\label{extension-01}
\sigma^{2-n}\int_{B_\sigma(y)} |Dw|^2\leq \epsilon \rho^{2-n}\int_{B_\rho(y)}|D u|^2+\epsilon^{-1} C\rho^{-n}\int_{B_\rho(y)} |u-\lambda_{x, \rho}|^2
\end{equation}

\item  If $\epsilon \in (0, \delta_0]$ and $\Omega=B_{(1+\epsilon)\rho}(y)\backslash B_{\rho}(y)$, and  if $u, v\in W^{1, 2}(\Omega; N)$ satisfies $\rho^{2-n}\int_{\Omega}(|D u|^2+|Dv|^2)\leq \Lambda$ and
$\rho^{-n}\int_{\Omega} |u-v|^2\leq \delta_0^2\epsilon^{2n}$, then  there exists a function $w=w_\epsilon \in W^{1, 2}(\Omega; N)$ such that  $w=u$ in a neighborhood of $\p B_\rho (y)$ and  $w=v$ in a neighborhood of $\p B_{(1+\epsilon)\rho}(y)$, and
\begin{equation}\label{extension-02}
\rho^{2-n}\int_{\Omega} |Dw|^2\leq C \rho^{2-n}\int_{\Omega}(|D u|^2+|Dv|^2)+\epsilon^{-2} C\rho^{-n}\int_{\Omega} |u-v|^2
\end{equation}
\end{enumerate}
\end{cor}

We want to use the extension lemmas above to construct comparison almost complex structures. We need some preparations. First we consider almost complex structures in $W^{1, 2}$ and are compatible with the Euclidean metric (over a unit ball $B$), which we denote as $W^{1, 2}(\cJ_{Eu})$. 
We write an almost complex structure as a matrix-valued functions (in $\R^{m\times m}$) over $B$ (the unit ball).  Denote $N$ to be the submanifold of $\R^{m\times m}$ given by the restrictions, $J\in \R^{m\times m}$, $J^2=-id, J+J^t=0, \det(J)=1$. Then $W^{1, 2}(\cJ_{Eu})=W^{1, 2}(B, N)$. For such almost complex structures $J\in W^{1, 2}(\cJ_{Eu})$, we can then apply Corollary \ref{extension}. 
Now we consider almost complex structures in $W^{1, 2}(\cJ_g)$.  Suppose $J$ is an almost complex structure in $W^{1, 2}(\cJ_{Eu})$, then by Proposition \ref{matrix1}, we can construct a unique almost complex structure $\bar J$ out of $J$, such that $\bar J$ is $g$-compatible. Then by Proposition \ref{matrix1} again, we construct a unique almost complex structure $\bar{\bar J}$, out of $\bar J$, such that $\bar{\bar J}$ is compatible with the Euclidean metric. By the construction, it would be clear that $\bar {\bar J}=J$. Indeed, if $g=GG^t$ locally, then $\bar J=G^{-1}JG$. In particular, $J$ is $g$-skewsymmetric if and only if $\bar J=G^{-1}JG$ is skew-symmetric. 
With this relation between $W^{1, 2}(\cJ_{Eu})$ and $W^{1, 2}(\cJ_g)$, we can then apply Corollary \ref{extension} to get an extension for $J\in W^{1, 2}(\cJ_g)$. This will allow us to obtain an improved version of $\epsilon$-regularity as follows.
We focus locally over a unit ball $B$.  The metric $g$ is close to the Euclidean metric, as in \eqref{delta}. 

\begin{thm}\label{improve}Given $\Lambda>0$, there exists a constant $\epsilon_0=\epsilon_0(m, \Lambda)$ such that if $E_1(J)\leq \Lambda$ and
\[
\int_B |J-\lambda_{0, 1}|^2 \leq \epsilon_0,
\]
then $J$ is H\"older continuous on $B_{3/16}$. 
\end{thm} 

\begin{proof}Given $J$, we can then construct $\bar J$ (using Proposition \ref{matrix}), which is a compatible almost complex structure with the Euclidean metric, such that
\begin{equation}\label{s1}
\int_B|\bar J-\bar\lambda_{0, 1}|^2\leq c_n( \epsilon_0+ \delta).\end{equation}
We can obtain \eqref{s1} as follows, 
\[
\int_B |\bar J-\bar\lambda_{0, 1}|^2\leq 3\int_B (|\bar J-J|^2+|J-\lambda_{0, 1}|^2)+3(\lambda_{0, 1}-\bar \lambda_{0, 1})^2,
\]
where
\[
|\lambda_{0, 1}-\bar \lambda_{0, 1}|^2=\left|\int_B (J-\bar J)\right|^2\leq c_n \int_B|J-\bar J|^2
\]
It then follows that
\[
\int_B |\bar J-\bar\lambda_{0, 1}|^2\leq c_n \int_B (|\bar J-J|^2+|J-\lambda_{0, 1}|^2)
\]
Since $g$ is close to the Euclidean metric as in \eqref{delta}, we can then estimate as usual, that
\[
\int_B |\bar J-J|^2\leq c_n \delta. 
\]
This establishes \eqref{s1}. Moreover we have $E_1(\bar J)\leq 2(\Lambda+1)$, similarly as in \eqref{comparison3}. Now we apply Luckhaus' extension lemma to $\bar J$, as in Corollary \ref{extension} to $\bar J$. For $\epsilon\in (0, 1]$, we obtain an almost complex structure $\bar K$, which is compatible with the Euclidean metric and which agrees with $\bar J$ near $\p B_\sigma$ for some $\sigma\in (3/4, 1)$, such that
\[
\sigma^{2-m}\int_{B_\sigma}|D\bar K|^2\leq \epsilon E_1(\bar J)+\epsilon^{-1} C\int_B |\bar J-\bar \lambda_{0, 1}|^2
\]
Now we need to construct $K$, out of $\bar K$ such that $K$ is an almost complex structure which is compatible with $g$, using Proposition \ref{matrix1} again. Similarly as in \eqref{comparison3}, we have
\[
\int_{B_{\sigma}}|\nabla K|^2 dv\leq (1+c_n \delta) \int_{B_\sigma}|D\bar K|^2+c_n \delta^2 \sigma^m
\]
It then follows that,
\[
\sigma^{2-m}\int_{B_{\sigma}}|\nabla K|^2 dv\leq c_n \epsilon (\Lambda+1)+C \epsilon^{-1} (\epsilon_0+\delta). 
\]
Since $\bar K$ agrees with $\bar J$ in a neighborhood of $\p B_\sigma $, then $K$ agrees with $J$ in a neighborhood of $\p B_\sigma$, by the discussion above Theorem \ref{improve}. The energy minimizing then implies that
\[
\sigma^{2-m}\int_{B_{\sigma}}|\nabla J|^2 dv\leq \sigma^{2-m}\int_{B_{\sigma}}|\nabla K|^2 dv\leq c_n \epsilon (\Lambda+1)+C \epsilon^{-1} (\epsilon_0+\delta). 
\]
If $\epsilon_0$ and $\delta$ are sufficiently small, then we can choose $\epsilon$ accordingly, such that $\sigma^{2-m}\int_{B_{\sigma}}|\nabla J|^2 dv$ is sufficiently small. By the $\epsilon$-regularity, we can then conclude that $J$ is H\"older continuous in $B_{3/16}$ ($\sigma>3/4$). 
\end{proof}

\subsection{Tangent almost complex structure and stratification of singular set}
We first need a local notion for energy-minimizing almost complex structures. Let $\Omega\subset \R^{2n}$ be an open set with a metric $g=g_{ij}dx^i\wedge dx^j$. We assume that $g$ is closed to the Euclidean metric as in \eqref{delta}. Suppose $g$ is almost Hermitian, we can then define $g$-compatible almost complex structures in $\Omega$, denoted by $W^{1, 2}(\Omega, \cJ_g)$. The energy functional in a ball $B_\rho(y)$ such that $\bar B_\rho(y)\subset \Omega$, is given by
\[
E_{B_\rho(y)}(J)=\int_{B_\rho(y)} |\nabla J|^2 dv,
\]
where $\nabla, dv$ are corresponding notions defined by $g$. We also use $D, \p$ to denote the derivatives with respect to coordinates $x=(x_1, \cdots, x_{2n})$. Suppose $g$ is $\lambda$-closed to the Euclidean metric as in \eqref{delta} (replacing $\delta$ by $\lambda$), then in particular,
\[
\int_{B_\rho(y)}|\nabla J|^2 dv=\int_{B_{\rho(y)}} |DJ|^2+o(\lambda). 
\]

\begin{defn}An almost complex structure $J\in W^{1, 2}(\Omega, \cJ_g)$ is energy-minimizing if for each ball $B_\rho(y)\subset \Omega$, 
\[
E_{B_\rho(y)}(J)\leq E_{B_\rho(y)}(K)
\]
for any $K\in W^{1, 2}(\Omega, \cJ_g)$ such that $K=J$ in a neighborhood of $\p B_\rho(y)$. 
\end{defn}

Let $J$ be an energy-minimizing almost complex structure in $W^{1, 2}(\cJ_g)$. 
For any point $p\in M$, we can define a tangent almost complex structure as follows. Choose a local coordinate around $p$, which we identify with the unit ball $B$.
Again we assume in $B$, the metric $g$ is close the Euclidean metric as in \eqref{delta}. 
For any $x\in \R^{m}=\R^{2n}$ such that $\l_i x\in B$, define
\[
J_i(x)=J(x \l_i),
\]
where $\l_i$ is a sequence of positive numbers with limit zero. We will also need the notation $g_i(x)=g(\lambda_i x)$. When $\lambda_i\rightarrow 0$, $g_i(x)$ converges uniformly on compact subsets to the Euclidean metric on $\R^{2n}$. 
A point to make is that $J_i(x)$ is energy-minimizing locally for $B_R\subset \R^{2n}$, when $i$ is sufficiently large such that $\lambda_i R<1$. 
If $p$ is a smooth point, it is evident that $J_i(x)$ converges (uniformly on compact subsets) to a constant almost complex structure $J_p=J(0)$. When $p$ is a singular point, we have the following,

\begin{thm}\label{tangent}By passing to a subsequence if necessary, $J_i(x)$ converges to an almost complex structure $T_0$ on $\R^{2n}$, which is compatible with the Euclidean metric. Moreover, the convergence is in the strong topology in $W^{1, 2}$ (on compact subsets) and the limit $T_0$ is an energy-minimizing almost complex structure on $\R^{2n}$.   
\end{thm}

\begin{proof}For each fixed $R>0$, we consider $J_i(x)$ defined in $B_R$, centered at a fixed point. We compute, denoting $y=\lambda_i x$, 
\[
R^{2-m} \int_{B_R} |D J_i|^2 dx=(\lambda_i R)^{2-m}\int_{B_{\lambda_i R}} |D J|^2 dy
\]
The monotonicity formula then implies that there exists a uniform constant $\Lambda$, such that, 
\[
R^{2-m} \int_{B_R} |D J_i|^2 dx\leq \Lambda
\]
Hence $J_i(x)$ converges (by subsequence) weakly to $J_0$ in $W^{1, 2}$ on compact subsets of $\R^{2n}$ and the convergence is strong in $L^2$. Since $J_i(x)$ is an almost complex structure compatible with $g_i(x)$, then $T_0$  satisfies the constraints, $T_0^2=-id$, $T_0^t+T_0=0$ a. e. In other words, $T_0\in W^{1, 2}(\cJ_{Euc})$. 
Next we want to show that $J_i$ converges by subsequence to $T_0$ strongly and $T_0$ is also energy-minimizing. Our argument follows closely \cite{Simon}[Section 2.9, Lemma 1]. 

Fix $\theta\in (0, 1)$ and $\delta_1>0$. Choose an integer $N\in \N$ such that
\[
R^{2-m}\int_{B_R} |DJ_i|^2< N\delta_1. 
\]
Let $\varepsilon\in (0, (1-\theta)N^{-1})$. There must be an integer $2\leq l\leq N$ such that
\[
R^{2-m}\int_{B_{R(\theta+l\varepsilon)}\backslash B_{R(\theta+(l-2)\varepsilon)}} |DJ_i|^2<\delta_1
\]
holds for infinitely many $i$. Choose such $l$ and let $r=R(\theta+(l-2)\varepsilon)$, then $r(1+\varepsilon)<R(\theta+l\varepsilon)<R$. There exists a subsequence, which we still denote by $J_i$,  such that
\[
R^{2-m}\int_{B_{r(1+\varepsilon)}\backslash B_{r}} |DJ_i|^2<\delta_1
\]
We can assume $J_i$ converges to $T_0$ weakly in $W^{1, 2}$ and strongly in $L^2$, by passing to a subsequence if necessary. We also have
\[
R^{2-m}\int_{B_{r(1+\varepsilon)}\backslash B_{r}} |DT_0|^2\leq \liminf_{i\rightarrow \infty} R^{2-m}\int_{B_{r(1+\varepsilon)}\backslash B_{r}} |DJ_i|^2\leq  \delta_1. 
\]
Now we apply first Proposition \ref{matrix1} to construct $\bar J_i$, out of $J_i$, such that $\bar J_i$ is compatible with the Euclidean metric. In particular, we have the following estimates, argued similarly as in Proposition \ref{c3},
\[
\int_{B_{R}}|J_i-\bar J_i|^2=o(\lambda_i),\;\; R^{2-m}\int_{B_{r(1+\varepsilon)}\backslash B_{r}} |D\bar J_i|^2=R^{2-m}\int_{B_{r(1+\varepsilon)}\backslash B_{r}} |D J_i|^2+o(\lambda_i)\]
We can then apply Corollary \ref{extension} (2) to $\bar J_i$ and $T_0$. 
We can get an almost complex structure $K_i$, which agrees with $\bar J_i$ in a neighborhood of $\p B_{r(1+\varepsilon)}$, and which agrees with $T_0$ in a neighborhood of $\p B_{r}$, and 
\[
 R^{2-m}\int_{B_{r(1+\varepsilon)}\backslash B_{r}} |D\bar K_i|^2\leq C  R^{2-m}\int_{B_{r(1+\varepsilon)}\backslash B_{r}}\left(|D\bar J_i|^2+|DT_0|^2+\varepsilon^{-2}R^{-2}|\bar J_i-T_0|^2
\right)\]

Now suppose $K$ is an arbitrary almost complex structure compatible with the Euclidean metric, defined in $B_{\theta R}$, and  it agrees with $T_0$ in a neighborhood of $\p B_{\theta R}$. We extend $K$ over $B_{R}$ by taking $K=T_0$ on $B_{R}\backslash B_{\theta R}$. 
Define $\tilde K_i$ as follows, $\tilde K_i=K$ over $B_r$, $\tilde K_i=\bar K_i$ over $B_{r(1+\varepsilon)}\backslash B_{r}$, and $\tilde K_i=\bar J_i$ over $B_R\backslash B_{r(1+\varepsilon)}$. 

Next we want to use the fact that $J_i$ is energy-minimizing (with fixed boundary data on $\p B_{r(1+\varepsilon)}$). However $\bar K_i$ is compatible with Euclidean metric, instead of $g_i$. We need to apply Proposition \ref{matrix1} again, to construct a unique almost complex structure, out of $\tilde K_i$, which we denote as $\tilde J_i$. Again we have the estimate
\[
\int_{B_{r(1+\varepsilon)}} |D\tilde J_i|^2=\int_{B_{r(1+\varepsilon)}}|D \tilde K_i^2|+o(\lambda_i). 
\]
It is important to note that $\tilde K_i$ agrees with $\bar J_i$ near $\p B_{r(1+\varepsilon)}$, hence $\tilde J_i$ agree with $J_i$ near $\p B_{r(1+\varepsilon)}$. Hence by minimizing property of $J_i$, we have
\[
\int_{B_{r(1+\varepsilon)}} |\nabla_i J_i|^2 dv_i\leq \int_{B_{r(1+\varepsilon)}} |\nabla_i \tilde J_i|^2dv_i, 
\]
where $\nabla_i, dv_i$ are geometric quantities associated with the metric $g_i(x)$. Now we need the estimates
\[
\int_{B_{r(1+\varepsilon)}} |\nabla_i J_i|^2 dv_i= \int_{B_{r(1+\varepsilon)}} |D J_i|^2+o(\lambda_i) , \; \;\int_{B_{r(1+\varepsilon)}} |\nabla_i \tilde J_i|^2dv_i=\int_{B_{r(1+\varepsilon)}} |D\tilde J_i|^2+o(\lambda_i). 
\]
Hence, we have the following,
\begin{equation}\label{c4}
\begin{split}
r^{2-m}\int_{B_r}|DJ_0|^2\leq& \liminf r^{2-m}\int_{B_r}|DJ_i|^2\\
\leq &\liminf r^{2-m}\int_{B_{r(1+\epsilon)}} |D\tilde J_i|^2\\
\leq &\liminf r^{2-m}\int_{B_{r(1+\epsilon)}} |D\tilde K_i|^2\\
\leq &r^{2-m}\int_{B_r}|D K|^2+\liminf r^{2-m}\int_{B_{r(1+\varepsilon)}\backslash B_{r}} |D\bar K_i|^2\\
\leq &r^{2-m}\left(\int_{B_r}|D K|^2+C\delta_1\right)
\end{split}
\end{equation}
Since $\delta_1$ is arbitrary, and $\theta R<r$, we have the following,
\[
\int_{B_{\theta R}}|DJ_0|^2\leq \int_{B_{\theta R}} |D K|^2. 
\]
In other words, $T_0$ is energy minimizing on $B_{\theta R}$. Since $\theta\in (0, 1)$ and $R$ are also arbitrary, this proves that $T_0$ is energy-minimizing. 

To prove the strong convergence in $W^{1, 2}$, we can take $K=T_0$ in \eqref{c4}, then we obtain, 
\[
\liminf \int_{B_r} |DJ_i|^2\leq \int_{B_r} |D T_0|^2. 
\]
Hence we have
\[
\liminf \int_{B_r} |DJ_i|^2= \int_{B_r} |D T_0|^2. 
\]
By passing to a subsequence if necessary, we can assume $J_i$ converges to $T_0$ weakly and
\[
\lim  \int_{B_r} |DJ_i|^2= \int_{B_r} |D T_0|^2. 
\]
This would imply particularly the strong convergence of $J_i$ to $T_0$ in $W^{1, 2}$. 
\end{proof}

There are nice results about  tangent maps for energy-minimizing harmonic maps and the stratification of singular sets using the properties of tangent maps , see \cite{SU, White97} . These results can be viewed as  adaption of corresponding results for minimal surfaces of F. Almgren \cite{Ag83}.  Given the results we have derived for energy-minimizing almost complex structures, these corresponding results for tangent maps of energy-minimizing harmonic maps can be directly carried over to tangent almost harmonic complex structures. For completeness we include some details. Our presentation follows \cite{Simon}[Section 3].  

Suppose $J$ is an energy-minimizing almost complex structure and $p\in \cS$ is a singular point. 
Denote $T$ to denote a tangent almost complex structure of $J$ at $p\in M$, constructed by Theorem \ref{tangent}. 
Then we have the following,

\begin{prop}\label{tanmon}The energy density satisfies \[\Theta_J(p)=\Theta_T(0)=r^{2-m}\int_{B_r} |DT|^2, \forall r>0.\] Moreover, we have the formula
\begin{equation}\label{tanmon1}
0=R^{2-m}\int_{B_R}|DT|^2-r^{2-m}\int_{B_r}|DT|^2=\int_{B_R\backslash B_r} |\p_r T|^2
\end{equation}
Hence $T(\lambda x)=T(x)$ for all $\lambda>0$, $x\in \R^{2n}$, and $\p_r T=0$. 
\end{prop}

\begin{proof}By definition, there exists $\l_i\rightarrow 0$, such that $J_i(x)=J(\l_i x)$ converges to $T$ strongly in $W^{1, 2}$ on compact subsets of $\R^{2n}$. Write $g_i(x)=g(\l_i x)$.
Then 
\[
\Theta_J(p)=\lim_{r\rightarrow 0} r^{2-m}\int_{B_r(p)}|\nabla J|^2dv=\lim_{i\rightarrow \infty}\int_{B_1} |\nabla_i J_i|^2 dv_i=\int_{B_1}|DT|^2. 
\]
A similar argument gives, for $\sigma>0$ fixed
\[
\Theta_J(p)=\lim_{r\rightarrow 0} (\sigma r)^{2-m}\int_{B_{\sigma r}(p)}|\nabla J|^2dv=\lim_{i\rightarrow \infty}\sigma^{2-m}\int_{B_\sigma} |\nabla_i J_i|^2 dv_i=\sigma^{2-m}\int_{B_\sigma}|DT|^2.
\]
This completes the proof. 
\end{proof}

Proposition \ref{tanmon} gives a nice description of a regular point of $J$: $p\in M$ is regular if and only if there exists a constant tangent almost complex structure at $p$, since $T$ is constant if and only if  $\Theta_T(0)=\Theta_J(p)=0$.
We can also have the following version of compactness result, which allows we consider a sequence of energy-minimizers (even at different points). Suppose $J_i$ is a sequence of energy-minimizing almost complex structures in $W^{1, 2}(\cJ_g)$. Fix a sequence $p_i\in M$ and identify $B_1(p_i)$ with $B_1\subset \R^m$ via the exponential map $\exp_{p_i}$. Let $r_i$ be a sequence such that $r_i\rightarrow 0$, define 
\[\tilde J_i(x):=J_i(\exp_{p_i}(r_ix)), x\in B_{r_i^{-1}}(0)\subset \R^m.\]
By the identical proof as in Theorem \ref{tangent}, we have the following
\begin{thm}\label{com1}By passing to a subsequence if necessary, $\tilde J_i$ converges to an admissible almost complex structure $T_0$ on $\R^{2n}$ strongly in $W^{1, 2}$ (on compact subsets). Moreover $T_0$ is also energy-minimizing.
\end{thm}

\begin{rmk}
We can even assume the background metrics $g$ varies, once we have uniform $C^3$ control on the metrics. Since we do not need these results, we skip the details.  We should emphasize a tangent almost complex structure satisfies the homogeneity condition $\p_r T=0$. But in Theorem \ref{com1}, $\p_rT_0=0$ is not true if we allow the base points vary. 
\end{rmk}

\begin{rmk}
We have shown that energy-minimizing almost complex structures in $W^{1, 2}(\cJ_g)$ are compact. With a little more effort, we can also prove a compactness theorem for local energy-minimizing almost complex structures, based on Luckhaus' extension theorem (see Corollary \ref{extension}). Such a nice compactness result for energy-minimizing harmonic maps is due to Luckhaus \cite{Luckhaus1, Luckhaus2}. 
\end{rmk}

A compatible  almost complex structure on $\R^{2n}$ (or bounded domain $\Omega$) (which induces a fixed orientation of $\R^{2n}$) is a map from $\R^{2n}$ to $SO(2n)/U(n)$, endowed with the natural metric 
\[(J_1, J_2)=-\frac{1}{2n}\tr(J_1J_2).\]
An energy-minimizing almost complex structure on $\R^{2n}$ is then an energy-minimizing harmonic map, see such an correspondence in \cite{Davidov} for example. Theorem \ref{tangent}  implies that a tangent almost complex structure of an energy-minimizing almost complex structure is an energy-minimizing, homogeneous (degree zero) almost complex structure on $\R^{2n}$, hence it is also an energy-minimizing harmonic map in $W^{1, 2}(\R^{2n}, SO(2n)/U(n))$ of homogeneous degree zero (even though an energy-minimizing almost complex structure might not be a harmonic map itself). We consider an example when the dimension is four ($n=2$).

\begin{exam}When $n=2$, $SO(4)/U(2)$ can be identified with $S^2$ (with the standard metric, up to a scaling factor). A tangent almost complex structure can be identified with an energy-minimizing, homogeneous harmonic map (of degree zero) in $W^{1, 2}_{loc}(\R^4, S^2)$. We can compute this in coordinates directly. A compatible almost complex structure satisfies $J^2=-id, J^t+J=0$. It is easy to compute that $J$ is of the form, 
\[
J=\begin{pmatrix} 0&a &b& c\\
-a&0&-c&b\\
-b&c&0&-a\\
-c&-b&a&0
\end{pmatrix} \;\text{or,}\;\; 
J=\begin{pmatrix} 0&a &b& c\\
-a&0&c&-b\\
-b&-c&0&a\\
c&b&-a&0
\end{pmatrix}
\]
with the constraint $a^2+b^2+c^2=1$.
The energy functional is then given by $E(J)=\int |\nabla J|^2=4\int|\nabla u|^2$, if we denote $u=(a, b, c)$. Hence an energy-minimizing almost complex structure on $\R^4$ is equivalent to an energy-minimizing harmonic map in $W^{1, 2}_{loc}(\R^4, S^2)$. Moreover, we know such an energy-minimizing map must be homogeneous (of degree zero), meaning $T(\lambda x)=T(x)$ for $\lambda>0$. Such harmonic maps were studied extensively, in particular see Hardt-Lin \cite{HL2}. 
\end{exam}

Now we consider the natural stratification of the singular set $\cS$. 
Suppose $T$ is a tangent almost complex structure of $J$ at $p\in \cS$. Denote
\[
S(T)=\{z\in \R^{2n}: \Theta_T(z)=\Theta_T(0)\}. 
\]
Then we have the following,
\begin{prop}\label{tan1}For any $z\in \R^{2n}$, $\Theta_T(z)\leq \Theta_T(0)$.  If the equality holds, then $T(z+\lambda x)=T(z+x)$ for any $\lambda>0, x\in \R^{2n}$. Moreover $S(T)$ is a linear subspace of $\R^{2n}$ and for any $z\in S(T)$, $x\in \R^{2n}$, $\lambda>0$,
\[
T(\lambda x+z)=T(x+z)=T(x).
\]
We also have $\text{dim}(S(T))\leq 2n-1$ and $S(T)\subset \text{Sing}(T)$. 
\end{prop}
\begin{proof}The proof in the setting of harmonic maps can be directly carried over. See Section 3.3 in \cite{Simon}. 
\end{proof}

This leads the following definition,
\begin{defn}\label{k-homogeneous}An admissible almost complex structure on $\R^{2n}$  is $k-$homogeneous at $z\in \R^{2n}$ with respect to the $k$-plane $V^k\subset \R^{2n}$ if
\begin{enumerate}
\item $T(z+\lambda x)=T(z+x)$ for any $\lambda>0, x\in \R^{2n}$.
\item $T(x+v)=T(x)$, for any $x\in \R^{2n}$, $v\in V^k$.  
\end{enumerate}
\end{defn}
If $z=0\in R^{2n}$, then we say $T$ is $k$-homogeneous.  By \emph{admissible almost complex structure} we mean the almost complex structure is in $W^{1, 2}$ and it is compatible with the Euclidean metric almost everywhere.  
We denote $\text{dim}(S(T))=k$, for some $k\in \{0, \cdots, 2n-1\}$. Proposition \ref{tan1} then means that $T$ is $k$-homogeneous at any $z\in S(T)$ with respect to $S(T)$. 
We define, for $j=0, 1, \cdots, 2n-1$,
\[
\cS^j=\{p\in \cS: \text{dim}\;S(T)\leq j\; \text{for all tangent almost harmonic complex structures}\; T\; \text{at}\; p\}
\]
In other words, $p\in \cS^j$ if and only if no tangent almost complex structure at $p$ is $(j+1)$-homogeneous. 
We can then get a stratification of the singularities, 
\[
\cS^0\subset \cS^1\subset \cdots \subset \cS^{2n-3}\subset \cS^{2n-2}\subset \cS^{2n-1}=\cS. 
\]
We have the following, 
\begin{thm}We have $\cS=\cS^{2n-3}$. Moreover, for each $j=0, \cdots, 2n-3$, the Hausdorff dimension $\text{dim}\;\cS^j\leq j$. 
\end{thm}
\begin{proof}The proof follows the case for energy-minimizing harmonic maps; in particular we follow \cite{Simon}[Section 3.4]. For completeness, we sketch a proof briefly.

We first show  $\cS=\cS^{2n-3}$. By definition, if a tangent map $T$ at $p$ satisfies $\text{dim}S(T)=2n$, then $T$ will be a constant map with $\Theta_T(0)=0$. This would imply that $p$ is a regular point. Hence $\cS=\cS^{2n-1}$.  Since a tangent almost complex structure is energy-minimizing, this implies that $\text{dim}S(T)\leq 2n-3$ since $S(T)\subset \text{Sing}(T)$ and the singular set of $T$ has Hausdorff dimension $\cH^{2n-2}=0$. This implies that $\cS^{2n-1}=\cS^{2n-2}=\cS^{2n-3}$.  

Next we want to show that the Hausdorff dimension of $\cS^j\leq j$, for $0\leq j\leq 2n-3$. 
We cover the manifold $M$ by finitely many of geodesic balls $B_{p_i}$ of radius $1$. Note that we assume in these balls, the metric $g$ satisfies the bounds \eqref{delta}; in particular the exponential maps are one-to-one and onto (on unit balls of tangent space $T_{p_i}(M)$ to the geodesic balls $B_{p_i}$) and the exponential maps are isometry on each ball, given by
\[\exp_{p_i}: (B_1, \exp^*g)\rightarrow B_{p_i}.\] Hence we need to show that the Hausdorff dimension of $\cS^j\cap B_{p_i}$ is at most $m-3$ for each $i$. From now on we fix $i$ and consider $\cS^j\cap (B_{p_i}, g)$.

For each $p_i\in M$, we identify the geodesic ball $B_{p_i}$ with the ball $B_1\subset T_{p_i}M$ through the exponential map. For any $p\in B_1$, we denote
\[
\eta_{p, \rho}(x)=\rho^{-1}(x-p), \forall x\in B_\rho(p).
\]  
We claim that for each $p\in \cS^j\cap B_{p_i}$, and each $\delta>0$, there exists an $\varepsilon_0=\varepsilon_0(J, p, \delta)$ such that for each $\rho\in (0, \varepsilon_0]$,
\begin{equation}\label{claim1}
\eta_{p, \rho}\{x\in B_\rho(p): \Theta_J(x)\geq \Theta_J(p)-\varepsilon_0\}\subset \text{the $\delta$-neighborhood of $L_{p, \rho}$} \cap B_1
\end{equation}
for some $j$-dimensional subspace $L_{p, \rho}$ of $\R^{m}=T_pM$, where $B_1$ is the unit ball in $\R^{m}=T_pM$. 

Suppose this is false, then there exists $\delta_0>0$ and $p\in \cS^j$ and a sequence of $\rho_k\rightarrow 0$, $\varepsilon_k\rightarrow 0$ such that
\[
\{x\in B_1: \Theta_{J_{\rho_k}}(x)\geq \Theta_J(y)-\varepsilon_k\}
\]
is not contained in the $\delta_0$ neighborhood of $L$ for every $j$-dimensional subspace $L$ of $\R^m$, where we denote
\[
J_k(x):=J_{\rho_k}(x)=J\circ \exp_p(\rho_k x)
\]
By passing to a subsequence if necessary, $J_k(x)$ converges to a tangent almost complex structure $T$ at $p$, such that $\Theta_J(p)=\Theta_T(0)$. Since $p\in \cS^j$, we have $\text{dim}S(T)\leq j$ by definition. Hence there exists a $j$-dimensional subspace $L_0$ such that $S(T)\subset L_0$ ($L_0=S(T)$ if $\text{dim}S(T)= j$) and a positive number $\alpha$ (depending on $\delta_0, T$ in particular), such that
\begin{equation}\label{density1}
\Theta_T(x)<\Theta_T(0)-\alpha
\end{equation}
for all $x\in \bar B_1$ with $\text{dist}(x, L_0)\geq \delta_0$. Thus, for all $k$ sufficiently large, we have
\[
\{x\in B_1: \Theta_{J_{\rho_k}}(x)\geq \Theta_T(0)-\alpha\}\subset \{x: \text{dist}(x, L_0)<\delta_0\}
\]
Otherwise we would have a subsequence $\tilde k$ with $\Theta_{J_{\rho_{\tilde k}}(x_{\tilde k})}\geq \Theta_T(0)-\alpha$ for some $x_{\tilde k}\in B_1$ with $\text{dist}(x_{\tilde k}, L_0)\geq \delta_0$. 
By taking another subsequence if necessary and using the upper semicontinuity that ($x_\infty=\lim x_{\tilde k}$), we get
\[
\Theta_T(x_\infty)\geq \limsup \Theta_{J_{\rho_{\tilde k}}(x_{\tilde k})}\geq \Theta_T(0)-\alpha,
\]
and $\text{dist}(x_\infty, L_0)\geq \delta_0$. This contradicts \eqref{density1}. 

Now we denote $\cS_{j, k}$ to be the set of points in $\cS^j\cap B_{p_i}$ such that \eqref{claim1} holds with $\varepsilon_0=k^{-1}$, for each $k\in \N$. By the claim \eqref{claim1}, we know that
\[
\cS^j\cap B_{p_i}=\cup_k \cS_{j, k}
\]
Now we denote, for each $q\geq 1$, 
\[
\cS_{j, k, q}=\{x\in \cS_{j, k}: \Theta_J(x)\in (\frac{q-1}{k}, \frac{q}{k}]\}
\]
Hence we have
\[
\cS^j\cap B_{p_i}=\cup_{k, q}\cS_{j, k, q}
\]
For each $p\in \cS_{j, k, q}$, then 
\[
\cS_{j, k, q}\subset \{x: \Theta_J(x)\geq \Theta_J(y)-\frac{1}{k}\}.
\]
Hence by the claim \eqref{claim1} with $\varepsilon_0=k^{-1}$, $\rho\leq k^{-1}$, 
\[
\eta_{p, \rho}\left(\cS_{j, k, q}\cap B_{\rho}(p)\right)
\]is contained in the $\delta$-neighborhood of $L_{p, \rho}$ for some $j$-dimensional subspace of $\R^m$. Thus each set $\cS_{j, k, q}$ has the ``$\delta$-approximation property" as in \cite{Simon}[Section 3.4]. Then by a careful covering argument (see \cite{Simon}[Section 3.4, Lemma]), it follows that $\cS_{j, k, q}$ has Hausdorff dimension at most $j$. This completes the proof. 
\end{proof}

\section{Quantitative stratification and the regularity scale}
In this section we apply the method used by Cheeger-Naber \cite{CN, CN1}, called \emph{quantitative stratification}, to study the regularity of energy-minimizing almost complex structures. Instead of focusing on the tangent almost complex structures (or tangent map, tangent cone), the method of quantitative stratification focuses on the almost tangent complex structure in an effective way. The techniques in \cite{CN1} can be applied in a straightforward way to our setting to obtain quantitative results. Our arguments follow closely \cite{CN1} for harmonic maps and we include details for completeness.  Moreover,  by combining the regularity results in Section 4 and quantitative stratification, we can improve regularities results and give an effective estimate on the \emph{regularity scale}  for energy-minimizing harmonic almost complex structures, similar as in \cite{CN1}. Such an estimate will imply the Minkowski dimension bound for the singular set. As an application, we will get uniform $L^p$ estimate of $|\nabla J|$ for any $p<3$. We should also mention the estimates of Minkowski dimension (the volume estimate of $r$-tubular neighborhood) of singular sets were developed by Chen-Donaldson \cite{CD1, CD2} in the context of K\"ahler-Einstein metrics.
\subsection{Quantitative stratification for energy-minimizing almost complex structures}
In this section we consider the quantitative stratification for energy-minimizing almost complex structures. We should mention that in \cite{CN1} for the case of harmonic maps, the discussion hold for more general stationary harmonic maps. In our setting, however, there seems to be no good notion of stationary harmonic almost complex structures, see discussions in Appendix. Hence we focus on minimizing case.

Consider a geodesic ball $B_2(p)\subset M$, we identify $B_2(p)$ with an Euclidean ball $B_2(0)$ via the exponential map, $\exp: B_2(0)\rightarrow M$ with $\exp(0)=p$. Again we assume that the injectivity radius of $(M, g)$ is larger than $2$, by a fixed scaling if necessary. We also assume that the metric $g$ is universally close to an Euclidean metric in $B_2(0)$, given by \eqref{delta}. 

We consider the almost complex structure defined on $B_{2r^{-1}}(0)\subset \R^{2n}$ by
\[
T_{p, r} J(z)=J_{p, r}(z):=J\circ \exp_p(rz)
\]
\begin{rmk}We have also used the notation $T_{p, r}(z)=J(rz)$ in Section 4, when we choose a coordinate of $B_{2}(p)$ via $\exp_p$ (hence $p=0$). We may use both notations, in particular when $p$ is fixed. 
\end{rmk}

By the discussion in Section 4, for any sequence $r_i\rightarrow 0$, $T_{p, r_i}J(z)$ converges by subsequence to a tangent almost complex structure $T$ on $\R^{2n}$ which is compatible with the Euclidean metric. Moreover $T$ is energy minimizing and homogeneous of degree zero and  this gives a stratification of the singular set $\cS$ based on homogeneity of its tangent map.
For a quantitative version of such stratification, we consider almost homogeneity of almost tangent almost complex structures of $J$. 

\begin{defn}An admissible almost complex structure $J$ is $(\varepsilon, r, k)$-homogeneous ($r\leq 1$) at $p$, if there exists a $k$-homogeneous admissible almost complex structure $T$ on $\R^{2n}$ such that \[
\int_{B_1(0)} |T_{p, r}J-T|^2\leq \varepsilon.
\]
\end{defn}

If $T$ is $k$-homogeneous with respect to the $k$-plane $V^k\subset\R^{2n}$ (see Definition \ref{k-homogeneous}), then we write 
\[
V^k_{J, p}:=B_r(p)\cap \exp_p(V^k)
\]
and it is called the defining plane of $T$. Now we introduce the quantitative singular sets.

\begin{defn}For each $\eta>0$ and $r\in (0, 1)$, the $k^{th}$ effective singular stratum $\cS^k_{\eta, r}(J)\subset M$ is the set
\[
\cS^k_{\eta, r}(J):=\left\{p\in M: \int_{B_1(0)}|T_{p, s}J-T|^2>\eta \right\}
\]
where the inequality holds for all $s\in [r, 1]$ and all $(k+1)$-homogeneous admissible almost complex structures $T$ on $\R^{2n}$. 
\end{defn}
In most cases we will simply write $\cS^k_{\eta, r}$ when $J$ is clear and fixed. 
By definition, $p\in \cS^k_{\eta, r}$ if and only if $J$ is not $(\eta, s, k+1)$-homogeneous for every $s\in [r, 1]$. It also follows immediately from the definition that
\[
\cS^{k}_{\eta, r}\subset \cS^{k^{'}}_{\eta^{'}, r^{'}}\; \;\text{if}\;\; k\leq k^{'}, r\leq r^{'}
, \eta\geq \eta^{'} \]
Moreover, we have the following description of the singular stratum,
\[
\cS^k_\eta=\bigcap_r \cS^k_{\eta, r}
\]
and
\[
\cS^k=\bigcup_{\eta} \bigcap_r \cS^k_{\eta, r}
\]

Now we can state the main result in this section, 

\begin{thm}\label{quantitative}Let $J\in W^{1, 2}(\cJ_g)$ be an energy minimizing (weakly harmonic) almost complex structure, then for any $\eta>0$ there exists $C=C(M, g, \eta)$ such that for any $r<1$, and for any $p\in M$, 
\begin{equation}\text{Vol}(\cT_r(\cS^k_{\eta, r})\cap B_1(p))\leq C r^{m-k-\eta},
 \end{equation}
 where $\cT_r(A)$ is the $r$-tubular neighborhood of a set $A$. 
\end{thm}

Note that Theorem \ref{quantitative} and its proof are purely local, hence given a point $p\in M$, we can focus on the ball $B_2(0)\subset \R^{2n}$, which is identified with $B_2(p)$ through the exponential map.  The proof follows the general strategy as in Cheeger-Naber \cite{CN, CN1}. In particular the authors proved in \cite{CN1} the volume estimate holds for stationary harmonic maps . For the stationary harmonic maps, the volume estimate obtained in \cite{CN1} is of the order $r^{m-k-\eta}$. Later on Naber-Valtorta \cite{NV} strengthen the results to the order $r^{m-k}$ for stationary harmonic map and they also proved that the $\cS^k_\eta$ is $k$-rectifiable (among other results). All these arguments are based on quantitative technique and are rather complicated. We will defer the proof of Theorem \ref{quantitative} and discuss its application on the estimate of \emph{regularity scale}. 

\subsection{Regularity scale and improved estimates}
Given a point $p$ and an almost complex structure $J\in W^{1, 2}(\cJ_g)$, we define the regularity scale $r^J(p)$ to be the maximum of $r\in (0, 1]$ such that $J$ has uniform $C^2$ bound (in terms of $r^{-1}$) over $B_r(p)$. If such a $r>0$ does not exists, $r^J(p)=0$. More precisely, we define
\[
r^J(p):=\max\{r\in [0, 1]: \sup_{B_r(p)} r|\nabla J|+r^2 |\nabla^2 J|\leq 1\}
\]
For an energy minimizer $J$, then $r^J(p)=0$ if and only if $p\in \cS$. We need an effective estimate on the set where $r^J$ is small. We denote, for $r\in (0, 1)$,
\[
\cB_r(J):=\{p\in M, r^J(p)\leq r\}.
\]
\begin{thm}\label{volumeestimate}Let $J$ be an energy-minimizing almost complex structure in $W^{1, 2}(\cJ_g)$. For all $\eta>0$, there exists a constant $C=C(\eta, M, g)$ such that for any $r\in (0, 1)$, we have
\begin{equation}\label{pestimate}
\text{Vol}(\cT_r(\cB_r(J))\cap B_1(p))\leq Cr^{3-\eta}.
\end{equation}
In particular, we get the Minkowski dimension bound
\[\text{dim}_{\text{Min}}\cS\leq 2n-3.\]
\end{thm}
 
As a direct corollary, we have the following $L^p$ bound on  $(r^J)^{-1}$ for $p<3$. 

\begin{cor}Let $J\in W^{1, 2}(\cJ_g)$ be an energy minimizer. Then for any $p\in [1, 3)$, there exists a constant $C=C(M, g, p)$ such that
\[
\int_M |\nabla J|^2 dv\leq \int_M (r^J)^{-p}dv \leq C. 
\]
Moreover, for any $p\in [2, 3),$ we have the following
\[
\int_M |\nabla^2 J|^{p/2} dv\leq \int_M (r^J)^{-p}dv \leq C.
\]
\end{cor}

\begin{proof}By covering $M$ with unit geodesic balls, the proof is local. We divide $B_1$ into the subsets $A_i$ such that $2^{-i-1}\leq r^J \leq 2^{-i}$. For any $1\leq p<3$ fixed, we can choose $\eta>0$ such that $p< 3-\eta$. With such a fixed $\eta$, we apply \eqref{pestimate} to obtain,
\[
\int_{B_1} (r^J)^{-p} dv\leq C \sum_i 2^{(i+1)p} \text{Vol}(A_i)\leq C \sum_i 2^p 2^{i (p-3+\eta)}<\infty.
\]
\end{proof}

To prove Theorem \ref{volumeestimate}, we need a new version of $\epsilon$-regularity, based on the homogeneity of the energy-minimizing almost complex structure. 

\begin{thm}\label{improve1}Let $J$ be an energy-minimizing almost complex structure in $W^{1, 2}(\cJ_g)$. For any $p\in M$, there exist $\epsilon=\epsilon(M, g)$ and $r_0=r_0(M, g)$ such that if $J$ is $(\epsilon, r, m-2)$-homogeneous over $B_{r}(p)$, then $r^J(p)\geq r/8$, for all $r\in (0, r_0]$. 
\end{thm}

The proof of Theorem \ref{improve1} relies on the improved regularity Theorem \ref{improve} and a quantitative description of a $(m-2)$-homogeneous almost complex structure. First we state this general result for homogeneous objects, see \cite{CN1}[Lemma 2.5].  

\begin{exam}
A prototype example is as follows, the map $(x, y)\rightarrow (x/\sqrt{x^2+y^2}, y/\sqrt{x^2+y^2})$ is a homogeneous map of degree zero from $\R^2\rightarrow S^1$. A straightforward computation shows that the map has infinite energy over $B_1(0)$. Indeed any homogeneous map of degree zero on $B_1(0)\subset \R^2$ has infinite energy, unless it is a constant map. Similarly, a $(m-2)$-homogeneous map with finite energy has to be $m$-homogeneous, namely a constant map. Lemma 2.5 in \cite{CN1} is an effective version. Since this is a general phenomenon for homogeneous objects,  it also holds for our setting. 
\end{exam}

We need the following version in our setting, for an admissible almost complex structure $T$ on $B_2(0)\subset \R^{2n}=\R^m$ with energy bounded by $\Lambda$.
\begin{lemma}If $T$ is $(m-2)$ homogeneous admissible almost complex structure with bounded energy over $B_2(0)$, then $T$ is a constant almost complex structure.  We also have the following effective version, for all $\epsilon>0$, there exist $\delta=\delta(n, \epsilon, \Lambda)>0$  such that if $T$ is $(\delta, 2, m-2)$-homogeneous with energy bounded by $\Lambda$, 
then $T$ is $(\epsilon, 2, m)$-homogeneous. 
\end{lemma}

\begin{proof}The proof is similar to \cite{CN1}[Lemma 2.5]. Suppose $T$ is $(m-2)$-homogeneous, we need to show $D T=0$, or in other words,
\[
\int_{B_2} |DT|^2=0
\]
If this were to fail, by radial invariance of $T$, we can assume there exists some point $e\in \R^m$ with norm one, and some $c>0$, 
\[
\int_{B_{1/4}(e)} |DT|^2=c>0.
\]
For $j=1, 2, \cdots,$ the balls $B_{4^{-j-1}}(4^{-j} e)$ are mutually disjoint. By radial invariance again,
\begin{equation}\label{energy-1}
\int_{B_{4^{-j-1}}(4^{-j} e)} |DT|^2=c 4^{-j(m-2)}
\end{equation}
Suppose $T$ is $(m-2)$-homogeneous with respect to $V^{m-2}$, then there exists $c_1>0$, we can construct, for each $j$, at least a collection of at least $c_1 4^{j(m-2)}$ mutually disjoint balls contained in $B_2(0)$, with centers in $V^{m-2}$ and radius $4^{-j-1}$. By translation invariance, for each of these balls, \eqref{energy-1} holds. Note that for different $j$, all corresponding balls are also mutually disjoint. But this would implies that  $\int_{B_2} |DT|^2\geq \sum_j c c_1=\infty$, contradicting finite energy assumption.

Now we prove the quantitative version. Suppose this were not true, then there exists $\epsilon>0$ and a sequence of $(i^{-1}, 2, m-2)$-homogeneous admissible almost complex structures $T_i$ defined over $B_2(0)\subset \R^m$ with uniformly bounded energy, but none of them is $(\epsilon, 2, m)$-homogeneous. We can assume that $T_i$ is $(i^{-1}, 2, m-2)$-homogeneous with respect to a $(0, 2, m-2)$-homogeneous admissible almost complex structure $K_i$. By orthogonal rotation if necessary, we can assume $K_i$ are $(m-2)$-homogeneous with respect to the same $(m-2)$ plane in $\R^m$ (Correspondingly, we need to do the orthogonal rotation in coordinates for $T_i$). By passing to a subsequence,  $T_i$ converges to $T_\infty$ in $L^2$, and $T_\infty$ has bounded energy. In particular, $T_\infty$ is then $(2i^{-1}, 2, m-2)$-homogeneous with respect to $K_i$ for each $i$ hence $T_\infty$ is $(0, 2, m-2)$-homogeneous, but $T_\infty$ is not $(\epsilon, 2, m)$-homogeneous.   This contradicts the fact above. 
\end{proof}

Now we are ready to prove Theorem \ref{improve1}.
\begin{proof}[Proof of Theorem \ref{improve1}] We argue by contradiction. Fix a point $p\in M$, for a sequence of $\epsilon_i$ and $r_i$ both converging to zero, there exists an energy-minimizing almost complex structure $J_i$, such that $J_i$ is $(\epsilon_i, r_i, m-2)$-homogeneous over $B_{r_i}(p)$, but $r^{J_i}(p)<r_i/8$. Then  by Theorem \ref{com1} $T_{p, r_i}J_i(z)$ converges to (in $W^{1, 2}$ norm over compact sets of $\R^{2n}$), by subsequence if necessary, an energy-minimizing almost complex structure $T$ on $\R^{2n}$. Moreover, since $J_i$ is $(\epsilon_i, r_i, m-2)$-homogeneous, then there exists $T_i$, a $(m-2)$ homogeneous almost complex structure, such that
\[
\int_{B_2(0)} |T_{p, r_i}J_i-T_i|^2<\epsilon_i
\]
When $i$ is sufficiently large, this implies that
\[
\int_{B_2(0)}|T-T_i|^2 <2\epsilon_i. 
\]
Since $T_i$ is $(m-2)$-homogeneous, we know that $T$ is $(2\epsilon_i, 1, m-2)$-homogeneous for all $i$, hence it is $(0, 1, m-2)$-homogeneous. This implies that $T$ is a constant almost complex structure. Hence it follows that
\[
\int_{B_2} |T_{p, r_i}J_i- T|^2\rightarrow 0
\]
when $i$ sufficiently large, where $T$ is a constant almost complex structure. 
Let $\l_i$ denote the average of $T_{p, r_i}J_i$, then clearly 
\[
\int_{B_2} |T_{p, r_i}J_i- \l_i|^2\leq \int_{B_2} |T_{p, r_i}J_i- T|^2
\]

Note that $T_{p, r_i}J_i(z)=J_i(r_iz))$, we can then apply the improved $\epsilon$-regularity theorem (see Theorem \ref{improve}), 
to $\tilde J_i(z):=T_{p, r_i}J_i(z)$, 
for sufficiently large $i$, $\tilde J_i(z)$ has uniformly bounded $C^\alpha$ norm (hence smooth with uniformly $C^k$ norm) in $B_{3/16}$. In particular we have the uniform estimate, for $i$ sufficiently large, 
\[
\sup_{B_{3/16}} |\nabla \tilde J_i|+|\nabla^2 \tilde J_i|\leq 1.\]
This implies that $r^{\tilde J_i}(0)\geq 3/16$. By scaling this implies $r^{J_i}(p)\geq 3r_i/16$ for $i$ sufficiently large. This contradiction completes the proof. 
\end{proof}

Now we can finish the proof of the main theorem in this section,

\begin{proof}[Proof of Theorem \ref{volumeestimate}]We only need to prove for $r<r_0$, where $r_0$ is the uniform constant in Theorem \ref{improve1}. Given $r<r_0$, then for $p\in \cB_{r/8}(J)$, $r^J(p)<r/8$,  $J$ is not $(\epsilon, 2r, m-2)$ homogeneous by Theorem \ref{improve1}. In other words, $p\in \cS^{m-3}_{\eta, 2r}$ for any $\eta\in (0, \epsilon]$. Therefore 
\[
\cT_{r/8}(\cB_{r/8}(J))\subset \cT_{r/8}(\cS^{m-3}_{\eta, 2r}). 
\]
Taking $\tilde r=r/8$, then we get, by Theorem \ref{quantitative},  
\[
\text{Vol}(\cT_{\tilde r}\cB_{\tilde r}(J)\cap B_1)\leq \text{Vol} (\cT_{\tilde r}(\cS^{m-3}_{\eta, 16 \tilde r}))\leq C \tilde r^{3-\eta}.
\]
Given the volume estimate above, the Minkowski dimension bound is a routine argument by definition. 
\end{proof}

\subsection{Proof of Theorem \ref{quantitative}}
We now prove Theorem \ref{quantitative} in this section. The proof follows the arguments in \cite{CN, CN1} and it is technically involved. The essential point is to decompose the set $\cS^k_{\eta, r}\cap B_1$ into subsets under various homogeneity scales. The growth rate of the numbers of these subsets (at various scales) is of polynomial growth, using a quantitative version of the fact that tangent cones are homogeneous degree zero together with the uniform energy bound. In particular, points in $\cS^k_{\eta, r}$ with same (small) homogeneity scales line up in a tubular neighborhood of some $k$-plane. This would allow us to find an effective covering of $\cS^k_{\eta, r}$ using small balls.

First we recall a cone-splitting principle for admissible almost complex structures on $\R^m$. 
\begin{prop}
Suppose an admissible almost complex structure $T$ on $\R^m$ is $k$-homogeneous at $y$ with respect to $k$-plane $V^k$, and $T$ is homogeneous of degree zero at $z$. Suppose $z-y\notin V^k$. Then $T$ is $(k+1)$-homogeneous at $y$ with respect to $(k+1)$-plane spanned by $z-y$ and $V^k$. 
\end{prop}
\begin{proof}We can assume $y=0$ by a translation. For any $t\in \R$, $x\in \R^m$, $v\in V^k$, we need to show that
\[
T(x+tz+v)=T(x)
\]
Write $t=\l-\l^{-1}$ for some $\l>0$. We compute
\[
\begin{split}
T(x+tz+v)=&T(x+tz)=T(\l (z+\l^{-2}(\l x-z)))=T(z+\l^{-2}(\l x-z))\\
=&T(z+\l x-z)=T(\l x)=T(x). 
\end{split}
\]
\end{proof}

We now state a quantitative version of cone splitting lemma.

\begin{lemma}[Cone splitting I]\label{cone1}
For any $\epsilon, \tau>0$, there exists $0<\delta=\delta(M, g, \epsilon, \tau)$ and $0<\theta=\theta(M, g, \epsilon, \tau)$ such that the following holds. For any energy-minimizing almost complex structure $J$ in $W^{1, 2}(\cJ_g)$ such that $J$ is $(\delta, r, k)$-homogeneous at $p\in M$, for $r\leq \theta$ and there exists $y\in B_r(p)\backslash \cT_{\tau r}(V^k_{J, p})$ such that $J$ is $(\delta, r, 0)$-homogeneous at $y$. Then $J$ is $(\epsilon, r, k+1)$-homogeneous at $p$.
\end{lemma} 

\begin{proof}We argue by contradiction. Suppose for some given $\epsilon, \tau>0$, taking a sequence $\delta_i, \theta_i\rightarrow 0$, there exists a sequence of minimizing almost complex structures $J_i$ violating the claim: $J_i$ is $(\delta_i, r_i, k)$-homogeneous at $p$ for $r_i\leq \theta_i$ and there exists $y_i\in B_{r_i}(p)\backslash \cT_{\tau r_i}(V^k_{J_i, p})$ such that $J_i$ is $(\delta_i, r_i, 0)$-homogeneous at $y_i$, but $J_i$ is not $(\epsilon, r_i, k+1)$-homogeneous at $p$. Use the exponential map $\exp_p$ to identify $B_1(p)$ with $B_1(0)$. We also assume the defining plane of $J_i$ are all the same, otherwise we perform a rotation. 
Consider dilation by $z\rightarrow r_iz$ and the sequence $\tilde J_i(z)=J_i(r_i z)$. Then $\tilde J_i$ is $(\delta_i, 1, k)$-homogeneous at $0$ and there exists $\tilde y_i\in B_1(0)\backslash \cT_{\tau}(V^k_{\tilde J_i, 0})$ such that $\tilde J_i$ is $(\delta_i, 1, 0)$-homogeneous at $\tilde y_i$, but $\tilde J_i$ is not $(\epsilon, 1, k+1)$-homogeneous. By passing to subsequence, we can assume $\tilde J_i$ converges strongly to an energy-minimizing almost complex structure $T_\infty$ on $\R^{2m}$ of homogeneous degree zero. In particular $T_\infty$ will be $(0, 1, k)$-homogeneous at $0$ (that is $k$-homogeneous by definition). Moreover, we can assume $\tilde y_i$ converges to $\tilde y_\infty\subset B_1(0)\backslash \cT_{\tau}(V^k_{T, 0})$ such that $T_\infty$ is $(0, 1, 0)$-homogeneous at $y_\infty$. By cone-splitting principle, we know that $T$ is $(k+1)$-homogeneous with respect to $\text{span}\{y_\infty, V^k_{T, 0}\}$. Since $\tilde J_i$ converges to $T$ strongly, this contradicts that $\tilde J_i$ is not $(\epsilon, 1, k+1)$-homogeneous. 
\end{proof}

\begin{lemma}[Cone splitting II]\label{cone2}Given $\eta, \tau>0$, there exist $\epsilon=\epsilon(M, g, \eta, \tau)$ and $\theta=\theta(M, g, \eta, \tau)$ such that the following holds. Let $J\in W^{1, 2}(\cJ_g)$ be an energy-minimizing almost complex structure. Suppose $x\in \cS^k_{\eta, r}$ and $J$ is $(\epsilon, r, 0)$ homogeneous at $x$ for $r\in (0, \theta]$. Denote
 \[
 H_{x, r}:= \{y\in B_r(x): J\; \text{is}\; (\epsilon, r, 0)\text{-homogeneous at}\; y\}.
 \]
 Then there exists some $k$-plane $V^k$, $H_{x, r}\subset \cT_{\tau r}(V^k_{J, x})$. 
\end{lemma}

\begin{proof}
Given some fixed $\eta, \tau>0$, we can find a sequence $J_i$, $x_i$ and $r_i\rightarrow 0$,  such that $x_i\in \cS^k_{\eta,  r_i}$, but for any $k$-plane $V^k$, $H_{x_i, r_i}$ contains at least a point outside $\cT_{\tau r_i}(V^k_{J_i, x_i})$. Then we can construct inductively a collection of $k+1$ points $(y^1_i, \cdots, y^{k+1}_i)$ in $H_{x_i, r_i}$, with $\tilde y^j_i=\exp_{x_i}^{-1}(y^j_i)$ such that the following holds. 
Denote $V^0=\{0\}$ and $V^j=\text{span}\{0, \tilde y^1_i, \cdots, \tilde y^j_i\}$ for $j=1, \cdots k+1$, then $y^j_i\in B_{r_i}(x)\backslash \cT_{\tau r_i}(V^{j-1})$, for $j=1, \cdots, k+1$. 
We claim that $J_i$ is $(\eta, r_i, k+1)$-homogeneous at $x_i$ for $i$ sufficiently large, hence contradicts $x_i\in \cS^k_{\eta, r_i}$. 

To prove the claim, we apply Lemma \ref{cone1} inductively. 
For $\epsilon, \tau$ fixed in Lemma \ref{cone1}, consider $\delta=\delta(M, g, \epsilon, \tau)$ determined by Lemma \ref{cone1}. Set $\delta^{[0]}=\eta$ and define inductively $\delta^{[j+1]}=\delta(M, g, \delta^{[j]}, \tau)$ for $j=0, \cdots, k$.  Take $\epsilon$ to be $\delta^{[k+1]}$ (the smallest number among $\delta^{[j]}$ for $j=0, 1, \cdots, k+1$). 
We can then apply Lemma \ref{cone1} inductively to $x_i, V^j$ for $j=0, \cdots, k$ with corresponding $\epsilon=\delta^{[k-j]}$ and $\tau$, to conclude that $J_i$ is $(\delta^{[k-j]}, r_i, j+1)$-homogeneous at $x_i$, for $j=0, 1, \cdots, k$. That is, $J_i$ is $(\delta^{[k+1-j]},  r_i, j)$-homogeneous at $x_i$ (with respect to $j$-plane $V^j$), and $J_i$ is $(\delta^{[k+1-j]}, r_i, 0)$ homogeneous at $y^{j+1}_i$, we can conclude by Lemma \ref{cone1} that $J_i$ is $(\delta^{[k-j]}, r_i, j+1)$-homogeneous at $x_i$ (with respect to $(j+1)$-plane $V^{j+1}$). 
This completes the proof. 
\end{proof}

Lemma \ref{cone2} follows Breiner-Lamm \cite{BL}, where the authors apply the quantitative technique in \cite{CN1} to biharmonic maps. 
We will also need the following quantitative rigidity, which is a quantitative description of the fact that any tangent almost complex structure is homogeneous degree zero. 

\begin{prop}[Quantitative rigidity]\label{qr}
For every $\epsilon>0$, there exists $\zeta=\zeta(M, g, \epsilon)$ and $\tau_1=\tau_1(M, g, \epsilon)$  such that for any energy minimizing almost complex structure $J\in W^{1, 2}(\cJ_g)$, and if for any $0<r\leq 1$, $\beta\leq \tau_1$, 
\[
\tilde \Theta_{J}(x, r)-\tilde \Theta_{J}(x, \beta r)\leq \zeta, 
\]
then $J$ is $(\epsilon, r, 0)$-homogeneous at $x$. 
\end{prop}

\begin{proof}This is a straightforward argument by contradiction. \end{proof}

\begin{rmk}We use the energy function $\tilde \Theta_J(x, r)$ instead of $\Theta_{J}(x, r)$ since $\tilde \Theta_J(p, r)$ is monotone in $r$, satisfying \eqref{m23}. But $\Theta_J(x, r)$ is only almost monotone (using the monotonicity of $\tilde \Theta_J(x, r)$) and it is more convenient to use $\tilde \Theta_J(x, r)$ by its monotonicity. \end{rmk}

From now on we consider an energy-minimizer $J$ and consider $J$ over $B_1(p)$ for a fixed point $p\in M$. Given $\eta>0$, we want to decompose the points $B_1(p)$ into subsets at various homogeneity scales. For any $x\in B_1(p)$, denote for $0\leq s<t\leq 1$,
\[
W(x, s, t):=\tilde \Theta_J(x, t)-\tilde \Theta_{J}(x, s)\geq 0.
\]
Note by the monotonicity of $\tilde \Theta_J(x, r)$, $W(x, s_1, t_1)+W(x, s_1, t_2)\leq W(x, s_1, t_2)$ for $t_2\leq s_1$. 
In particular, let $(s_i, t_i)$ denote a sequence of intervals $[s_i, t_i]$ with $t_{i+1}\leq s_i$ and $t_1=1$, then
\begin{equation}\label{number1}\sum W(x, s_i, t_i)\leq \tilde \Theta_J(x, 1)\leq \Lambda,\end{equation}
where $\Lambda=\Lambda(M, g, n)$ is a uniformly bounded constant. 
Given a fixed number $\gamma\in (0, 1)$ and a fixed positive integer $q$, for a fixed positive number $\zeta$, consider all positive integers $j$
\[
W(x, \gamma^{j+q}, \gamma^j)>\zeta
\]
Let $K$ be the number of all such $j$s, then 
\begin{equation}\label{number2}K\leq (q+1)\zeta^{-1} \Lambda\end{equation} 
Otherwise, there exist at least $\zeta^{-1} \Lambda$ \emph{disjoint} closed intervals of the form $[\gamma^{j+q}, \gamma^j]$ with $W(x, \gamma^{j+q}, \gamma^j)>\zeta$, contradicting \eqref{number1}. 

Now we decompose the points in $B_1(p)$ into subsets as follows. 
For each $x\in B_1(0)$, we can also associate a $j$-tuple $T^j(x)=(a_1, \cdots, a_j)$ as follows. For its $i$-th entry $a_i=0$, if 
$W(x, \gamma^{j+q}, j)\leq \zeta$, otherwise $a_i=1$. Note that by \eqref{number2} $T^j(x)$ has at most $K$ entries of value $1$. 
For a fixed $j$-tuple $T^j$ with entires either zero or one, we put
\[
E_{T^j}=\{x\in B_1(p): T^j(x)=T^j\}.
\]
A priori there could be $2^j$ possible sets $E_{T^j}$. However, by \eqref{number2}, we see that if $E_{T^j}\neq \emptyset$, then $T^j$ contains at most $K$ entries of value $1$, hence we have at most $\begin{pmatrix}j\\ K\end{pmatrix}<j^K$ many nonempty subsets. As a consequence, we can decompose $\cS^k_{\eta, \gamma^j}$ accordingly to the subsets $\cS^k_{\eta, \gamma^j}\cap E_{T^j}$.
Next we show we can cover the set $\cS^k_{\eta, \gamma^j}\cap E_{T^j}$ in an effectively way.

We first specify the choice of constants. Given $\eta>0$ and $\gamma\in (0, 1)$, as in Lemma \ref{cone2}, with $\eta, \tau=\gamma$ fixed, we choose $\epsilon=\epsilon(M, g, \eta, \gamma)$ and $\theta=\theta(M, g, \eta, \gamma)$ such that Lemma \ref{cone2} holds.  With this $\epsilon$, we choose $\zeta$ and $\tau_1$ in Lemma \ref{qr} such that Lemma \ref{qr} holds. Now choose $q$ such that $\gamma^q\leq \tau_1$, hence if \[\tilde \Theta_{J}(x, \gamma^j)-\tilde \Theta_{J}(x, \gamma^{j+q})\leq \zeta\] then $J$ is $(\epsilon, \gamma^j, 0)$ homogeneous at $x$. With this $\zeta, q$ fixed, we can then define $E_{T^j}$ accordingly.  We also denote $K_0$ to be the smallest positive integer such that $\gamma^{K_0}\leq \theta$.
Then we have the following, 

\begin{lemma}[Covering lemma] Given $\eta>0$ and $\gamma\in (0, 1)$, there exists $c_0=c_0(m)$ and $j_0=j_0(M, g, \eta, \gamma)$ such that each set $\cS^k_{\eta, \gamma^j}\cap E_{T^j}$ can be covered by at most $ (c_0\gamma^{-m})^{j_0}(c_0\gamma^{-k})^{j-j_0}$ many balls of radius $\gamma^j$ (for $j\geq j_0$). 
\end{lemma}

\begin{proof} We use $T^{j-1}$ to denote the $(j-1)$-tuple by dropping the last entry of $T^j$ and we can get $T^0, \cdots, T^j$. Clearly with such a choice $E_{T^j}\subset E_{T^{j-1}}$, for $j=1, \cdots$. 

Now we construct the covering of $\cS^k_{\eta, \gamma^j}\cap E_{T^j}$ inductively.
Suppose we have covered $\cS^k_{\eta, \gamma^j}\cap E_{T^j}$ by balls of radius $\gamma^j$. Then for each such balls $B_{\gamma^j}$ of radius $\gamma^j$, we choose a minimal covering of  $\cS^k_{\eta, \gamma^{j+1}}\cap E_{T^{j+1}}\cap B_{\gamma^j}$ by balls of radius $\gamma^{j+1}$, centered at the points in $\cS^k_{\eta, \gamma^{j+1}}\cap E_{T^{j+1}}\cap B_{\gamma^j}$ (note that $\cS^k_{\eta, \gamma^{j+1}}$ is contained in $\cS^k_{\eta, \gamma^j}$.)

For each step, we can choose a minimal covering such that we need at most the number of balls by $c_0(m)\gamma^{-m}$. This would lead to the trivial bound $(c_0\gamma^{-m})^j$. To get the improved bound,  we observe that for $j\geq j_0=\max{(K, K_0)}$, one can get a much better estimate of the number of balls. First note that for $j\geq j_0$,  the $j$-th entry of $T^j$ has to be zero. This implies that for any $y\in \cap E_{T^j}$, $W(y, \gamma^{j+q}, \gamma^j)\leq \zeta$. Applying Lemma \ref{qr} (with $\gamma^j\leq \theta, \gamma^q\leq \tau_1$), $J$ is $(\epsilon, \gamma^j, 0)$-homogenous at $y$. Now for $x\in \cS^k_{\eta, \gamma^j}\cap E_{T^j}$, we apply Lemma \ref{cone2} (with $\tau=\gamma, r=\gamma^j$) to conclude that for any $y\in E_{T^j}\cap B_{\gamma^j}(x)$, it is contained in a $\cT_{\gamma^{j+1}}(V^k_{J, x})$, for some $k$-plane $V^k$. In this case, we can then cover $\cS^k_{\eta, \gamma^j}\cap E_{T^j}\cap B_{\gamma^j}(x)$ by the radius $\gamma^{j+1}$ balls with the number bounded by $c_0(m)\gamma^{-k}$.  This completes the proof. 
\end{proof}

Now we can complete the proof of Theorem \ref{quantitative}.

\begin{proof}
Fix $\eta>0$ and we choose $\gamma= \min{(c_0^{-\eta/2}, 1/2)}$. For $j\geq j_0$, we can cover $\cS^k_{\eta, \gamma^j}\cap B_1(p)$ by at most $j^K$ subsets, with each covered by at most $(c_0\gamma^{-m})^{j_0}(c_0\gamma^{-k})^{j-j_0}$ many balls of radius $r=\gamma^j$. We should note that since for each such ball, its $\gamma^j$-tubular neighborhood is simply the radius $2\gamma^{j}$ ball of the same center.  The covering gives the volume estimate,
\begin{equation}\label{final1}
\text{Vol}(\cS^k_{\eta, \gamma^j}\cap B_1(p))\leq j^K (c_0\gamma^{-m})^{j_0}(c_0\gamma^{-k})^{j-j_0} \omega_m (\gamma^j)^m
\end{equation}
Note that by our choice, we have
\[
c_0^j\leq (\gamma^j)^{-\eta/2}, j^K\leq c(\eta, K) (\gamma^j)^{-\eta/2}, 
\]
Then by \eqref{final1} we can get the volume estimate,
\begin{equation}\label{final2}
\text{Vol}(\cS^k_{\eta, \gamma^j}\cap B_1(p))\leq C(M, g, \eta) (\gamma^j)^{m-k-\eta}
\end{equation}
This proves the result for $r=\gamma^j$, $j\geq j_0$. For $r\leq \gamma^{j_0}$, this is obvious by choosing a different constant $C$. For general $r\in (\gamma^{j+1}, \gamma^j)$, the estimate follows by the fact that 
$\cS^{k}_{\eta, r}\subset \cS^k_{\eta, \gamma^j}$. 
Regarding the $r$-tubular neighborhood,  we can replace the covering $r$-balls by $2r$-balls. This gives the volume estimate, 
\[
\text{Vol}(\cT_r(\cS^k_{\eta, r})\cap B_1(p))\leq  2^mC(M, g, \eta) r^{m-k-\eta}
\]
This completes the proof. 
\end{proof}

\begin{rmk}
We note that Naber-Valtorta \cite{NV, NV1} recently improve the volume estimate to the order $r^{m-k}$, and prove that $\cS^k_\eta$ is $k$-rectifiable and  some other deep results for energy-minimizers in the setting of (almost) harmonic maps. It is of no doubt that their methods and corresponding results hold in our setting. Since these results are very technically involved, we will discuss Naber-Valtorta type results elsewhere. 
\end{rmk}

\section{Energy-minimizing almost complex structure on Sobolev spaces}
In this section we study energy-minimizing almost complex structures on some other Sobolev spaces, which are subspaces of  \emph{admissible} almost complex structures. First we discuss an example on $S^{2n}$ for $n\geq 2$. 

\subsection{Examples of energy-minimizing admissible almost complex structures}

Consider the sphere $S^m=S^{2n}$ with the round metric, which is conformally flat. It is well-known that $S^{2n}$ admits no smooth almost complex structures for $n\neq 1, 3$. We describe the round sphere $S^{2n}$ in the coordinate on $\R^{2n}$ by stereographic projection, and hence the metric is of the form
\[
g_0=\frac{4}{(1+|x|^2)^2} ds^2. 
\]
Let $J_0$ be the standard  almost complex structure on $\R^{2n}$, given by
\[
J_0\p_i=\p_{n+i}, J_0\p_{n+i}=-\p_i, i=1, \cdots, n.
\]
Clearly $J_0$ is compatible with $g_0$ on $\R^{2n}$. We need to check that $J_0\in W^{1, 2}$.  We compute
\[
\nabla_iJ_j^k=\p_{i}J^k_j-\Gamma_{ij}^l J_l^k+\Gamma_{il}^kJ^l_j=-\Gamma_{ij}^l J_l^k+\Gamma_{il}^kJ^l_j. 
\]
A direct computation gives 
\[
|\nabla J_0|^2=c_n|x|^2
\]
Hence for any $p<2n$, we have
\[
\int_{S^{2n}} |\nabla J_0|^pdv\sim \int_0^\infty r^{p+2n-1}(1+r^2)^{-2n}dr<\infty.\]
Similarly we can compute \[\int_{S^{2n}}|\nabla^2 J_0|^{p/2}<\infty.\]
In particular $J_0$ is an admissible almost complex structure on $(S^{2n}, g_0)$ with an isolated singularity. 
A similar approximation as in Lemma 4.2 gives the following, 
\begin{prop}\label{obstruction1}A compact Riemannian manifold of even dimension $m=2n$ admits a smooth almost complex structure if and only if it admits an admissible almost complex structure in $W^{1, m}$ (or $W^{2, m/2}$). 
\end{prop}

Hence it might be preferably to consider the energy-minimizing problem over $W^{1, p}(\cJ_g)$ for the $p$-energy functional, in particular for $p=m$,
\[
 E_p(J)=\int_M |\nabla J|^pdv.
\] 
Note that the obstruction to the existence of smooth almost complex structures can be detected by its Chern classes. For admissible almost complex structure in $W^{2, m/2}$, one can consider an associate connection $\tilde \nabla$ such that its curvature form $F_{\tilde \nabla}$ is in $L^{m/2}$. Chern-Weil theory can work for such objects, and hence 
all the obstructions involved with Chern classes apply.  This gives an intuitive interpretation that the same obstruction applies to admissible almost complex structure in $W^{2, m/2}$. 
In particular the example on $S^{m}$ shows that $p=m$ in Proposition \ref{obstruction1} is optimal. 
From the point view of almost Hermitian geometry, $\nabla J$ only contains information of ``torsion" of an almost complex structure, while the second derivative contains information of its curvature and hence all necessary topological information, such as Chern classes. 
Hence it is natural to study energy functional of second order, such as
\[
\cE_2(J)=\int_M|\Delta J|^2 dv
\] 

We continue to discuss $(S^{2n}, g_0, J_0)$. A less obvious fact is that $J_0$ on $(S^{2n}, g_0)$ constructed above is indeed energy-minimizer. This follows essentially Bor-Hern\'andez-Lamoneda-Salvai  \cite{BorL, BorSalvai}. We introduce some notations. 
It would be more convenient to use the two form $\omega=g(\cdot, J\cdot)$ instead of $J$. Note that we have $|\nabla \omega|^2=|\nabla J|^2$ holds pointwise. 
One applies the well-known Gray-Hervella decomposition \cite{GH} to the bundle $\cW:=\Lambda^1\otimes [\Lambda^{2, 0}]$ on almost Hermitian manifolds into four subbundles
\[\cW=\cW_1\oplus \cW_2\oplus \cW_3\oplus \cW_4
\]
corresponding to the decomposition of $\cW$ into irreducible $U_n$ modules. One can then decompose $\nabla \omega$ accordingly 
\[
\nabla \omega=\sum_{i=1}^4(\nabla \omega)_i
\] 
Denote $e_i(\omega)=\int_M |(\nabla\omega)_i|^2$, then we have $E(J)=\sum e_i(\omega)$.
\begin{thm}(\cite{BorSalvai})\label{BLS}Let $(M^{2n}, g)$ be a compact Riemannian manifold such that it is conformally flat, or anti-self dual (when $n=2$), then for any $J\in \cJ_g$, 
\begin{equation}\label{bls3}
2\int_M |(\nabla \omega)_1|^2dv-\int_M |(\nabla \omega)_2|^2dv+(n-1)\int_M |(\nabla \omega)_4|^2dv=\frac{n-1}{2n-1}\int_M R_g dv,
\end{equation}
where $R_g$ is the scalar curvature of $g$. 
\end{thm}
The above theorem holds  if $J\in W^{1, 2}(\cJ_g)\cap W^{2, 1}(\cJ_g)$ (hence $\omega\in W^{1, 2}\cap W^{2,1}\cap L^\infty$). 
 Hence by the regularity we have obtained, it holds for any minimizer in $W^{1, 2}(\cJ_g)$ (even when $\cJ_g$ is empty). The proof of Theorem \ref{BLS} involves a Bochner formula, 
 \[
 \Delta_d\omega=-\Delta \omega+E(\omega),
 \]
where $E$ is a curvature operator acting on $\omega$ linearly, which is certainly applicable when $J\in W^{1, 2}(\cJ_g)\cap W^{2, 1}(\cJ_g)$. 
The curvature condition appears since in this case, 
\begin{equation}\label{bls1}
(E(\omega), \omega)=c_n R_g,
\end{equation}
where $R_g$ is the scalar curvature of $g$. 
The other ingredient is that under Gray-Hervella decomposition, there are linear relations (with suitable constants $a_i, b_i$) 
\begin{equation}\label{bls2}
|d\omega|^2=\sum a_i |(\nabla \omega)_i|^2,\; |\delta \omega|^2=\sum_i b_i|(\nabla \omega)_i|^2
\end{equation}
Note that \eqref{bls1} and \eqref{bls2} are algebraic identities which hold pointwise, hence are applicable for $J\in W^{1, 2}(\cJ_g)$. In other words, \eqref{bls3} holds when $J\in W^{1, 2}(\cJ_g)\cap W^{2, 1}(\cJ_g)$. In particular the main theorem \cite{BorSalvai}[Theorem 1] holds for $J\in W^{1, 2}(\cJ_g)\cap W^{2, 1}(\cJ_g)$. Bor-Hern\'andez-Lamoneda-Salvai constructed many examples of non-K\"ahler energy minimizers (in $\cJ_g$) using these results, see \cite{BorSalvai}[Theorem 1]. Hence we have the following, 

\begin{prop}All examples constructed in \cite{BorSalvai}[Section 3] are energy-minimizers in $W^{1, 2}(\cJ_g)$. In other words, $E_0(g)=E^s_0(g)$ (there is no energy gap).  \end{prop}
Applying \cite{BorSalvai}[Theorem 1], we have this special example, 
\begin{prop}The standard complex structure $J_0$ on $\R^{2n}$ is an energy-minimizer in $W^{1, 2}(\cJ_g)$ on $(S^{2n}, g_0)$. When $n\neq 3$, all energy minimizers are of type $\cW_4$ and when $n=3$, all minimizers are of type $\cW_1\oplus \cW_4$. 
\end{prop}

\begin{proof}Note that $(g_0, J_0)$ is locally conformally K\"ahler, hence  $\nabla \omega_0\in \cW_4$. Let $J$ be an energy-minimizer in $W^{1, 2}(\cJ_g)$, then $J\in W^{1, p}\cap W^{2, p/2}$ for $p\in[2, 3)$. Hence Theorem \ref{BLS} is applicable and we have, for $n\geq 3$
\[
E(J)=\sum e_i(\omega)\geq \frac{2}{n-1}e_1(\omega)-\frac{1}{n-2}e_2(\omega)+e_4(\omega)=C(M, g)=e_4(\omega_0)=E(J_0)
\]
When $n=2$, $e_1=0$, hence the argument also works. 
\end{proof}

\begin{rmk}When $n=3$, the standard nearly K\"ahler structure on $S^6$ gives a smooth minimizer \cite{BorSalvai}. 
\end{rmk}

\subsection{Energy-minimizing almost complex structure in smooth category}
We have  studied the minimizing problem over  \emph{admissible almost complex structures} denoted by $W^{1, 2}(\cJ_g)$, which can be viewed as a companion of Sobolev maps between Riemannian manifold $W^{1, 2}(M, N)$.
The Sobolev embeddings between Riemannian manifolds have been extensively studied by many authors, whose study were motivated by the study of harmonic maps. We refer the readers to  in particular the papers by Bethuel \cite{B}, Brezis-Li \cite{BLi}, Hang-Lin \cite{HangLin1, HangLin2} and the references therein for the theory of Sobolev mappings.
 
\begin{exam}
Suppose $(M, g)$ is an Euclidean ball $B^{2n}$ with its flat metric (in this example $M$ is not complete nor compact), then almost complex structures which are compatible with the flat metric is modeled on $SO(2n)/U(n)$. Hence in this case $W^{1, 2}(\cJ_g)=W^{1, 2}(B^{2n}, SO(2n)/U(n))$, where $SO(2n)/U(n)$ is endowed with the standard metric. When $n=2$ for example, this coincides the space $W^{1, 2}(B^4, S^2)$, the Sobolev mappings. 
\end{exam}

Admissible almost complex structure is a most general extension of the concept of (smooth) almost complex structures. It is well-known that there are obstructions to the existence of a smooth almost complex structures. For example, only $S^2$ and $S^6$ admit smooth almost complex structures among all spheres. But we have seen $S^{2n}$  admits admissible almost complex structures for any $n$.  In particular $\cJ_g$ is not dense in $W^{1, 2}(\cJ_g)$ in general (for $m\geq 4$). 

\begin{pl}Given a compact Riemannian manifold $(M, g)$ of even dimension, is there always an admissible almost complex structure in $W^{1, 2}(\cJ_g)$?
\end{pl}

We believe the answer should be affirmative and it should be interesting to study these objects. 
On the other hand, to study the minimizing problem on the smooth category, we define the following two spaces,
\[\begin{split}H^{1, 2}_s(\cJ_g)&=\{\text{the closure of}\; \cJ_g\;\text{ in}\; W^{1, 2}(\cJ_g)\; \text{in the \emph{strong}}\; W^{1, 2} \; \text{topology}\}\\
H^{1, 2}_w(\cJ_g)&=\{\text{the closure of}\; \cJ_g\;\text{ in}\; W^{1, 2}(\cJ_g)\; \text{in the \emph{weak}}\; W^{1, 2} \; \text{topology}\}
\end{split}\]
Clearly we have the inclusions
\[
H^{1, 2}_s(\cJ_g)\subset H^{1, 2}_w(\cJ_g)\subset W^{1, 2}(\cJ_g). 
\]
When the manifold is $S^4$, we have seen that $H^{1, 2}_s(\cJ_g)=H^{1, 2}_w(\cJ_g)=\emptyset$ but $W^{1, 2}(\cJ_g)$ is not empty (this holds for any smooth Riemannian metric using Proposition \ref{matrix1}). 
When  the inclusions are strict or not is difficult problem and we refer the readers to \cite{HangLin2} and the references therein for results for Sobolev mappings.

We roughly discuss energy-minimizing problem over the spaces  $H^{1, 2}_w(\cJ_g)$. Note that $H^{1, 2}_s(\cJ_g)$ does not have the nice compactness properties. Hence it would be a rather hard problem in general to consider whether there exists a minimizer in $H^{1, 2}_s(\cJ_g)$. Regarding energy-minimizing problem over $H_w^{1, 2}$, we have the following result similar to Theorem \ref{existence1},

\begin{thm}Let $(M, g)$ be a compact almost Hermitian manifold with a compatible almost complex structure. Then there exists at least one minimizer of $E(J)$ over $H^{1, 2}_w(\cJ_g)$ and it is a weak solution of \eqref{H1} in the sense of \eqref{H00}. Hence if $J$ is continuous, then it is smooth. Moreover, for any sequence of energy-minimizers $\{J_i\}\in H^{12, }_w(\cJ_g)$, there exists a subsequence $\{J_{k_i}\}$ and an energy-minimizing harmonic almost complex structure $J$ such that $J_{k_i}\rightarrow J$ strongly in $W^{1, 2}$-norm. \end{thm}

\begin{proof}Note that we assume that $\cJ_g$ and hence $H^{1, 2}_w(\cJ_g)$ is not empty. Denote
\begin{equation}
E^w_0(g)=\inf \{E(J): J\in H^{1, 2}_w(\cJ_g)\}. 
\end{equation}
Find a minimizing sequence $J_k\in H^{1, 2}_w(\cJ_g)$. Then $J_k$ has uniformly bounded $W^{1, 2}$-norm, hence it converges weakly to $J$ in $W^{1, 2}(M, TM\otimes T^*M)$ by subsequence. The convergence is weak in $W^{1, 2}$ and strongly in $L^2$ (by Rellich's compactness). It implies the pointwise convergence almost everywhere. In particular, $J$ satisfies the constraint $J^2=-id,$ and $g(J\cdot, J \cdot)=g$ almost everywhere. Hence $J\in W^{1, 2}(\cJ_g)$.

By the weak convergence we know that $E(J)\leq \liminf E(J_k)=E^w_0(g)$. Since $J_k\in H^{1, 2}_w(\cJ_g)$, it follows that $J$ is also in $H^{1, 2}_w(\cJ_g)$ by weak compactness. 
Hence this forces $E(J)=E^w_0(g)$. It then follows that $J_k$ converges to $J$ strongly in $W^{1, 2}$ and $J$ is an energy-minimizer in $H^{1, 2}_w(\cJ_g)$. Next we need to show that $J$ satisfies $\Delta J-J\nabla J\nabla J=0$ in the weak sense, namely satisfying \eqref{H00} for $T\in L^\infty \cap W^{1, 2}$. The argument is more involved since we need to construct variation path in $H^{1, 2}_w(\cJ_g)$.  

We claim that for any $S\in L^\infty\cap W^{1, 2}$ satisfying \eqref{C3} a.e, $J(t)=J(id+t S)(id-tS)^{-1} \in H^{1, 2}_w(\cJ_g)$ when $|t|$ is sufficiently small. Let $J_k\in \cJ_g$ such that $J_k$ converges to $J$ weakly in $W^{1, 2}$ and strongly in $L^2$, and let $T_k\in C^\infty(M, TM\otimes T^*M)$ converge to $S$ strongly in $W^{1, 2}$-norm (hence pointwise to $S$ a.e). We can assume that $T_k$ has uniform $L^\infty$ norm since $S\in L^\infty$. Denote
$S_k=(T_k+J_kT_kJ_k)-g(T_k+J_kT_kJ_k)g^{-1}$. Then $S_k$ is smooth and satisfies the constraint $J_kS_k+S_kJ_k=0, g(S_kx, y)+g(x, S_ky)=0$. Clearly $S_k$ has uniformly bounded $L^\infty$ norm and its pointwise limit is 
\[
(S+JSJ)-g(S+JSJ)g^{-1}=4S. 
\]
We need to show $S_k$ converges to $4S$ weakly in $W^{1, 2}$. This is a routine check. For example, for any $R\in L^2(TM\otimes T^*M)$, we write
\[
(\nabla J_k T_k J_k-\nabla J S J, R)=((\nabla J_k-\nabla J) SJ, R)+(\nabla J_k (T_k-S)J, R)+(\nabla J_k T_k(J_k-J), R)
\] 
It converges to zero when $k\rightarrow \infty$. Now consider $J_k(t)=J_k(id+tS_k/4)(id-tS_k/4)^{-1}$. Then for $|t|$ small $J_k(t)\in \cJ_g$ and it converges to $J(t)$ pointwise a.e Similarly we can check that $J_k(t)$ converges to $J(t)$ weakly in $W^{1, 2}$ when $|t|$ is sufficiently small. Since $J$ is an minimizer, we have
\[
\p_t E(J(t))=0. 
\]
It follows the same argument as in Theorem \ref{existence1} that $J$ satisfies \eqref{H00}. 
\end{proof}

We should emphasize that the regularity theory for energy-minimizing almost complex structure in $W^{1, 2}(\cJ_g)$ cannot directly be applied to energy-minimizers in $H^{1, 2}_w(\cJ_g)$ since the comparison almost complex structures constructed do not belong in $H^{1, 2}_w(\cJ_g)$. 

\begin{rmk}
B. White \cite{White86, White88} discussed these problems for $H^{1, p}_w(M, N)$ in the setting of Sobolev mappings (among certain fixed homotopy class).
F.H. Lin \cite{Lin1} has studied the regularity problem for energy-minimizing $p$-harmonic maps in continuous (smooth) category. 
\end{rmk}

Another interesting problem is to study energy-minimizing directly in $\cJ_g$. Recall 
\[
E^s_0(g)=\inf\{E(J): J\in \cJ_g\}.
\]
Suppose $J_k\in \cJ_g$ is a minimizing sequence, then $J_k$ converges to $J$ weakly in $W^{1, 2}$ (by subsequence) for some $J\in H^{1, 2}_w(\cJ_g)\subset W^{1, 2}(\cJ_g)$. A main point is that 
\[
E(J)\leq \liminf E(J_k)=E^s_0(g)
\]
and $J$ might not realize $E^s_0(g)$ in general, due to the weak convergence. Following F. H. Lin \cite{Lin1, Lin2}, we consider a sequence of Radon measure
\[
\mu_i=|\nabla J_i|^2dv
\]
We may assume $\mu_i\rightharpoonup \mu$ as Radon measures (weak convergence). By Fatou's theorem, we may write $\mu=|\nabla J|^2+\nu$, where $\nu\geq 0$ is also a Radon measure, called the \emph{defect measure} in \cite{Lin1}.  It seems that the method in \cite{Lin1} can be applied to study the regularity of the weak limit $J$ and the \emph{defect measure} $\nu$ in our setting. We summarize as the following question, 
\begin{pl}Let $J_k$ be a sequence of smooth energy-minimizing almost complex structure in $\cJ_g$. By taking subsequence we have the $J_k\rightharpoonup J$ for some $J\in H^{1, 2}_w(\cJ_g)$ and 
\[
|\nabla J_k|^2dv\rightharpoonup |\nabla J|^2 dv+\nu
\]
where $\nu\geq 0$ is also a Radon measure. It should be true that $J$ is a weakly harmonic map. 
The question is whether $J$ is an energy-minimizer in $H^{1, 2}_w(\cJ_g)$, namely $E(J)=E^w_0(g)$.  How about the regularity of energy-minimizer in $H^{1, 2}_w(\cJ_g)$ with respect to $E^w_0(g)$? How about the regularity of a weak limit $J$ and defect measure $\nu$ coming from a minimizing sequence of smooth almost complex structure? Is the defect measure canonical (depending on $g$)? \end{pl}

These questions seem to be subtle. We will discuss  these problems elsewhere.

\section{Appendix}
In Appendix we discuss various topics related to harmonic almost complex structures. In our setting, we do not fix a homotopy class of almost complex structures. It would be interesting to discuss minimizing over a fixed homotopy class. First we show that in general metric compatibility does not impose any homotopy restriction on almost complex structures, contrary to the setting that we fix a symplectic form. Then we roughly summarize the twistor geometry, where an almost complex structure can be viewed as a section of a twistor bundle. We briefly discuss the harmonic almost complex structure via the approach of twistor geometry, used by C. Wood \cite{Wood1} and others, see Davidov \cite{Davidov}. 
We also  derive a Bochner formula for smooth harmonic almost complex structures. In the end we give a self-contained proof of the fact that a continuous weakly harmonic almost complex structure is smooth. 

\subsection{Homotopy classes and the heat flow}
Suppose $(M, g)$ is a compact almost Hermitian manifold and denote $\cJ_g$ to be the space of smooth compatible almost complex structures which induce the orientation of $(M, g)$. A basic problem is whether $\cJ_g$ is path connected. The answer is no in general. As a comparison, it is well-known that if we fix a non-degenerate two form $\omega$, then the space of all compatible almost complex structures is contractible. Let $J$ be an almost complex structure. By Proposition \ref{matrix1}, we can construct a unique almost complex structure $\bar J$ (out of $J$) such that $\bar J$ is compatible with $g$. Then $J$ and $\bar J$ are in the same homotopy class.

\begin{prop}There exists a smooth path of $J(t)$ such that $J(0)=J$ and $J(1)=\bar J$.
\end{prop}
\begin{proof}
We can construct $J(t)$ pointwise since the construction is canonical. Given $p\in M$, we choose a coordinate such that $g=id$ at $p$. Then we write $J=A+S$, the skew and symmetric decomposition. Then we have $AS+SA=0$ and $A^2=-id-S^2$. Denote $Q$ to be the square root of $id+S^2$, then we have $\bar J=Q^{-1}A$. 
We should emphasize that $Q$ is a polynomial of $id+S^2$, hence it commutes with both $S$ and $A$ (since $A$ commutes with $id+S^2$). Moreover $Q\geq id$ since $id+S^2\geq id$. 
Now for $t\in [0, 1]$ let
\[N(t)=tJ+(1-t)\bar J\]
We compute
\[
N^2=(tJ+(1-t)\bar J)^2=-(t^2+(1-t)^2)id+t(1-t)(J\bar J+\bar J J).
\]
Now we have
\[
J\bar J+\bar J J=(A+S) Q^{-1}A+Q^{-1}A(A+S)=2Q^{-1}A^2=-2Q. 
\]
It follows that
\[
N^2=-id-2t(1-t)(Q-id).
\]
Since $id+2t(1-t)(Q-id)\geq id$ is a positive definite symmetric matrix, there exists a unique square root $P(t)$ such that $P(t)^2=id+2t(1-t)(Q-id)$. Note that $P(t)$ is a polynomial of $id+2t(1-t)(Q-id)$, hence it commutes with $Q, A$ and $S$. In particular $P(t)$ commutes with $N(t)$. Denote $J(t)=P^{-1}(t)N(t)$, then $J(0)=J$, $J(1)=\bar J$ and $J(t)^2=-id$ by construction. 
\end{proof}

In other words, $\cJ_g$ contains the same connected components as the space of almost complex structures itself on $M$ (inducing the same orientation of $(M, g)$). In particular the space $\cJ_g$ may not be path-connected. A typical example is $S^2\times S^2$ or $S^1\times S^3$. Hence it makes sense to seek energy-minimizers in a fixed homotopy class. In the setting of harmonic map, Eells-Sampson \cite{ES} proved that there exists a smooth harmonic map $u: M\rightarrow N$ in each homotopy class (of smooth maps) if the curvature of $N$ is nonpositive, using the harmonic map heat flow. We will study the heat flow for almost complex structures in a separate paper,
\[
\p_t J=\Delta J-J\nabla_p J\nabla_pJ. 
\]
\subsection{Admissible harmonic almost complex structure on $S^4$}

It is well-known that $S^4$ does not support a smooth almost complex structures. Equivalently given a Riemannian structure on $S^4$, its twistor bundle does not admit a smooth global section. However, the $W^{1, 2}$ sections exist and these objects are just admissible almost complex structures. By our results, there always exist energy-minimizing admissible almost complex structures on $S^4$ with respect to any Riemannian metric and they form a compact space. 

Twistor geometry was introduced by R. Penrose, and later the Riemannian setting was considered by Atiyah-Hitchin-Singer \cite{AHS}. 
We briefly summarize the approach by C. Wood \cite{Wood1} using the twistor theory to discuss harmonic almost complex structures. The survey paper \cite{Davidov} is a good reference. 

Denote $\cJ(\R^{2n})$ to be the set of almost complex structures on $\R^{m}=\R^{2n}$ compatible with the Euclidean metric. Then $J\in \cJ(\R^{2n})$ satisfies $J^2=-id, J+J^t=0$. It is easy to see that $\cJ(\R^{2n})$ can be described as $O(2n)/U(n)$ and it is a symmetric space with a canonical metric.  It contains two connected components $\cJ_{+}(\R^{2n})$ and $\cJ_{-}(\R^{2n})$, corresponding to the almost complex structures on $\R^{2n}$ inducing two orientations of $\R^{2n}$.  Each of them has a homogeneous representation  $SO(2n)/U(n)$. Given a Riemannian manifold $(M, g)$ of dimension $m=2n$, the twistor bundle is then described as the associated bundle
\[
\mathbb{J}=O(TM)\times_{O(2n)} \cJ(\R^{2n})
\]
For us the more relevant one would be the positive twistor bundle
\[
\mathbb{J}_{+}=SO(TM)\times_{SO(2n)}\cJ_{+}(\R^{2n}).
\]
Then all compatible smooth almost complex structure can be identified with the smooth sections of the positive twistor bundle,
\[
\cJ_g=\Gamma(M, \mathbb{J}_{+})
\]
With the Levi-Civita connection, the tangent bundle $T\mathbb{J}_{+}$ has a splitting $\cV\oplus \cH$, the vertical and horizontal subspace. The vertical part $\cV$ can be identified with the tangent space of the fiber, and hence has a canonical metric induced from the fiber $SO(2n)/U(n)$, and the horizontal part $\cH$ is the lift of the metric $g$ on $TM$. Denote such a metric on $\mathbb{J}_{+}$ as $\cG$, then 
\[
\pi: (\mathbb{J}_{+}, \cG)\rightarrow (M, g)
\]
is a Riemannian submersion with a totally geodesic fiber. Given a section (a map) $J: M\rightarrow \mathbb{J}_{+}$, denote $dJ$ to be its tangent map. Then the energy can be defined as (in the setting of harmonic map) 
\[
\tilde E(J):=\int_M \|dJ\|^2_{g, \cG} dv
\]
It turns out that $\tilde E(J)=\int_M |\nabla J|^2dv+n\text{Vol}(M)$. Similarly, $W^{1, 2}(\cJ_g)$ can then be described as the $W^{1, 2}$ sections $s: M\rightarrow \mathbb{J}_{+}$.

Using this correspondence, we know that an energy-minimizing almost complex structure on $(S^4, g_0)$ will be of type $\cW_4$ (Proposition 6.3). Hence on its smooth locus, this gives a complex structure on $S^4\backslash \cS$, where $\cS$ is the singular set of Minkowski dimension one. In this case one should be able to use the technique of Hardt-Lin \cite{HL2} (or rather their results directly) to prove that the singular set if a union of a finite set and a finite family of H\"older continuous closed curves having at most a finite number of crossings.               

One might speculate whether these objects are interesting and have applications to study the smooth structure of $S^4$.

\subsection{Bochner formula for smooth harmonic almost complex structures}    
In this section we derive a Bochner formula for smooth harmonic almost complex structures, which we have used in smooth setting to derive a priori estimates. 
Suppose now we have a smooth harmonic almost complex structure $J\in \cJ_g$. 
\begin{prop}We derive a Bochner formula for $J$,
\begin{equation}
\frac{1}{2}\Delta |\nabla J|^2=|\nabla^2 J|^2-|\nabla_p J\nabla_pJ|^2+Ric(\nabla J, \nabla J)+2R^k_{pil, p} J^l_j \nabla_iJ^k_j+4R^k_{pil}\nabla_pJ^l_j \nabla_iJ^k_j
\end{equation}
\end{prop}        
 \begin{proof}We use the following curvature convention, 
 \[
 R(\p_i, \p_j)\p_k=R_{ijk}^l\p_l, R_{ijlk}=R_{ijp}^lg_{kl}, R_{ij}=R^l_{lij}=g^{kl}R_{kilj}
 \]
 and the following commutative derivatives,
 \[
 [\nabla_i, \nabla_j] T^k=R_{ijl}^kT^l, [\nabla_i, \nabla_j] S_k=-R_{ijk}^l S_l. 
 \]
 Choose a coordinate such that $g_{ij}=\delta_{ij}$ at one point, we compute
 \begin{equation}
 \begin{split}
 \Delta|\nabla J|^2=&\nabla_p\nabla_p \left(\nabla_iJ^k_j \nabla_i J^k_j\right)\\
 =&2(\nabla_p\nabla_p\nabla_i J^k_j)  \nabla_i J^k_j+2(\nabla_p\nabla_iJ^k_j)^2.
 \end{split}
 \end{equation}
We compute
 \begin{equation}
 \begin{split}
 \nabla_p\nabla_p\nabla_i J^k_j=&\nabla_i  \nabla_p\nabla_pJ^k_j+\nabla_p \left([\nabla_p, \nabla_i]J^k_j\right)+[\nabla_p, \nabla_i]\left(\nabla_pJ^k_j\right)\\
 =&\nabla_i \Delta J^k_j+\nabla_p(R_{pil}^kJ^l_j-R_{pij}^lJ^k_l)-R_{pip}^l\nabla_lJ^k_j-R_{pij}^l\nabla_pJ^k_l+R_{pil}^k\nabla_pJ^l_j\\
 =&\nabla_i \Delta J^k_j+R^k_{pil, p} J^l_j-R^l_{pij, p}J^k_l+2R^k_{pil}\nabla_pJ^l_j-2R^l_{pij}\nabla_pJ^k_l+R^l_i\nabla_l J^k_j\\
 =&\nabla_i \Delta J^k_j+\cR^k_{ij}
 \end{split}
 \end{equation}
 where $\cR$ denotes the remain terms involved with the curvature. 
 We compute
 \begin{equation}
 \begin{split}
 \cR^k_{ij}\nabla_iJ^k_j=&\left(R^k_{pil, p} J^l_j-R^l_{pij, p}J^k_l+2R^k_{pil}\nabla_pJ^l_j-2R^l_{pij}\nabla_pJ^k_l+R^l_i\nabla_l J^k_j\right)\nabla_iJ^k_j\\
 =&2R^k_{pil, p} J^l_j \nabla_iJ^k_j+4R^k_{pil}\nabla_pJ^l_j \nabla_iJ^k_j+R^l_i\nabla_l J^k_j \nabla_iJ^k_j.
 \end{split}
 \end{equation}
 We compute, using the equation \eqref{H1}, 
 \begin{equation}
 \nabla_i \Delta J^k_j=\nabla_iJ^a_j\nabla_pJ^b_a\nabla_pJ^k_b+J^a_j\nabla_i\nabla_pJ^b_a\nabla_pJ^k_b+J^a_j\nabla_pJ^b_a\nabla_i\nabla_pJ^k_b
 \end{equation}      
 It follows that
\begin{equation}
\begin{split}
\langle \nabla \Delta J, \nabla J\rangle=&\nabla_i \Delta J^k_j\nabla_iJ^k_j\\
=&\nabla_iJ^a_j\nabla_pJ^b_a\nabla_pJ^k_b\nabla_iJ^k_j+J^a_j\nabla_i\nabla_pJ^b_a\nabla_pJ^k_b\nabla_iJ^k_j+J^a_j\nabla_pJ^b_a\nabla_i\nabla_pJ^k_b\nabla_iJ^k_j\\
=&\nabla_iJ^a_j\nabla_pJ^b_a\nabla_pJ^k_b\nabla_iJ^k_j\\
=& -\nabla_iJ^a_j\nabla_iJ^j_k \nabla_pJ^b_a\nabla_pJ^k_b\\
=& -|\nabla_pJ\nabla_pJ|^2
\end{split}
\end{equation} 
where we use the fact that, (at the point, $J^j_i+J^i_j=0$ since $\delta_{ij}=g_{ij}$)
\[\begin{split}J^a_j\nabla_i\nabla_pJ^b_a\nabla_pJ^k_b\nabla_iJ^k_j+J^a_j\nabla_pJ^b_a\nabla_i\nabla_pJ^k_b\nabla_iJ^k_j=&J^a_j\nabla_i\nabla_pJ^b_a\nabla_pJ^k_b\nabla_iJ^k_j+J^k_j\nabla_pJ^b_k\nabla_i\nabla_pJ^a_b\nabla_iJ^a_j\\
=&J^a_j\nabla_i\nabla_pJ^b_a\nabla_pJ^k_b\nabla_iJ^k_j-\nabla_iJ^k_j\nabla_pJ^b_k\nabla_i\nabla_pJ^a_bJ^a_j\\
=&0\end{split}
\]
Combining all the computations above, we have
\begin{equation}
\Delta |\nabla J|^2=2|\nabla^2 J|^2-2|\nabla_pJ\nabla_pJ|^2+4R^k_{pil, p} J^l_j \nabla_iJ^k_j+8R^k_{pil}\nabla_pJ^l_j \nabla_iJ^k_j+2R^l_i\nabla_l J^k_j \nabla_iJ^k_j.
\end{equation}
This completes the proof.
 \end{proof}     

\begin{cor}Denote $u=|\nabla J|^2$, then it satisfies the following,
\begin{equation}\label{bochner}
\Delta u\geq -C(u^2+\sqrt{u}),
\end{equation}
where $C=C(n, |Rm|, |\nabla Rm|)$ depends on the metric. \end{cor}

In applications, we  assume that $|Rm|, |\nabla Rm|$ are both bounded above by a small dimensional constant  by a fixed scaling. Hence the constant $C$ in \eqref{bochner} is assumed to be a dimensional constant.

\subsection{Stationary harmonic almost complex structures}
In this section we discuss a possible approach for \emph{stationary} weakly harmonic almost complex structures following \cite{P} (for Yang-Mills equation and harmonic maps). We would see that there are nontrivial difficulties to carry the theory for stationary harmonic map to our setting. 

Consider $\phi_t: M\rightarrow M$ a  parameter family of diffeomorphisms and we consider $J_t:=(\phi^*_t) J=(\phi^{-1}_t)_{*} \circ J \circ (\phi_t)_{*}$. Suppose $J_0=J\in W^{1, 2}(\cJ_g)$ is a weakly harmonic almost complex structure, satisfying \eqref{H00}. Note that $J_t$ has to be compatible with $g$, 
\begin{equation}\label{R1}
g(J_t, J_t)=g(\phi_t^* J, \phi_t^* J)=g
\end{equation}
We consider the variation of the energy functional
\[
E(J_t)=\int_M |\nabla J_t|^2 dv. 
\]
The key point in \cite{P} is to derive a variational formula not directly on $E(J_t)$, but on the pair of $(g_t=(\phi^{-1}_t)^* g, J)$ which is obviously diffeomorphic to $(g, J_t)$ via $\phi_t$. Hence we have
\begin{equation}\label{stationary1}
E_t(J)=\int_M |\nabla_t J|^2 dv_t=E(J_t). 
\end{equation}
The main advantage to consider $E_t(J)$ over $E(J_t)$ is that when taking variation, the derivatives will be taken on $g_t$ hence no second derivatives on $J$ would appear. 
This would be crucial for $J\in W^{1, 2}(\cJ_g)$.
Denote  $g_t=\psi_t^* g$, where $\psi_t=\phi_t^{-1}$ is generated by a vector field $X$. Then it is straightforward to see that the restriction \eqref{R1} is equivalent to the following condition on $X$,
\begin{equation}\label{R2}
L_Xg(J \cdot , J \cdot )-L_Xg(\cdot, \cdot)=0=g(L_XJ \cdot, \cdot)+g(\cdot, L_XJ\cdot). 
\end{equation}

Similar to Price's approach, it is then attempting to define the following, 
\begin{defn}Let $J\in W^{1, 2}(\cJ_g)$ be a weakly harmonic almost complex structure. We call $J$ is a stationary (weakly) harmonic almost complex structure if it is a stationary critical point of the energy function $E_t(J)$ in \eqref{stationary1}, for any $g_t=\psi_t^* g$ with $\psi_0=id$. 
\end{defn}

\begin{prop}Let $X$ be a $C^2$ vector field satisfying \eqref{R2} and let $J$ be a stationary harmonic almost complex structure. Then we have
\begin{equation}\label{H6}
\int_M |\nabla J|^2 div(X)-2(\nabla_iJ, \nabla_j J) g^{ki}\nabla_kX^j +4 R^l_{j\beta p}J^p_k\nabla_\beta J^l_k X^j=0.
\end{equation}
\end{prop}

\begin{proof}Denote $g_t=\psi^*_t g$, with $\psi_t$ generated by $X$. 
We need to compute, at $t=0$,
\[
\p_t |\nabla_t J|^2=2\left(\p_t \nabla_t J, \nabla J\right)-2(\nabla_i J, \nabla_j J) \left(g^{ik}X^j_{,k}\right)
\]
The direct computation then gives 
\[
\p_t E_t(J)=\int_M |\nabla J|^2 div(X)-2(\nabla_iJ, \nabla_j J) g^{ki}\nabla_kX^j +2(\p_t (\nabla_t J), \nabla J)
\]
Hence we only need to show 
\begin{equation}\label{em0}
A:=\int_M  2(\p_t (\nabla_t J), \nabla J)=\int_M 4 R^l_{j\beta p}J^p_k\nabla_\beta J^l_k X^j 
\end{equation}
We compute
\[
\begin{split}
2\p_t (\nabla_t J)=&-2\p_t \Gamma^l_{ij}J^k_l+2\p_t \Gamma^k_{il}J^l_j
\end{split}
\]
Then we have
\begin{equation}\label{em1}
\begin{split}
A=&\int_M \left(-2\p_t \Gamma^l_{ij}J^k_l+2\p_t \Gamma^k_{il}J^l_j, \nabla_iJ^k_j\right)\\
=&
-4\int_M \left(\p_t \Gamma^l_{ij}J^k_l, \nabla_iJ^k_j\right)\\
\end{split}
\end{equation}
To proceed, we write at $t=0$, 
\[
h_{ij}:=\p_t g_{ij}=(L_X g)_{ij}=g_{kj}{X^k}_{,i} +g_{ki}{X^k}_{,j} 
\]
We compute, at $t=0$, 
\begin{equation}\label{em2}
\begin{split}
2\p_t \Gamma^k_{ij}(t)=&g^{kl}\left(\nabla_i h_{jl}+\nabla_j h_{il}-\nabla_l h_{ij}\right)\\
=&\left(\nabla_i\nabla_j+\nabla_j\nabla_i\right)X^k+g^{kl}\left(g_{aj}[\nabla_i, \nabla_l]+g_{ai}[\nabla_j, \nabla_l]\right)X^a\\
=&\left(\nabla_i\nabla_j+\nabla_j\nabla_i\right)X^k+g^{kl}\left(g_{aj}R^a_{ilb}+g_{ai}R^a_{jlb}\right)X^b\\
=&\left(\nabla_i\nabla_j+\nabla_j\nabla_i\right)X^k-\left(R^k_{jbi}+R^k_{ibj}\right)X^b
\end{split}
\end{equation}
By \eqref{em1} and \eqref{em2}, we have
\begin{equation}\label{em3}
\begin{split}
A=&2\int_M \left(\left(\left(R^l_{jbi}+R^l_{ibj}\right)X^b
 -X^l_{, ji}-X^l_{, ij}\right)J^k_l, \nabla_iJ^k_j\right)\end{split}
\end{equation}
We also have
\begin{equation}\label{em4}
\nabla_j\nabla_iX^l=\nabla_i\nabla_j X^l-R_{ijb}^l X^b,
\end{equation}
Hence we compute, by \eqref{em3} and \eqref{em4}, 
\begin{equation}\label{em5}
\begin{split}
A=&-4\int_M (\nabla_i\nabla_j X^l J^k_l, \nabla_iJ^k_j)+2\int_M \left(R^l_{jbi}+R^l_{ibj}+R_{ijb}^l\right)X^b J^k_l \nabla_iJ^k_j\\
=&-4\int_M  (\nabla_i\nabla_j X^l J^k_l, \nabla_iJ^k_j)+4\int_M \left(R^l_{ibj}X^bJ^k_l, \nabla_iJ^k_j\right)
\end{split}
\end{equation}
Using the fact that $J$ is a weakly harmonic almost complex structure, satisfying \eqref{H00}, we get 
\[
\int_M  (\nabla_i\nabla_j X^l J^k_l, \nabla_iJ^k_j)=0
\]
This together with \eqref{em5} give \eqref{em0}, hence it completes the proof. 
\end{proof}

If there is no restriction condition \eqref{R2},  \eqref{H6} would be enough to derive a monotonicity formula (with a suitable choice of $X$), see \cite{P} for details. However, \eqref{H6} is clearly necessary, due to the $g$-compatibility. It seems to be possible to explore the restriction equation \eqref{R2} locally, by choosing the vector field $X$ in a careful way.  But we are not able to derive a monotonicity formula with \eqref{R2}. We would like to ask, 

\begin{pl}Does a stationary harmonic almost complex structure satisfy a monotonicity formula?
\end{pl}

\subsection{Continuous weakly harmonic almost complex structure  is smooth}
It is well-known that a continuous weakly harmonic map is  smooth. Our argument follows closely  Borchers-Garber \cite{BG}.  However we believe there is a technical overlook in the proof and the statement of their Lemma 2.2; see Lemma \ref{L-1} below and the remark afterwards.  Since the proof in \cite{BG} does not apply directly to our system (even though the proof is very similar), we include a detailed argument for completeness. After we finish our proof, we find that Chang-Yang-Wang \cite{CWY} have given a nice proof that a continuous weakly harmonic map is smooth. Their proof works for a general elliptic system of the form
\[
\Delta u=f(x, u, \nabla u)
\]
such that $f\leq C|\nabla u|^2$. With a minor modification, their proof would cover the system for harmonic almost complex structures. Nevertheless  our arguments should have its own interest, we present them for completeness.

We focus on an Euclidean ball $B=B(0, 1)\subset \R^{2n}$ with a smooth almost Hermitian metric $g=(g_{ij})$. Let $J=(J_j^k)$ be an almost complex structure in $B$ with $L^\infty$ coefficients, satisfying the constraints, for $x\in B$ a. e.
\begin{equation}\label{C1}
J_j^k(x) J_k^l(x)=-\delta_{jl},  g_{kl}(x)J^k_i(x)J_j^l(x)=g_{ij}(x).
\end{equation}
We also require that $\p_iJ_j^k$ is in $L^2_{loc}(B)$. For such \emph{admissible} almost complex structure, we use the notation $J\in L^\infty(M_{2n})\cap W^{1, 2}_{loc}(M_{2n})$.  Note that the constraints for compatible almost complex structure demand that $J\in L^\infty(M_{2n})$.

For the system \eqref{H1}, the essential difficulty comes from nonlinearity in $J\nabla_p J\nabla_p J$, which reads $g^{pq}J_j^k\nabla_p J_k^l \nabla_q J_l^m$. In many situations we need to consider the system with additional lower order terms. Hence we consider the system
\begin{equation}\label{H3}
\p_i\left(g^{ij}\p_j J_k^l\right)=g^{pq}J_k^a\p_p J_a^b \p_q J_b^l+(A_0+A_1*J)*\p J+B_0*J*J*J,
\end{equation}
where $A_0, A_1, B$ involve first derivatives of the metric coefficients, which are smooth and uniformly bounded and the $*$ denotes the possible contractions (finite many summations over indices) and we do not need to understand the precise expression. We shall also use the notation $A*\p J=(A_0+A_1*J)*\p J$ and $B=B_0*J*J*J$
The  weak solution of the system \eqref{H3} is then defined as,  
\begin{defn}
For $J\in L^\infty(M_{2n})\cap W^{1, 2}_{loc}(M_{2n})$ satisfying the constraints \eqref{C1}, we call $J$  a weak solution of the system \eqref{H3} if for any matrix $T\in L^\infty(M_{2n})\cap W^{1, 2}_{loc}(M_{2n})$ and  $\eta\in L^\infty(B)\cap W^{1, 2}_0(B),$ we have
\begin{equation}\label{E1}
\int_B g^{ij}\p_j J_k^l \p_i (\eta^2 T_k^l) = \int_B \eta^2 T_k^l g^{pq}J_k^a\p_p J_a^b \p_q J_b^l  +\int_B \eta^2  T (A*\p J+B).
\end{equation}
where the integration is taken with respect to the Euclidean measure, if not specified otherwise. 
\end{defn}

\begin{rmk}The definition of weak solution above applies even if $J$ does not satisfy the constraints \eqref{C1}. We allow the cut-off function $\eta\in L^\infty\cap W^{1, 2}_0(B)$. The definition of weak solution can allow that the testing matrices are taken in  $L^\infty(M_{2n})\cap W^{1, 2}_0(M_{2n})$. Clearly when $\eta\in L^\infty\cap W^{1, 2}_0(B)$ and $T\in L^\infty(M_{2n})\cap W^{1, 2}_0(M_{2n})$, then $\eta^2 T\in L^\infty(M_{2n})\cap W^{1, 2}_0(M_{2n})$.
\end{rmk}

\begin{thm}\label{main1}Let $J$ be a weak solution of the system \eqref{H3}, and if in addition $J$ is continuous, then $J$ is smooth. 
\end{thm}
The result will be proved by several steps; the main point is to prove that $\p J\in L^p_{loc}$ for any $p>1$. First note that the argument is purely local. We only need to prove for example, in a small neighborhood of $0$, $J$ is smooth. To make the dependence more explicit in the estimates below, we require
\[
\frac{1}{2}\delta_{ij}\leq g_{ij}\leq 2\delta_{ij}, |A|, |B|\leq \Lambda
\]

\begin{lemma}\label{L-1}For any $\epsilon>0$, there exists $\delta=\delta(\epsilon, \Lambda)$ such that if $|J(x)-J(y)|<\delta$ holds  for any  $x, y\in B_r=B(0, r)$, then  for any $\zeta\in L^\infty(B_r)\cap W^{1, 2}_0(B_r)$, 
\[
\int_{B_r}|\p J|^2 \zeta^2\leq \epsilon \left(\int_{B_r} |\p \zeta|^2+\int_{B_r} \zeta^2\right)
\]
We should emphasize $r>0$ is any positive number. In particular, when $J$ is continuous at $x=0$, then for any $\delta>0$, there exists $r$ such that when $x, y\in B_r$, $|J(x)-J(y)|<\delta$. 
\end{lemma}

\begin{proof}Choose the matrix $T(x)=J(x)-J(0)$ in \eqref{E1} and cutoff function $\zeta$, we get
\[
\int_B g^{ij}\p_j J_k^l \p_i (\zeta^2 T_k^l) = \int_B \zeta^2 T_k^l g^{pq}J_k^a\p_p J_a^b \p_q J_b^l  +\int_B \zeta^2  T (A*\p J+B).
\]
It follows that
\[
\int_Bg^{ij}\p_iJ_k^l\p_jJ_k^l \zeta^2-\int_Bg^{ij}\p_j J_k^l \p_i (\zeta^2) T_k^l  =\int_B \zeta^2 T_k^l g^{pq}J_k^a\p_p J_a^b \p_q J_b^l  +\int_B \zeta^2  T (A*\p J+B)
\]
Now assume $|T|\leq \delta$ for $x\in B_r$, then we get
\[
\left(\frac{1}{2}-\delta\right) \int_{B_r}|\p J|^2 \zeta^2\leq 4\delta \int_{B_r}|\p J||\zeta| |\p \zeta|+\delta \Lambda \int_{B_r}\zeta^2 (|\p J|+1)
\]
We estimate
\[
2\int_{B_r}|\p J||\eta| |\p \zeta|\leq \int_{B_r}|\p J|^2\zeta^2+\int_{B_r}|\p\zeta|^2, 2\int_{B_r}\zeta^2 |\p J|\leq \int_{B_r}|\p J|^2\zeta^2+\int_{B_r}\zeta^2
\]
Then it follows that
\[
\left(\frac{1}{2}-3\delta-\Lambda \delta\right) \int_{B_r}|\p J|^2 \zeta^2\leq 2\delta \int_{B_r}|\p\zeta|^2+2\delta\Lambda \int_{B_r}\zeta^2
\]
It completes the proof by choosing $\delta=\delta(\epsilon, \Lambda)$ sufficiently small. 
\end{proof}

\begin{rmk}Our argument clearly follows Lemma 2.2 in \cite{BG}. We believe that there is an overlook in their statement. One should require that $\xi\in L^\infty\cap W^{1, 2}_0$ (or $\xi \in C^0(B)\cap W^{1, 2}_0$) in Lemma 2.2, \cite{BG}. Otherwise the testing function $\varphi=(n(x)-n(y))\xi^2$ in the argument will not be in $W^{1, 2}_0$. Indeed the integrand on the left hand in the inequality in Lemma 2.2 simply does not make sense for general functions $\xi \in W^{1, 2}_0$ (since $|\nabla n|^2$ is merely in $L^1$). This overlook will cause some technical issues in the proof of Theorem 2.4 in \cite{BG}. On the other hand, Lemma 2.2 \cite{BG} is indeed sufficient for the regularity,  with a careful modification. \end{rmk}

We will use a lot difference quotients in the proofs. We go over several important properties. We denote the difference quotient by \[\delta_\alpha^h S=\frac{S(x+he_\alpha)-S(x)}{h},\] where $S$ can be taken as a function or a matrix-valued function. The derivatives and difference quotients are related by
\begin{lemma} Suppose $f\in L^p(\Omega)$.
\begin{enumerate}\item If $\p_\alpha f\in L^p(\Omega)$, then for all $h$ and $\Omega^{'}$ satisfying $|h|<d(\p\Omega, \Omega^{'})$), we have \[\|\delta^h_\alpha f\|_{L^p(\Omega^{'})}\leq \|f\|_{L^p(\Omega)}. \]

\item Suppose there exists a constant $K$ such that $\delta^h_\alpha f\in L^p(\Omega^{'})$ for $p\in (1, \infty)$ and 
\[
\|\delta^h_\alpha f\|_{L^p(\Omega^{'})}\leq K
\]
for all $h$ and $\Omega^{'}\subset \subset \Omega$ such that $|h|<d(\p\Omega, \Omega^{'})$. Then weak derivative $\p_\alpha f$ exists and satisfies 
$\|\p_\alpha f\|_{L^p(\Omega)}\leq K$. Moreover for any $\Omega^{'}\subset \subset \Omega$, $\delta^h_\alpha f$ converges strongly in $L^p(\Omega^{'})$ to $\p_\alpha f$ for $h\rightarrow 0$. 
\end{enumerate}
\end{lemma} 
\begin{proof}These facts are well-known, one can find the proofs  in \cite{GT}, Section 7.11 except the strong convergence part. To prove the strong convergence, see Lemma 2.1 \cite{BG}. 
\end{proof}

We also need the ``integration by parts" and the ``product rule". For any two functions $u, v$, suppose one of them has compact support in $B$, then
\begin{equation}
\int_{B} u\left(\delta^{-h}_\alpha v \right)=\int_B \left(\delta^h_\alpha u\right) v
\end{equation}
The product rule works as follows, 
\begin{equation}
\begin{split}
\delta_\alpha^h\left(u S\right)=&\frac{u(x+he_\alpha)S(x+he_\alpha)-u(x)S(x)}{h}\\
=&\delta^h_\alpha(u) S(x+he_\alpha)+u\delta^h_\alpha (S)\\
=&\delta^h_\alpha(u) S+u(x+he_\alpha)\delta^h_\alpha(S)
\end{split}
\end{equation}
Note that the product rule in the difference quotients is not as clean as the normal derivatives and it does cause some technical complications; but as in Proposition \ref{L-2}, it does not cause essential difficulties. 
Now we are ready to prove
\begin{prop}\label{L-2}$|\p J|\in L^4_{loc}(B)$ and $|\p^2 J|\in L^2_{loc}(B)$.
\end{prop}
\begin{proof}We only need to show the statement holds near zero. 
Consider the following matrix $\delta^{-h}_\alpha\left(\eta^2 \delta^h_\alpha J\right)$ for $\eta\in C^1_0(B)$, where $|h|$ is small enough such that $\delta^{\pm h}_\alpha \eta\in C^1_0(B)$. Then we have
\[
\int_B g^{ij}\p_iJ^l_k\p_j\left(\delta^{-h}_\alpha\left(\eta^2 \delta^h_\alpha J\right)\right)=\int_B \delta^{-h}_\alpha\left(\eta^2 \delta^h_\alpha J\right)\left(g^{pq}J_k^a\p_p J_a^b \p_q J_b^l  +(A*\p J+B)\right).
\]
We will use the notations $\delta J=(\delta^h_\alpha J)$, $\alpha=1, \cdots, 2n$ and $|\delta J|^2=\sum_{k, l, \alpha}(\delta^h_\alpha J^l_k)^2$. 
We compute
\[
\begin{split}
\int_B g^{ij}\p_iJ^l_k\p_j\left(\delta^{-h}_\alpha\left(\eta^2 \delta^h_\alpha J\right)\right)=&\int_B\delta^h_\alpha\left(g^{ij}\p_iJ^l_k\right)\p_j\left(\eta^2 \delta^h_\alpha J\right)\\
=&\int_B\left(g^{ij}(x+he_\alpha)\p_i(\delta^h_\alpha J^l_k)+\delta^h_\alpha(g^{ij})\p_iJ^l_k(x+he_\alpha)\right)\p_j\left(\eta^2 \delta^h_\alpha J^l_k\right)\\
=&\int_B(g^{ij}(x+he_\alpha)\p_i(\delta^h_\alpha J^l_k) \p_j\left(\delta^h_\alpha J^l_k\right)\eta^2\\&+\int_B g^{ij}(x+he_\alpha)\p_i(\delta^h_\alpha J^l_k)\p_j(\eta^2) \delta^h_\alpha J^l_k\\
&+\int_B \delta^h_\alpha(g^{ij})\p_iJ^l_k(x+he_\alpha)\left( \p_j(\delta^h_\alpha J) \eta^2+2\p_j(\eta)\eta (\delta^h_\alpha J) \right)\\
\geq&\frac{1}{2}\int_B |\p\delta J|^2\eta^2-C\int_B |\p \delta J||\delta J| \eta|\p\eta|\\&-C\int_B |\p J|(x+he_\alpha) \left(\eta^2|\p\delta J|+\eta|\p \eta||\delta J|\right)\\
\geq &\frac{1}{3}\int_B |\p\delta J|^2\eta^2-C\int_B |\delta J|^2(\eta^2+|\p \eta|^2)-C\int_B |\p J|^2 \eta^2_\alpha,
\end{split}
\]
where we estimate in particular
\[
C\int_B |\p J|(x+he_\alpha) \left(\eta^2|\p\delta J|+\eta|\p \eta||\delta J|\right)\leq C\int_B \eta^2|\p J|^2(x+he_\alpha)+\frac{1}{10}\int_B |\p\delta J|^2\eta^2+C\int_B|\delta J|^2|\p\eta|^2
\]
and we can write the integral, with $\eta_\alpha(x)=\eta(x-he_\alpha)$, 
\[
\int_B \eta^2|\p J|^2(x+he_\alpha)=\int_B \eta^2(x-he_\alpha)|\p J|^2(x)=\int_B \eta^2_\alpha |\p J|^2. 
\]
The constant $C$ is uniform but can vary line by line.
Next we need to deal with
\[
\int_B \delta^{-h}_\alpha\left(\eta^2 \delta^h_\alpha J\right)\left(g^{pq}J_k^a\p_p J_a^b \p_q J_b^l  +(A*\p J+B)\right)=\int_B\eta^2 \delta^h_\alpha J\delta^h_\alpha \left(g^{pq}J_k^a\p_p J_a^b \p_q J_b^l  +(A*\p J+B)\right)
\]
We need be very careful about the product rule of difference quotients, as noted above. 
Note that $A=A_0+A_1*J$, $B=B_0*J*J*J$, where $A_0, A_1, B_0$ involve with first derivatives of the metric. The difference quotients of these terms are uniformly bounded.  We will use a uniform constant $C$ to denote these bounds, which can again vary line by line. We compute
\[
\begin{split}
\int_B\eta^2 \delta^h_\alpha J\delta^h_\alpha \left(g^{pq}J_k^a\p_p J_a^b \p_q J_b^l \right) \leq&C \int_B \eta^2|\delta J|\left(|\p J|^2(x+he_\alpha)+|\delta J||\p J|^2(x+he_\alpha)\right)\\
&+C\int_B \eta^2 |\delta J||\p J||\p \delta J|+\eta^2 |\delta J||\p \delta J| |\p J|(x+he_\alpha).
\end{split}
\]
We need to estimate the terms on the right hand side. First we have,
\[
\begin{split}
 \int_B \eta^2|\delta J|\left(|\p J|^2(x+he_\alpha)+|\delta J||\p J|^2(x+he_\alpha)\right)\leq 2\int_B \eta^2 (|\delta J|^2+1)|\p J|^2(x+he_\alpha)\\
\end{split}
\]
We also estimate 
\[
C\int_B \eta^2 |\delta J||\p J||\p \delta J|\leq \frac{1}{100}\int_B |\p\delta J|^2\eta^2+C\int_B\eta^2|\delta J|^2|\p J|^2
\]
and
\[
C\int_B \eta^2 |\delta J||\p \delta J| |\p J|(x+he_\alpha)\leq \frac{1}{100}\int_B |\p\delta J|^2\eta^2+C\int_B \eta^2 |\delta J|^2 |\p J|^2(x+he_\alpha)
\]
The lower order terms involved with $(A_0+A_1*J)*\p J+B_0 J*J*J$ are actually easier to estimate. Similarly we have
\[
\int_B\eta^2 \delta^h_\alpha J\delta^h_\alpha (A*\p J+B)\leq \frac{1}{100}\int_B \eta^2 |\p \delta J|^2+C\int_B \eta^2|\delta J|^2+\int_B\eta^2 |\delta J|^2 |\p J|^2 (x+he_\alpha)
\]
Combine the estimates above together, we reach the following
\begin{equation}\label{E2}
\begin{split}
\int_B |\p\delta J|^2\eta^2\leq & C\int_B |\delta J|^2 |\p J|^2 \eta^2+C\int_B (\eta^2+|\p \eta|^2)(|\delta J|^2+|\p J|^2)\\
&+C\int_B \eta^2(|\delta J|^2+1) |\p J|^2(x+he_\alpha)
\end{split}
\end{equation}
Note that for the last term, we can also write it in the form
\[
\begin{split}
\int_B \eta^2(|\delta J|^2+1) |\p J|^2(x+he_\alpha)=&\int_B \eta^2(x)(|\delta^h_\alpha J|^2(x)+1)|\p J|^2(x+he_\alpha)\\
=&\int_B\eta^2(y-he_\alpha) (|\delta^{-h}_\alpha J|^2(y)+1) |\p J|^2(y)\\
=&\int_B \eta^2_\alpha (|\delta^{-h} J|^2+1)|\p J|^2
\end{split}
\]
We emphasize that for \eqref{E2}, we have not used the fact that $J$ is continuous and it holds for all weak solutions. We need to use Lemma \ref{L-1} in the following, which depends on the continuity of $J$ in an essential way. For a cutoff function $\eta\in C^1_0(B)$, first we take $\zeta=\eta |\delta J|$. When $h$ is small enough, $\zeta=\eta|\delta J|\in C^0(B_r)\cap W^{1, 2}_0(B_r)$. By Lemma \ref{L-1}, 
\[
\int_{B_r}|\p J|^2 |\delta J|^2 \eta^2\leq \epsilon\int_{B_r} \left(|\p (\eta|\delta J|)|^2+|\delta J|^2\eta^2\right)\leq 4\epsilon \int_{B_r} |\delta J|^2 (\eta^2+|\p \eta|^2)+4\epsilon\int_{B_r}|\p \delta J|^2\eta^2,
\]
where we have use the fact that $|\p |\delta J||\leq |\p \delta J|$. Next we need to take $\zeta=\eta_\alpha |\delta^{-h}| J$, then we get
\begin{equation}\label{T1}
\begin{split}
\int_B \eta^2|\delta J|^2 |\p J|^2(x+he_\alpha)=& \int_B \eta^2_\alpha |\delta^{-h} J|^2|\p J|^2\\
\leq & \epsilon \int_B |\p (\eta_\alpha \delta^{-h} J)|^2+\epsilon\int_B \eta_\alpha^2|\delta^{-h} J|^2\\
=& \epsilon \int_B |\p (\eta \delta J)|^2+\epsilon \int_B \eta^2 |\delta J|^2 ,
\end{split}
\end{equation}
where the last step follows by the invariance of translation of integral, provided that $\eta_\alpha$ and $\eta$ are both supported in $B$.

By taking $\epsilon$ small, depending on the constant $C$ in \eqref{E2}, we get the following two inequalities, 
\begin{equation}\label{E-3}
\int_{B_r} |\p\delta J|^2\eta^2\leq C\int_{B_r} (\eta^2+|\p \eta|^2)(|\delta J|^2+|\p J|^2)
\end{equation}
and 
\begin{equation}\label{E-4}
\int_{B_r}|\p J|^2 |\delta J|^2 \eta^2\leq C\int_{B_r} (\eta^2+|\p \eta|^2)(|\delta J|^2+|\p J|^2)
\end{equation}
With \eqref{E-3} and \eqref{E-4}, it follows that $|\p J|\in L^4_{loc}$ and $|\p^2 J|\in L^2_{loc}$ near zero. 
The argument clearly applies to any interior point in $B$, if we assume that $J$ is continuous. It completes the proof. 
\end{proof}

The proof of Proposition \ref{L-2} essentially follows the arguments in Theorem 2.3 in \cite{BG}. There are two technical differences. The first one  is that our cutoff functions are all in $L^\infty\cap W^{1, 2}_0$. We believe this carefulness is necessary. The second one is instead of using translation invariance of the equation as in Theorem 2.3 in \cite{BG}, we use \eqref{T1}, based on invariance of integrals under translation,  to handle the terms which appear as the product rule of difference quotients, such as $|\delta J|^2 |\p J|^2(x+he_\alpha)$.  Our arguments can be easily modified to extend to show that $\p J\in L^p$ for any $4\leq p<\infty$. We do not use the approach of Theorem 2.4 in 
\cite{BG} and we believe again one needs to be careful of the choice of the cutoff functions in Lemma 2.2 \cite{BG} or Lemma \ref{L-1}.  

\begin{thm}\label{main2}For any $p\geq 1$, $|\p J|\in L^p_{loc}(B)$. 
\end{thm}
\begin{proof}We only need to prove this statement near zero. 
The argument is similar to the argument in Proposition \ref{L-2}, and an induction process.  
Consider the following testing matrix $\delta^{-h}_\alpha\left(\eta^2|\delta J|^{2s}\delta^h_\alpha J\right)$, where $s\in \N$. When $s=0$, this coincides the testing matrix in Proposition \ref{L-2}. We have the following, by \eqref{E1} and ``integration by parts", 
\begin{equation*}
\int_B \delta^h_\alpha\left(g^{ij}\p_i J^l_k\right)\p_j\left(|\delta J|^{2s}\eta^2 \delta^h_\alpha J^l_k\right)=\int_B |\delta J|^{2s}\eta^2(\delta^h_\alpha J^l_k)\delta^h_\alpha\left(g^{pq}J_k^q\p_pJ^b_a\p_qJ^l_b+A*\p J+B\right)
\end{equation*}
We need to do the estimates as above, taking account of the product rule for difference quotients carefully. First we have, 
\[
\begin{split}
\int_B \delta^h_\alpha\left(g^{ij}\p_i J^l_k\right)\p_j\left(|\delta J|^{2s}\eta^2 \delta^h_\alpha J^l_k\right)=&\int_B g^{ij}(x+he_\alpha)\p_i(\delta^h_\alpha J^l_k)\p_j\left(|\delta J|^{2s}\eta^2 \delta^h_\alpha J^l_k\right)\\
&+\int_B \delta^h_\alpha (g^{ij}) \p_iJ^l_k\p_i(\delta^h_\alpha J^l_k)\p_j\left(|\delta J|^{2s}\eta^2 \delta^h_\alpha J^l_k\right)
\end{split}
\]
We estimate
\[
\begin{split}
\int_B g^{ij}(x+he_\alpha)\p_i(\delta^h_\alpha J^l_k)\p_j\left(|\delta J|^{2s}\eta^2 \delta^h_\alpha J^l_k\right)\geq& \frac{1}{2}\int_B |\p\delta J|^2|\delta J|^{2s}\eta^2-C\int_B |\p\delta J||\delta J|^{2s+1}|\p \eta|\eta\\
&+s\int_B g^{ij}_\alpha\p_i(\delta^h_\alpha J^l_k)(\delta^h_\alpha J^l_k) |\delta J|^{2s-2}\eta^2\p_j(\delta^h_\beta J^b_a)(\delta^h_\beta J^b_a),
\end{split}
\]
where we use the notation $g^{ij}_\alpha(x)=g^{ij}(x+he_\alpha)$. We have, 
\[
\int_B g^{ij}_\alpha\p_i(\delta^h_\alpha J^l_k)(\delta^h_\alpha J^l_k) |\delta J|^{2s-2}\eta^2\p_j(\delta^h_\beta J^b_a)(\delta^h_\beta J^b_a)\geq \frac{1}{2}\int_B |\p\delta J*\delta J|^2|\delta J|^{2s-2}\eta^2\geq 0.
\]
We need to estimate, using the product rule of difference quotients, 
\[
\begin{split}
\int_B |\delta J|^{2s}\eta^2(\delta^h_\alpha J^l_k)\delta^h_\alpha\left(g^{pq}J_k^q\p_pJ^b_a\p_qJ^l_b\right)\leq &C\int_B|\delta J|^{2s+1}\eta^2|\p\delta J|\left(|\p J|+|\p J|_\alpha\right)\\
&+C\int_B |\delta J|^{2s+2}|\p J|^2\eta^2+\int_B |\delta J|^{2s+1}\eta^2 |\p J|^2,
\end{split}
\]
where $|\p J|_\alpha=|\p J|(x+he_\alpha)$. We also estimate
\[
\int_B |\delta J|^{2s}\eta^2(\delta^h_\alpha J^l_k)\delta^h_\alpha\left(A*\p J+B\right)\leq C\int_B |\delta J|^{2s+1}\eta^2\left(|\p \delta J|+|\p J||\delta J|+|\p J|\right)
\]
We can estimate the following,
\[
C\int_B|\delta J|^{2s+1}\eta^2|\p\delta J||\p J|\leq \frac{1}{100}\int_B|\p\delta J|^2 |\delta J|^{2s}\eta^2+C\int_B |\delta J|^{2s+2}|\p J|^2\eta^2
\]
and
\[
C\int_B|\delta J|^{2s+1}\eta^2|\p\delta J||\p J|_\alpha\leq \frac{1}{100}\int_B|\p\delta J|^2 |\delta J|^{2s}\eta^2+C\int_B |\delta J|^{2s+2}|\p J|^2_\alpha\eta^2
\]
Put all above together, we reach the estimate
\begin{equation}\label{E5}
\begin{split}
\int_B |\p\delta J|^2|\delta J|^{2s}\eta^2\leq &C\int_B |\delta J|^{2s+2}\eta^2(|\p J|^2+|\p J|_\alpha^2)+C\int_B|\delta J|^{2s+2}(|\p \eta|^2+\eta^2)\\
&+C\int_B|\delta J|^{2s}|\p J|^2 \eta^2
\end{split}
\end{equation}
Note that when $s=0$, this coincides the estimates above. Next we need to estimate, by taking $\zeta=\eta|\delta J|^{s+1}$, 
\[
\begin{split}
\int_B |\delta J|^{2s+2}\eta^2|\p J|^2\leq \epsilon \int_B |\p (\eta |\delta J|^{s+1})|^2\leq 2\epsilon (s+1)\int_B |\delta J|^{2s}|\p\delta J|^2+C\epsilon \int_B|\p \eta|^2|\delta J|^{2s+2}.
\end{split}
\]
Note that we have,
\[
\int_B |\delta J|^{2s+2}\eta^2|\p J|_\alpha^2=\int_B |\delta^{-h}J|^{2s+2}\eta_\alpha^2|\p J|^2
\]
By taking $\zeta=\eta_\alpha |\delta^{-h}J|^{s+1}$, we compute
\[
\int_B |\delta^{-h}J|^{2s+2}\eta_\alpha^2|\p J|^2\leq \epsilon \int_B |\p (\eta_\alpha |\delta^{-h} J|^{s+1})|^2=\epsilon \int_B |\p (\eta|\delta J|^{s+1})|,
\]
which can be estimated in the same way. In other words, by taking $\epsilon=\epsilon (s, C)$ sufficiently small, we have
\begin{equation}\label{E6}
C\int_B |\delta J|^{2s+2}\eta^2(|\p J|^2+|\p J|_\alpha^2)\leq \frac{1}{100}\int_B |\p\delta J|^2|\delta J|^{2s}\eta^2+C_s\int_B |\p \eta|^2|\delta J|^{2s+2},
\end{equation}
where $C_s$ is to emphasize the dependence on $s$ ($C_s\rightarrow \infty$ when $s\rightarrow \infty$). 
By \eqref{E5} and \eqref{E6}, we have the following,
\begin{equation}\label{E7}
\begin{split}
&\int_B |\p\delta J|^2|\delta J|^{2s}\eta^2\leq C_s\int_B|\delta J|^{2s+2}(|\p \eta|^2+\eta^2)+C_s\int_B|\delta J|^{2s}|\p J|^2 \eta^2\\
&\int_B |\delta J|^{2s+2}\eta^2(|\p J|^2+|\p J|_\alpha^2)\leq C_s\int_B|\delta J|^{2s+2}(|\p \eta|^2+\eta^2)+C_s\int_B|\delta J|^{2s}|\p J|^2 \eta^2,
\end{split}
\end{equation}
With \eqref{E7} and Proposition \ref{L-2}, we prove that $|\p J|\in L^p_{loc}$  near zero for any $p$, by induction on $s=1, 2, \cdots$. 
\end{proof}

\begin{rmk}Our proof of $|\p J|\in L^p$ uses essentially the same idea in Theorem 2.3 \cite{BG} and Proposition \ref{L-2} above, with an induction argument. This not only takes care of the technical point of choosing test functions, but also makes the argument  more transparent and streamlined.  Our proof clearly works for other similar elliptic systems. \end{rmk}

Once $\p J\in L^p_{loc}(B)$ for any $p\in (1, \infty)$, then the weak equation implies that $\p_i\left(g^{ij}\p_j J\right)\in L^p$, or $g^{ij}\p^{2}_{ij} J\in L^p$, note that $\p^2_{ij} J$ is well-defined in $L^2$. The $L^p$ theory then implies that $J\in W^{2, p}$, hence in $C^{1, \alpha}$. It then follows from the standard bootstrapping argument that $J$ is smooth. This proves Theorem \ref{main1}.

\end{document}